\documentclass[twoside,11pt]{article}

\usepackage[margin=1in]{geometry}
\usepackage{natbib}
\usepackage{amsthm}
\usepackage{amssymb}
\usepackage{graphicx}
\usepackage{lastpage}
\usepackage{amsmath}
\usepackage{mathtools}
\usepackage{bm}
\usepackage{booktabs}
\usepackage{array}
\usepackage{enumitem}
\usepackage{microtype}
\usepackage{xcolor}
\usepackage{placeins}
\usepackage{capt-of}
\usepackage{hyperref}
\providecommand{\BlackBox}{\rule{1.3ex}{1.3ex}}
\newenvironment{keywords}
  {\par\smallskip\noindent\textbf{Keywords: }}
  {\par\smallskip}
\hypersetup{
  colorlinks=true,
  citecolor=blue,
  linkcolor=black,
  urlcolor=blue
}
\allowdisplaybreaks
\renewenvironment{proof}[1][Proof]
  {\par\noindent\textbf{#1.}\ }
  {\unskip\nobreak\hfill\penalty50\hskip.5em\hbox{}\nobreak\hfill
   \BlackBox\par\vspace{2mm}}

\newtheorem{theorem}{Theorem}
\newtheorem{assumption}{Assumption}
\newtheorem{lemma}{Lemma}
\newtheorem{proposition}{Proposition}
\newtheorem{remark}{Remark}
\newtheorem{corollary}{Corollary}
\newtheorem{definition}{Definition}
\newtheorem{conjecture}{Conjecture}
\newtheorem{axiom}{Axiom}

\newcommand{\R}{\mathbb{R}}
\newcommand{\E}{\mathbb{E}}
\newcommand{\Pp}{\mathbb{P}}
\newcommand{\cS}{\mathcal{S}}
\newcommand{\cA}{\mathcal{A}}
\newcommand{\cC}{\mathcal{C}}
\newcommand{\cX}{\mathcal{X}}
\newcommand{\cF}{\mathcal{F}}
\newcommand{\cN}{\mathcal{N}}
\newcommand{\cP}{\mathcal{P}}
\newcommand{\btheta}{\bm{\theta}}
\newcommand{\bx}{\bm{x}}
\newcommand{\by}{\bm{y}}
\newcommand{\bz}{\bm{z}}
\newcommand{\bu}{\bm{u}}
\newcommand{\bw}{\bm{w}}
\newcommand{\bq}{\bm{q}}
\newcommand{\be}{\bm{e}}
\newcommand{\bH}{\bm{H}}
\newcommand{\bg}{\bm{g}}
\newcommand{\bA}{\bm{A}}
\newcommand{\bB}{\bm{B}}
\newcommand{\bC}{\bm{C}}
\newcommand{\bD}{\bm{D}}
\newcommand{\bI}{\bm{I}}
\newcommand{\bK}{\bm{K}}
\newcommand{\bL}{\bm{L}}
\newcommand{\bP}{\bm{P}}
\newcommand{\bQ}{\bm{Q}}
\newcommand{\bR}{\bm{R}}
\newcommand{\bU}{\bm{U}}
\newcommand{\bmM}{\bm{M}}
\newcommand{\bG}{\bm{G}}
\newcommand{\bGamma}{\bm{\Gamma}}
\newcommand{\bOmega}{\bm{\Omega}}
\newcommand{\bS}{\bm{S}}
\newcommand{\bW}{\bm{W}}
\newcommand{\bcalW}{\bm{\mathcal{W}}}
\newcommand{\bPi}{\bm{\Pi}}
\newcommand{\bSigma}{\bm{\Sigma}}
\newcommand{\bPhi}{\bm{\Phi}}
\newcommand{\bphi}{\bm{\phi}}
\newcommand{\bvarphi}{\bm{\varphi}}
\newcommand{\bepsilon}{\bm{\epsilon}}
\newcommand{\bvarepsilon}{\bm{\varepsilon}}
\newcommand{\bzeta}{\bm{\zeta}}
\newcommand{\bnu}{\bm{\nu}}
\newcommand{\br}{\bm{r}}
\newcommand{\bDelta}{\bm{\Delta}}
\newcommand{\bzero}{\bm{0}}
\newcommand{\norm}[1]{\left\lVert #1\right\rVert}

\newcommand{\diag}{\operatorname{diag}}
\newcolumntype{P}[1]{>{\raggedright\arraybackslash}p{#1}}
\newcommand{\appendixsection}[1]{%
  \section{#1}%
  \addcontentsline{apc}{section}{\protect\numberline{\thesection}#1}}

\makeatletter
\newcommand{\appendixcontentsline}[2]{%
  \begingroup
  \parindent\z@ \rightskip 2em \parfillskip -2em
  \bfseries\leavevmode
  \setlength{\@tempdima}{1.8em}%
  \advance\leftskip\@tempdima
  \hskip-\leftskip
  #1\nobreak\hfill\nobreak\hb@xt@2em{\hfil #2}\par
  \endgroup
  \vspace{0.55em}}
\newcommand{\printappendixcontents}{%
  \begingroup
  \let\l@section\appendixcontentsline
  \@starttoc{apc}%
  \endgroup}
\makeatother

\begin{document}

\title{Online Statistical Inference for Nonlinear Stochastic Approximation with Markovian Data}

\author{Xiang Li\quad Jiadong Liang\quad Zhihua Zhang\\[0.5em]
  \small Department of Probability and Statistics,
  School of Mathematical Sciences, Peking University\\
  \small Beijing 100871, China\\[0.25em]
  \small\texttt{lx10077@pku.edu.cn},
  \texttt{jdliang@pku.edu.cn},
  \texttt{zhzhang@math.pku.edu.cn}}
\date{}

\maketitle

\begin{abstract}%
Many stochastic approximation (SA) algorithms evolve along a single trajectory, making uncertainty quantification challenging under nonlinear dynamics and Markov dependence. We develop an online inference framework for nonlinear SA with decreasing step sizes. Under local stability and verifiable conditions, we establish a functional central limit theorem for the partial-sum path. The proof uses a Poisson-equation decomposition to handle Markov dependence and a uniform bound to control endpoint-dependent remainders induced by decreasing step sizes. The resulting path limit yields self-normalized confidence intervals without estimating the asymptotic variance. Our primary construction uses a five-dimensional polynomial-series normalizer, has an asymptotic Student $t_5$ pivot, and requires only constant memory. We apply the framework to Q-learning, including asynchronous tabular, projected linear, and entropy-regularized updates, as well as SGD for generalized linear models with Markov data and inference for the identified product in low-rank adaptation (LoRA). Across four settings, the polynomial-series method achieves near-nominal coverage, with shorter confidence intervals and lower computational cost than online bootstrap.
\end{abstract}


\section{Introduction}
\label{sec:introduction}

Many modern machine learning and decision-making systems update their parameters sequentially as data arrive. Stochastic approximation (SA), originating with the Robbins--Monro algorithm, provides a general framework: each observation generates a recursive update toward a root of a population equation
\citep{robbins1951stochastic,borkar2009stochastic,kushner2003stochastic}. Its fixed memory and low per-iteration cost suit long or indefinitely generated data streams. Classical SA theory studies convergence to the target parameter. For decision-making, however, one must also quantify uncertainty, which is especially difficult when (i) the recursion is nonlinear and (ii) both learning and inference must rely on a single dependent trajectory.

These two challenges arise naturally in several applications. In Q-learning, an agent learns the long-run value of each state--action pair by repeatedly interacting with an environment. The resulting observations form a dependent trajectory because each new state depends on the current state and action, while the update is nonlinear because it involves maximizing over the actions available at the next state \citep{watkins1989learning}. Markov stochastic gradient descent (SGD) has a similar dependence structure when observations are generated sequentially by a time series, a simulator, or another evolving system \citep{bottou2018optimization,even2023markov}. Low-rank adaptation (LoRA) introduces a different form of nonlinearity: it represents a parameter update as the product of two low-rank factor matrices and updates both factors jointly \citep{hu2022lora}. These examples motivate our study of inference for nonlinear SA driven by a single Markov trajectory.

We use the following nonlinear SA model to capture these features. Let $\bx\in\R^d$ be the parameter, $(\xi_t)_{t\ge0}$ a Markov chain on a measurable space $\cX$ with a transition kernel $\cP$ and an invariant distribution $\mu$, and $\bH\colon\R^d\times\cX\to\R^d$ the update function. For a measurable scalar function $f$, we define $\cP f(\bx,\xi):=\int_{\cX}f(\bx,\xi')\cP(\xi,\mathrm{d}\xi')$, with the same convention for vector- and matrix-valued functions. The target $\bx^\star$ is a root of the mean field $\bg(\bx):=\int_{\cX}\bH(\bx,\xi)\,\mu(\mathrm{d}\xi)$, so $\bg(\bx^\star)=\bzero$.
Given  an initial  $\bx_0$, the SA algorithm follows
\begin{equation}
  \bx_{t+1}=\bx_t-\eta_t\bH(\bx_t,\xi_t),
  \label{eq:sa-recursion}
\end{equation}
where $\eta_t>0$ is a decreasing step size. Let $\cF_{t-1}:=\sigma(\xi_0,\ldots,\xi_{t-1})$ denote the information available before iteration $t$. Because $\bx_t$ is $\cF_{t-1}$-measurable, the Markov property yields $\E\{\bH(\bx_t,\xi_t)\mid\cF_{t-1}\}=\cP\bH(\bx_t,\xi_{t-1})$. This conditional expectation generally differs from $\bg(\bx_t)$, so $\bH(\bx_t,\xi_t)$ is conditionally biased relative to the mean field.

\paragraph{From CLT to FCLT.}
A natural starting point is a central limit theorem (CLT) for the Polyak--Ruppert average $\bar\bx_T:=T^{-1}\sum_{t=1}^T\bx_t$. Fix a nonzero vector $\bm v\in\R^d$ and consider the scalar linear functional $\bm v^\top\bx^\star$; choosing $\bm v$ as a coordinate vector selects an individual component. Under suitable regularity and stability conditions,
\[
\sqrt{T} \bm v^\top(\bar\bx_T-\bx^\star)
\ \Longrightarrow
\sigma_{\bm v} Z,
\qquad Z\sim\mathcal{N}(0,1),
\]
where the unknown scale $\sigma_{\bm v}$ captures both the local sensitivity of the recursion and temporal dependence along the trajectory \citep{polyak1992acceleration,liang2010trajectory}. Because $\sigma_{\bm v}$ is unknown, the CLT alone does not yield a directly computable confidence interval.

Two common online approaches address this unknown scale. Plug-in and batch-means methods estimate $\sigma_{\bm v}^2$ from the iterate path, but may converge slowly under persistent Markov dependence \citep{chen2020statistical,zhu2021online,roy2023online}. Multiplier bootstrap instead runs randomly perturbed copies of the same recursion, avoiding explicit variance estimation at the cost of additional parameter trajectories and evaluations of $\bH$. 
Existing validity theory for multiplier methods under Markov noise, however, has primarily been developed for linear SA, and its extension to nonlinear SA remains unclear
\citep{ramprasad2021online}.

Random scaling offers a third approach. It avoids both variance estimation and auxiliary recursions by constructing, from the same trajectory, a normalizer whose limit contains the same unknown scale $\sigma_{\bm v}$ as the numerator
\citep{lee2021fast}. The resulting ratio therefore has a limiting distribution that does not depend on $\sigma_{\bm v}$. Justifying this construction requires more than a CLT for $\bar\bx_T$, because the normalizer depends on fluctuations accumulated over intermediate portions of the trajectory. We therefore study the linearly interpolated partial-sum process
\begin{equation}
	\bPhi_T(r)=T^{-1/2}\!\left\{
	\sum_{t=1}^{\lfloor Tr\rfloor}(\bx_t-\bx^\star)
	+(Tr-\lfloor Tr\rfloor)(\bx_{\lfloor Tr\rfloor+1}-\bx^\star)\right\},
	\qquad r\in[0,1].
	\label{eq:partial-sum-process}
\end{equation}
At $r=1$, $\bm v^\top\bPhi_T(1)=\sqrt{T}\,\bm v^\top(\bar\bx_T-\bx^\star)$, the quantity appearing in the CLT. For $r<1$, the process records how the estimation error accumulates along the trajectory. An FCLT establishes that $\bPhi_T(\cdot)$ converges weakly in $C([0,1],\R^d)$, equipped with the uniform norm, to a linear transformation of standard $d$-dimensional Brownian motion.  It therefore describes the final average and all intermediate partial sums jointly.

\paragraph{Establishing the FCLT.}
Our first main result establishes such an FCLT for nonlinear, decreasing-step-size SA
with Markov data. To address the conditional bias identified above, we use the classical Poisson-equation approach
\citep{benveniste2012adaptive,andrieu2005stability}. Specifically, we assume that, for every $\bx\in\R^d$, the Poisson equation
\begin{equation}
  \bU(\bx,\xi)-\cP\bU(\bx,\xi)
  =\bH(\bx,\xi)-\bg(\bx).
  \label{eq:poisson-family}
\end{equation}
admits a solution $\bU(\bx,\cdot)$, normalized to have mean zero under $\mu$.
The Poisson solution accounts for dependence in the observations and decomposes the conditionally biased fluctuation $\bH(\bx_t,\xi_t)-\bg(\bx_t)$ into a martingale term, denoted by $\bepsilon_t$, and several remainder terms. 
Rather than tying the FCLT to a particular convergence theorem for the algorithm, we state it under general conditions that can be checked using application-specific results.

The central technical challenge is to control the remainder terms uniformly. 
The most delicate term is
$T^{-1/2}\sum_{j=1}^{\lfloor Tr\rfloor}
(\bm A_j^{\lfloor Tr\rfloor}-\bm A_j^T)\bepsilon_j$. 
Here $\bepsilon_j$ is the above Poisson martingale term, and $\bm A_j^n$ measures its total effect on the partial sum ending at $n$. Because these coefficients depend on the endpoint $\lfloor Tr\rfloor$, the resulting process is not a martingale, so standard martingale maximal inequalities do not apply directly. Lemma~\ref{lem:decreasing-gain} establishes the required uniform negligibility for stable decreasing-step-size recursions. Its proof is independent of the Poisson decomposition and may be useful for controlling similar endpoint-dependent remainders in other functional limit problems. 

\paragraph{The one-dimensional FCLT for inference.}
For inference on $\bm v^\top\bx^\star$, let $\mathcal W(\cdot)$ be a standard
one-dimensional Brownian motion and define
$B_T(r):=\bm v^\top\{\bPhi_T(r)-r\bPhi_T(1)\}$ and
$\mathcal B(r):=\mathcal W(r)-r\mathcal W(1)$.  Projecting the
$d$-dimensional FCLT along $\bm v$ gives
\begin{equation}
\bigl(\bm v^\top\bPhi_T(1),B_T(\cdot)\bigr)
\ \Longrightarrow\ \sigma_{\bm v}
\bigl(\mathcal W(1),\mathcal B(\cdot)\bigr).
\label{eq:intro-joint-fclt}
\end{equation}
Thus the scaled estimation error and the observed bridge carry the same unknown
scale $\sigma_{\bm v}$.

\paragraph{From the FCLT to a confidence interval.}
To convert~\eqref{eq:intro-joint-fclt} into an inference procedure, set
$C_0([0,1]):=\{f\in C([0,1],\R):f(0)=f(1)=0\}$ and consider a continuous
map $\mathcal D:C_0([0,1])\to[0,\infty)$ that is scale-equivariant,
$\mathcal D(cf)=|c|\mathcal D(f)$, and satisfies
$\Pp\{\mathcal D(\mathcal B)>0\}=1$.  The target $\bx^\star$ cancels from
$B_T(\cdot)$, so $\mathcal D(B_T)$ is observable.  Scale equivariance and
\eqref{eq:intro-joint-fclt} give
\begin{equation}
  \frac{\bm v^\top\bPhi_T(1)}{\mathcal D(B_T)}
  =\frac{\sqrt T\,\bm v^\top(\bar\bx_T-\bx^\star)}{\mathcal D(B_T)}
  \ \Longrightarrow\
  \frac{\mathcal W(1)}{\mathcal D(\mathcal B)}=:Z_{\mathcal D},
  \label{eq:pivot-limit}
\end{equation}
where the unknown scale $\sigma_{\bm v}$ has canceled.  If $c_{1-\alpha,\mathcal D}$ denotes the
$(1-\alpha)$ quantile of $|Z_{\mathcal D}|$, an asymptotic $(1-\alpha)$ confidence
interval for $\bm v^\top\bx^\star$ is
\begin{equation}
  \left[
  \bm v^\top\bar\bx_T-c_{1-\alpha,\mathcal D}
  \frac{\mathcal D(B_T)}{\sqrt T},
  \ \bm v^\top\bar\bx_T+c_{1-\alpha,\mathcal D}
  \frac{\mathcal D(B_T)}{\sqrt T}
  \right].
  \label{eq:main-CI}
\end{equation}

\paragraph{Contributions.}
We summarize our main contributions as follows.

\begin{enumerate}[leftmargin=*,label=(\roman*)]
  \item We establish an FCLT for the partial-sum process of nonlinear,
  decreasing-step-size SA with Markov data under a set of high-level regularity conditions.  We also provide primitive
  $V$-geometric conditions for verifying some requirements.

\item To establish the FCLT, we prove a lemma for decreasing-step-size recursions
(Lemma~\ref{lem:decreasing-gain}) that uniformly controls endpoint-dependent
remainders over all partial-sum endpoints. Its proof is independent of the
Poisson decomposition and may be useful in related functional limit arguments.

  \item We construct the single-trajectory confidence interval
  in~\eqref{eq:main-CI}, which avoids long-run covariance estimation, and study
  two families of $\mathcal D$: the $K$-dimensional polynomial-series
  normalizer $P_K$ and the $L_m$ bridge normalizers.

  \item We specialize the FCLT and confidence-interval construction to five
  applications and verify the required conditions under the stated
  application-specific assumptions.  We then compare the
  resulting inference procedures with two existing online baselines.  Across
  the experiments, the five-dimensional polynomial choice $P_5$ offers a
  favorable balance between
  coverage and interval length.
\end{enumerate}


\paragraph{Paper organization.}
The paper proceeds as follows.
Section~\ref{sec:related} reviews asymptotic theory for stochastic approximation and inference for SGD, Q-learning, and LoRA.
Section~\ref{sec:framework} introduces the assumptions and sufficient conditions for checking them.
Sections~\ref{sec:fclt} and~\ref{sec:inference} establish the functional limit
theorem and convert it into an online confidence interval.
Section~\ref{sec:examples} develops the applications, and Section~\ref{sec:experiments} presents the experiments.
Section~\ref{sec:proof-architecture} explains the proof architecture, Section~\ref{sec:conclusion} concludes, and the appendices give complete proofs and reproducibility details.
Code and numerical results are available at
\url{https://github.com/lx10077/nonlinear-sa-inference}.

\section{Related Work}
\label{sec:related}

\subsection{Central and Functional Limit Theorems for Stochastic Approximation}

Classical SA theory establishes convergence and asymptotic normality under
broad conditions
\citep{metivier1984applications,borkar2009stochastic,kushner2003stochastic,
benveniste2012adaptive}.  Within this theory, Polyak--Ruppert averaging attains
the optimal first-order covariance for broad classes of recursions
\citep{ruppert1988efficient,polyak1990new,polyak1992acceleration}, with the
corresponding local optimality characterized by convolution and minimax results
\citep{hajek1972local,van2000asymptotic,duchi2021asymptotic}.  Under Markov
sampling, serial dependence and conditional bias prevent a direct reduction to
independent-noise theory.  Poisson-equation decompositions provide the standard
route to a martingale approximation and hence to a Gaussian limit
\citep{andrieu2005stability}.

For SA with Markov noise, the closest classical antecedent is
\citet{liang2010trajectory}.  That work allows the transition kernel to depend
on the current iterate and establishes asymptotic normality and efficiency of
the Polyak--Ruppert average.  The object studied here is different: a
partial-sum FCLT for general nonlinear SA driven by a fixed Markov kernel,
which requires a uniform approximation over all partial-sum endpoints.

Classical functional limits study a different process: suitably rescaled
interpolations of the iterate trajectory itself
\citep{ljung1977analysis,borkar2021ode,liang2023asymptotic,
flamand2026functional}.  Such limits describe the local evolution of the
iterates rather than the cumulative estimation error captured by their
partial sums.  The present result concerns the latter for a single nonlinear,
decreasing-step-size recursion with Markov noise.

\subsection{Inference for Stochastic Gradient Descent}

Inference for the Polyak--Ruppert average requires handling its unknown
asymptotic covariance.  For independent data, recursive plug-in and online
batch-means estimators provide one solution
\citep{godichon2019online,chen2020statistical,zhu2021online,
singh2025utility}.  Overlapping batch means extend this approach to Markov
sampling \citep{roy2023online}.  A second group of procedures modifies or
replicates the SGD computation through fixed-step-size averages
\citep{li2018statistical}, perturbed or resampled recursion paths
\citep{fang2018online,lam2026cheap}, or hierarchical thread splitting
\citep{su2018uncertainty}.  For more general $\phi$-mixing observations,
\citet{liu2023mixing} combine mini-batch SGD with a mini-batch bootstrap.
Depending on the construction, these approaches require covariance estimation
and tuning, additional batching, or auxiliary recursion paths.

Another route constructs the scale directly from the observed iterate path.
\citet{zhu2021batchmeans} use a process-level FCLT and a fixed number of batches
to obtain covariance-free confidence regions for averaged SGD.  Building on
pivotal inference for time-series regression \citep{kiefer2000simple},
\citet{lee2021fast} use the full partial-sum path to develop random-scaling
inference for decreasing-step-size SGD with i.i.d. data.  This construction
avoids both covariance estimation and resampling.

FCLT-based inference has since been extended beyond standard single-machine SGD.  
\citet{li2025stream} establish an FCLT and online inference for smooth
stream SGD driven by parameter-dependent Markov data, focusing on queueing and inventory applications in which synchronous coupling verifies the required
Wasserstein-type regularity.  
\citet{li2021statistical} develop
communication-efficient random-scaling inference for Local SGD with
intermittent communication.  \citet{chen2021online} extend the approach to
gradient-free Kiefer--Wolfowitz optimization, where a plug-in covariance
estimator would require additional function evaluations to estimate a Hessian.
For nonsmooth problems, \citet{lee2022fast} establish the required path-level
limit for SGD in quantile regression.  These results exploit structure specific
to their respective optimization settings and do not cover general nonlinear
SA driven by a single Markov trajectory.

\subsection{Inference for Q-Learning}

Classical analyses of asynchronous Q-learning establish convergence from a
single trajectory but do not provide uncertainty quantification
\citep{tsitsiklis1994asynchronous,even2003learning}.  Early inferential work in
reinforcement learning focused primarily on policy evaluation
\citep{white2010interval,hanna2017bootstrapping,hao2021bootstrapping}.  In
particular, \citet{ramprasad2021online} establish distributional consistency of
the online bootstrap for temporal-difference and gradient temporal-difference
learning.  Both algorithms are linear SA.  Their validity
result therefore does not cover Q-learning, as the Bellman optimality operator
is nonlinear.

The closest prior random-scaling result for Q-learning is
\citet{li2021polyak}, who establish a partial-sum FCLT and covariance-free
inference for synchronous tabular Q-learning.  Their generative model supplies
independent observations for every state--action pair at each iteration.  In
contrast, the setting considered here updates only the visited state--action
pair and generates all observations along a single Markov trajectory.  Their
analysis therefore does not encounter coordinate asynchrony or serial
dependence.

More recently, \citet{bonnerjee2026sharp} derive strong Gaussian approximations
for synchronous tabular Q-learning and suggest a possible Gaussian-bootstrap
construction, but do not develop and validate a complete confidence-interval
procedure.  Their analysis also assumes synchronous generative-model access
and is specific to Q-learning.  Finally,
\citet{rubtsov2026gaussian} study Gaussian approximation for
entropy-regularized asynchronous Q-learning with linear function approximation,
but do not establish the partial-sum FCLT required for random-scaling
inference.

\subsection{Low-Rank Adaptation and Its Inference}

LoRA fine-tunes a pretrained model by representing a trainable low-rank change
to a weight matrix as the product $\bW=\bB\bA$ of two factors
\citep{hu2022lora}.  This parameterization creates two related difficulties.
The joint factor update is nonlinear, and the transformation
$(\bB,\bA)\mapsto(\bB\bQ,\bQ^{-1}\bA)$ leaves $\bW$ unchanged, so the individual
factors are not identified.  Existing theory has addressed these difficulties
from an optimization perspective.  \citet{mu2026lora} obtain a nonasymptotic
convergence analysis of the original LoRA gradient-descent update by lifting
the factors to an outer-product representation and controlling higher-order
discretization terms.  \citet{castin2026balancedlora} instead project the
factors onto a balanced manifold, preserving the low-rank product while
enforcing $\bB^\top\bB=\bA\bA^\top$ to improve conditioning.  A separate line
uses Bayesian approximations over the LoRA parameters to quantify predictive
uncertainty and improve calibration
\citep{yang2024bayesianlora}.

For the inferential goal considered here, the stochastic factor updates induce
a recursion for the identified product $\bW$, which is itself the target of
inference.  This representation enables confidence intervals for prespecified
linear functionals of the fitted low-rank update under Markov sampling.  Related
work by \citet{gao2026statlora} instead studies component testing and rank
allocation for LoRA optimizer trajectories.  To our knowledge, an
averaged-iterate FCLT and a covariance-free confidence interval for the fitted
product have not previously been established.

\section{General Assumptions}
\label{sec:framework}

We now state the assumptions used to establish the FCLT.  

\begin{assumption}[Local Mean Field]
\label{ass:local}
The target $\bx^\star$ satisfies $\bg(\bx^\star)=\bzero$.  On a neighborhood $\cN$ of $\bx^\star$,
\[
  \bg(\bx)=\bG(\bx-\bx^\star)+\br_g(\bx),
  \qquad
  \frac{\norm{\br_g(\bx)}}{\norm{\bx-\bx^\star}}\longrightarrow0
  \quad\text{as }\bx\to\bx^\star,
\]
where $\bG\in\R^{d\times d}$ is the local Jacobian and
$\br_g(\bx):=\bg(\bx)-\bG(\bx-\bx^\star)$ is its first-order remainder.
    Every eigenvalue of $\bG$ has strictly positive real part.
\end{assumption}

Assumption~\ref{ass:local} gives a first-order approximation to the mean dynamics
near $\bx^\star$.  With $\bm e_t=\bx_t-\bx^\star$, the mean recursion is locally approximated by $\bm e_{t+1}=(\bI-\eta_t\bG)\bm e_t$.
The eigenvalue condition makes this local dynamics stable, while $\br_g(\bx)$
is negligible to first order. In gradient-based applications, $\bG$ is the
Hessian of the population objective $\bm g$ at $\bx^\star$.

\begin{remark}[Local rather than global regularity]
Local differentiability and this eigenvalue condition are standard in asymptotic SA theory
\citep{polyak1992acceleration,benveniste2012adaptive,chen2002stochastic}.
The first-order approximation is needed only through the accumulated-remainder
condition in Assumption~\ref{ass:path}, so no corresponding regularity condition
is imposed on $\bg$ away from the target.
We also do not require $\bG$ to be symmetric or diagonalizable, thereby allowing the nonsymmetric dynamics that often arise in asynchronous updates and function approximation.
\end{remark}

\begin{assumption}[Poisson Regularity]
\label{ass:poisson}
The chain $(\xi_t)_{t \ge 0}$ is positive Harris recurrent and aperiodic, while its initial distribution need not be $\mu$.  The $\mu$-centered Poisson solution $\bU(\bx,\cdot)$ (defined in \eqref{eq:poisson-family}) exists for every $\bx\in\R^d$.
For some $p>2$, a measurable envelope function $L(\cdot):\cX\to[1,\infty)$, and a constant $C<\infty$, the following bound holds for all $\bx,\by\in\R^d$:
\begin{equation}
  \left\{\cP\norm{\bU(\bx,\cdot)-\bU(\by,\cdot)}^p(\xi)\right\}^{1/p}
  \le C L(\xi)\norm{\bx-\by}.
  \label{eq:u-weighted-lipschitz}
\end{equation}
Moreover,
\begin{equation}
  \sup_{t\ge0}\E\!\left[L(\xi_t)^p
  +\norm{\bU(\bx^\star,\xi_t)}^p\right]<\infty.
  \label{eq:poisson-moments}
\end{equation}
\end{assumption}

Assumption~\ref{ass:poisson} formalizes the Poisson-equation requirement
introduced in Section~\ref{sec:introduction}.  The solution $\bU$ accounts for
dependence in the observations.  We impose a global weighted Lipschitz bound in
$\bx$ and bounded $p>2$ moments at $\bx^\star$ along the observed trajectory.
The global bound makes the Poisson decomposition valid along the entire iterate
path, without assuming almost-sure convergence of the iterates or their eventual
localization in $\cN$. Finite-state applications can verify these conditions
directly. As an example, Proposition~\ref{prop:v-route} in Section~\ref{sec:verify} provides a $V$-geometric
route for unbounded state spaces.

\begin{remark}[Relation to earlier conditions]
Poisson equations (indexed by the parameter) are classical tools for SA with Markov
noise \citep{benveniste2012adaptive,andrieu2005stability}.
Earlier analyses commonly derive regularity of the Poisson solution from various drift,
mixing, and update conditions.  We instead state directly the global weighted
regularity used to compare the noise along the iterate path with its value at
$\bx^\star$.  Proposition~\ref{prop:v-route} gives a primitive sufficient route.
\end{remark}

\begin{assumption}[Slowly Decreasing Step Size]
\label{ass:gain}
The deterministic step sizes are positive and nonincreasing, and
\begin{equation}
  \eta_t\to0,\qquad t\eta_t\to\infty,\qquad
  \frac{\eta_{t-1}-\eta_t}{\eta_{t-1}^2}\to0,
  \qquad \sum_{t=0}^{\infty}\eta_t^2<\infty.
  \label{eq:gain-conditions}
\end{equation}
\end{assumption}

Assumption~\ref{ass:gain} balances two requirements on the step size.  On one hand, it must
vanish more slowly than $1/t$ and change gradually between successive
iterations.  On the other hand, the condition $\sum_t\eta_t^2<\infty$ limits the cumulative
effect of random fluctuations in the updates
\citep{polyak1992acceleration,chen2002stochastic,kushner2003stochastic}.
The canonical choice $\eta_t=\eta(t+1)^{-\kappa}$ satisfies the assumption for
every $\kappa\in(1/2,1)$.  

\begin{assumption}[Path Stability]
\label{ass:path}
With $\bm e_t=\bx_t-\bx^\star$ and the $p>2$ in Assumption~\ref{ass:poisson}, the following conditions hold:
\begin{enumerate}[leftmargin=*,label=(\roman*)]
  \item the following uniform moment bound holds:
  \begin{equation}
    \sup_{t\ge1}\E\!\left[\norm{\bH(\bx_t,\xi_t)}^p
    +L(\xi_{t-1})^p\norm{\bm e_t}^p\right]<\infty;
    \label{eq:path-moments}
  \end{equation}
  \item the iterates are $L^2$ consistent, that is,
  \begin{equation}
    \frac1T\sum_{t=1}^T\E\norm{\bm e_t}^2\longrightarrow0;
    \label{eq:time-averaged-l2}
  \end{equation}
  \item the accumulated local nonlinearity is negligible:
  \begin{equation}
    \frac1{\sqrt T}\sum_{t=1}^T\norm{\br_g(\bx_t)}\xrightarrow{p}0.
    \label{eq:integrated-nonlinearity}
  \end{equation}
\end{enumerate}
\end{assumption}

Assumption~\ref{ass:path} states only the trajectory properties needed by the
FCLT without prescribing conditions for global convergence.  In Condition~(i),
\eqref{eq:path-moments} controls the $p$th moments of the update and the weighted
error appearing in~\eqref{eq:u-weighted-lipschitz}.  The requirement $p>2$
provides the moment margin needed for the martingale CLT and for uniform control
over the partial-sum path.  Condition~(ii) requires the iterates
$\{\bx_t\}_{t\ge0}$ to approach $\bx^\star$ on average in mean square.
Together with Condition~(i) and Poisson regularity, this average
convergence allows the Markov-noise term along $\bx_t$ to be approximated
uniformly over partial sums by its counterpart at $\bx^\star$.
Condition~(iii) requires the accumulated nonlinear remainder
$\sum_{t\le T}\norm{\br_g(\bx_t)}$ to be $o_p(\sqrt T)$.  The next subsection
and the application propositions give concrete sufficient conditions.


\subsection{Primitive Verification Routes}
\label{sec:verify}

This subsection provides sufficient conditions for verifying
Assumptions~\ref{ass:poisson} and~\ref{ass:path}.  These are verification tools
rather than additional assumptions of the FCLT.  Proposition~\ref{prop:v-route} verifies
Poisson regularity from primitive conditions on the chain and update function.
The joint route that follows combines it with familiar convergence and
mean-square error bounds to verify path stability.  
Let $V:\cX\to[1,\infty)$ be measurable.  For a measurable function $f$ and a
signed measure $\lambda$, define
$\norm{f}_{V}:=\sup_\xi\norm{f(\xi)}/V(\xi)$ and
$\norm{\lambda}_{V}:=\sup_{|f|\le V}|\int f\,\mathrm{d}\lambda|$, respectively.

\begin{definition}[$V$-geometric ergodicity {\citep{meyn2012markov}}]
\label{def:v-geometric}
A Markov kernel $\cP$ with stationary distribution $\mu$ is $V$-geometrically
ergodic if, for some $C<\infty$ and $\rho\in(0,1)$,
\[
  \norm{\cP^k(\xi,\cdot)-\mu}_{V}
  \le C\rho^k V(\xi),
  \qquad \xi\in\cX,\quad k\ge0.
\]
\end{definition}

\begin{proposition}[$V$-geometric sufficient condition]
\label{prop:v-route}
Let $p>2$ and let $V:\cX\to[1,\infty)$ be measurable.  Suppose the chain
$\{\xi_t\}_{t\ge0}$ is
positive Harris recurrent, aperiodic, and $V^{1/p}$-geometrically ergodic in
the sense of Definition~\ref{def:v-geometric}.  Assume that
$\int V\,\mathrm{d}\mu<\infty$, $\sup_{t\ge0}\E V(\xi_t)<\infty$, and, for some
$C<\infty$ and all $\bx,\by\in\R^d$,
\begin{equation}
  \norm{\bH(\bx^\star,\cdot)}_{V^{1/p}}\le C,
  \qquad
  \norm{\bH(\bx,\cdot)-\bH(\by,\cdot)}_{V^{1/p}}
  \le C\norm{\bx-\by}.
  \label{eq:h-v-lipschitz}
\end{equation}
Then, for every $\bx\in\R^d$, the resolvent series
\begin{equation}
  \bU(\bx,\xi)=\sum_{k=0}^{\infty}
  \{\cP^k\bH(\bx,\xi)-\bg(\bx)\}
  \label{eq:poisson-resolvent}
\end{equation}
is absolutely convergent under $\norm{\cdot}_{V^{1/p}}$. Its sum $\bU(\bx,\xi)$ solves~\eqref{eq:poisson-family}, and the resulting $\bU(\cdot, \cdot)$ satisfies
\eqref{eq:u-weighted-lipschitz} with
$L(\xi)\le C\{1+(\cP V)(\xi)^{1/p}\}$.  Hence, the assumed uniform $V$ moment implies Assumption~\ref{ass:poisson}.  
\end{proposition}

The proof is in Appendix~\ref{app:poisson}.  Proposition~\ref{prop:v-route}
shows that $V$-geometric ergodicity and the target-envelope and global
parameter-Lipschitz bounds in~\eqref{eq:h-v-lipschitz} yield the Poisson
solution, the weighted Lipschitz bound~\eqref{eq:u-weighted-lipschitz}, and the
moment condition~\eqref{eq:poisson-moments} required by
Assumption~\ref{ass:poisson}.  The geometric bound
in Definition~\ref{def:v-geometric} ensures that the terms in the resolvent series \eqref{eq:poisson-resolvent} decay,
while~\eqref{eq:h-v-lipschitz} allows $\bH(\bx,\xi)$ to grow with the state
$\xi$ at a rate controlled by $V(\xi)^{1/p}$.  The uniform $V$-moment bound
then implies the required $p$th-moment control along the observed trajectory.
Finite-state applications may take $V\equiv1$, whereas a nonconstant $V$
allows unbounded Markov states and covariates.  

\paragraph{Joint Verification of Poisson Regularity and Path Stability.}
Proposition~\ref{prop:v-route} can be combined with familiar estimates for the
iterates to verify Assumptions~\ref{ass:poisson} and~\ref{ass:path}.  
If available
$p$th-moment estimates establish the uniform bound in
\eqref{eq:path-moments}, then Assumption~\ref{ass:poisson} and Condition~(i) of
Assumption~\ref{ass:path} follow.  A mean-square error bound
$\E\norm{\bx_t-\bx^\star}^2\le C\eta_t$ further implies
\eqref{eq:time-averaged-l2}.  Finally, Condition~(iii) follows whenever the
nonlinear remainder satisfies
$\norm{\br_g(\bx)}\le C\norm{\bx-\bx^\star}^{2}$ along the iterate path and
\begin{equation}
 T^{-1/2}\sum_{t=1}^T
 \norm{\bx_t-\bx^\star}^{2}\xrightarrow{p}0.
 \label{eq:remainder-occupation-route}
\end{equation}
Thus, verification reduces to moment, mean-square, and nonlinear-occupation
bounds for the iterate path.  Finite-sample SA analyses offer one way to obtain
such bounds \citep{chen2020finite,chen2021lyapunov,roy2023online}.  The
applications below use model-specific arguments based on contraction,
projection, finite-time bounds, or direct analysis.

\section{Functional Central Limit Theorem}
\label{sec:fclt}

We now state the main theoretical results.  Let $\bcalW(\cdot)$ be a standard
$d$-dimensional Brownian motion, and let $\bS$ denote the long-run covariance of
the update noise at $\bx^\star$ under the stationary Markov chain.  Since
$\bg(\bx^\star)=\bzero$, it can be written as
\begin{equation}
 \bS=\lim_{T\to\infty}\frac1T\E_\mu\!\left[
 \left\{\sum_{t=1}^T\bH(\bx^\star,\xi_t)\right\}
 \left\{\sum_{t=1}^T\bH(\bx^\star,\xi_t)\right\}^{\!\top}
 \right].
 \label{eq:long-run-covariance}
\end{equation}
The noise conditions in Assumption~\ref{ass:poisson} ensure that this limit
exists.  To state the pathwise approximation below, define the Poisson
martingale difference
\begin{equation}
 \bepsilon_t=\cP\bU(\bx_t,\xi_{t-1})-\bU(\bx_t,\xi_t).
 \label{eq:poisson-innovations}
\end{equation}
The Poisson solution in Assumption~\ref{ass:poisson} is defined for every
$\bx$, so $\bepsilon_t$ is well defined along the entire iterate path.
The first theorem upgrades the usual Polyak--Ruppert CLT to the entire
partial-sum path. All of the proofs for inference results are in Appendix \ref{app:fclt-proof}.

\begin{theorem}[Functional central limit theorem]
\label{thm:fclt}
For the linearly interpolated process $\bPhi_T(\cdot)$ defined
in~\eqref{eq:partial-sum-process}, under
Assumptions~\ref{ass:local}--\ref{ass:path},
\begin{equation}
  \bPhi_T(\cdot) \Rightarrow \bG^{-1}\bS^{1/2}\bcalW(\cdot)
  \quad\text{in } C([0,1],\R^d)
  \text{ equipped with the uniform norm.}
  \label{eq:fclt}
\end{equation}
More precisely, we have
\begin{equation}
  \sup_{0\le r\le1}
  \norm{\bPhi_T(r)-\bG^{-1}\frac{1}{\sqrt{T}}\left\{
  \sum_{t=1}^{\lfloor Tr\rfloor}\bepsilon_t
  +(Tr-\lfloor Tr\rfloor)\bepsilon_{\lfloor Tr\rfloor+1}\right\}}
  \xrightarrow{p}0.
  \label{eq:path-linear-representation}
\end{equation}
\end{theorem}

At $r=1$, Theorem~\ref{thm:fclt} recovers the Polyak--Ruppert CLT
\[
  \bPhi_T(1)  =
  \sqrt{T}(\bar\bx_T-\bx^\star)
  \Rightarrow \mathcal N\!\left(\bzero,\bG^{-1}\bS\bG^{-\top}\right).
\]
Under the classical regularity conditions for averaged-SA efficiency, this endpoint $  \bPhi_T(1)$
has the optimal first-order covariance
\citep{ruppert1988efficient,polyak1992acceleration,duchi2021asymptotic}.  The FCLT
is strictly stronger than this endpoint statement: it gives joint convergence at
all time fractions and asymptotic equicontinuity, thereby justifying continuous
bridge, integral, and supremum functionals.  

Many matrix and manifold algorithms add higher-order discretization terms to an
exact first-order SA update.  For any fixed integer $J\ge2$, the next result
covers perturbations of the form
$\sum_{k=2}^{J}\eta_t^k\bzeta_{k,t}$.  The random vectors
$\bzeta_{k,t}\in\R^d$ are $\cF_t$-measurable and satisfy
$\sup_{2\le k\le J,\,t\ge0}\E\norm{\bzeta_{k,t}}^p<\infty$.  Bounded step sizes reduce this
perturbation to the stated order-$\eta_t^2$ form.

\begin{proposition}[Second-order perturbations]
\label{prop:second-order-perturbation}
Consider the perturbed recursion
\begin{equation}
  \bx_{t+1}=\bx_t-\eta_t\bH(\bx_t,\xi_t)+\eta_t^2\bzeta_t,
  \label{eq:perturbed-sa-recursion}
\end{equation}
where $\bzeta_t$ is $\cF_t$-measurable and, for the exponent $p>2$ in
Assumption~\ref{ass:poisson}, $\sup_t\E\norm{\bzeta_t}^p<\infty$.  If
the sequence generated by~\eqref{eq:perturbed-sa-recursion} satisfies
Assumptions~\ref{ass:local}--\ref{ass:path}, then
all conclusions of Theorem~\ref{thm:fclt} hold.
\end{proposition}

The perturbation term $\bzeta_t$ need not be mean-zero.  Its contribution to the
averaged-iterate path is of order
$T^{-1/2}\sum_{t\le T}\eta_t\norm{\bzeta_t}$, which vanishes under
Assumption~\ref{ass:gain}.  We next extend the FCLT from linear contrasts of the
SA state to smooth scalar functionals.  This allows the inferential target to be
a smooth function of the coordinates used to analyze the recursion.

\begin{corollary}[Smooth scalar functionals]
\label{cor:smooth-functional}
Suppose the FCLT in~\eqref{eq:fclt} holds for $\bPhi_T$.  Let
$\psi:\R^d\to\R$ be differentiable at $\bx^\star$, put
$\bm v=\nabla\psi(\bx^\star)$, and assume
\begin{equation}
  \frac1{\sqrt T}\sum_{t=1}^T
  \left|\psi(\bx_t)-\psi(\bx^\star)
  -\bm v^\top(\bx_t-\bx^\star)\right|\xrightarrow{p}0.
  \label{eq:functional-occupation}
\end{equation}
Define the linearly interpolated functional path by
\[
 \Psi_T^\psi(r)=\frac1{\sqrt T}\!\left\{
 \sum_{t=1}^{\lfloor Tr\rfloor}[\psi(\bx_t)-\psi(\bx^\star)]
 +(Tr-\lfloor Tr\rfloor)
 [\psi(\bx_{\lfloor Tr\rfloor+1})-\psi(\bx^\star)]\right\}.
\]
If $\bm v^\top\bG^{-1}\bS\bG^{-\top}\bm v>0$, then
$\Psi_T^\psi(\cdot)\Rightarrow\sigma_\psi\mathcal W(\cdot)$ in
$C([0,1],\R)$ under the uniform norm, where
$\sigma_\psi^2=\bm v^\top\bG^{-1}\bS\bG^{-\top}\bm v$.
\end{corollary}

\begin{remark}[Linear or Smooth Functionals]
If $\psi$ is linear, i.e., $\psi(\bx)=\bm v^\top\bx$, the left-hand side
of~\eqref{eq:functional-occupation} is identically zero.  For a nonlinear
$\psi$, differentiability gives only the pointwise relation
$\psi(\bx)-\psi(\bx^\star)-\bm v^\top(\bx-\bx^\star)
=o(\norm{\bx-\bx^\star})$ as $\bx\to\bx^\star$, which need not control the sum
in~\eqref{eq:functional-occupation}.  If the gradient of $\psi$ is Lipschitz, it suffices that
$T^{-1/2}\sum_{t=1}^T\norm{\bx_t-\bx^\star}^2\xrightarrow{p}0$.
\end{remark}

\section{Online Self-Normalized Inference}
\label{sec:inference}

The FCLT provides the joint path information needed to construct an interval
from a single trajectory.  In view of Corollary~\ref{cor:smooth-functional}, we
present the construction for a prespecified linear functional.  The same
construction applies to smooth scalar functionals
satisfying~\eqref{eq:functional-occupation}.

To that end, we fix a nonzero vector
$\bm v\in\R^d$ and define
$Y_t=\bm v^\top\bx_t$, $\vartheta^\star=\bm v^\top\bx^\star$, and
$\bar Y_T=T^{-1}\sum_{t=1}^T Y_t$.  Denote a standard Brownian bridge by
$\mathcal B(r):=\mathcal W(r)-r\mathcal W(1)$.  The construction is based on
the observed bridge
\[
  B_T(r):=\bm v^\top\{\bPhi_T(r)-r\bPhi_T(1)\}
  =\frac{1}{\sqrt{T}}
  \bigg[\sum_{t=1}^{\lfloor Tr\rfloor}Y_t
  +(Tr-\lfloor Tr\rfloor)Y_{\lfloor Tr\rfloor+1}
  -r\sum_{t=1}^TY_t\bigg].
\]
Although $\bPhi_T$ is centered at the unknown target, $\bx^\star$ cancels from
$B_T$, making the bridge fully observable.  In the notation above, the joint
convergence in~\eqref{eq:intro-joint-fclt} reads
\begin{equation*}
	\bigl(\sqrt T(\bar Y_T-\vartheta^\star),B_T(\cdot)\bigr)
	\Rightarrow \sigma_{\bm v}\bigl(\mathcal W(1),\mathcal B(\cdot)\bigr).
	\tag{\ref{eq:intro-joint-fclt}}
\end{equation*}
Thus the estimation error and the observed bridge share the same unknown scale.

We next construct from $B_T$ a normalizer that carries this same scale and
cancels it when used as the denominator of the statistic.
Let $C_0([0,1])=\{f\in C([0,1],\R):f(0)=f(1)=0\}$.
Consider a continuous functional $\mathcal D:C_0([0,1])\to[0,\infty)$ that is
scale-equivariant, meaning $\mathcal D(cf)=|c|\mathcal D(f)$, and is
positive almost surely at $\mathcal B$.

\begin{corollary}[Bridge-Functional Coverage]
\label{cor:coverage}
Under the conditions of Theorem~\ref{thm:fclt}, assume the chosen vector has nonzero
limiting variance
$\sigma_{\bm v}^2:=\bm v^\top\bG^{-1}\bS\bG^{-\top}\bm v>0$.
Then
\begin{equation*}
  Z_{\mathcal D,T}(\vartheta^\star)
  :=\frac{\sqrt T(\bar Y_T-\vartheta^\star)}{\mathcal D(B_T)}
  \Rightarrow\frac{\mathcal W(1)}{\mathcal D(\mathcal B)}
  =:Z_{\mathcal D}.
  \tag{\ref{eq:pivot-limit}}
\end{equation*}
If $c_{1-\alpha,\mathcal D}$ satisfies
$\Pp(|Z_{\mathcal D}|\le c_{1-\alpha,\mathcal D})=1-\alpha$, then
the interval
\begin{equation*}
\left[
\bar Y_T-c_{1-\alpha,\mathcal D}\frac{\mathcal D(B_T)}{\sqrt T},
\bar Y_T+c_{1-\alpha,\mathcal D}\frac{\mathcal D(B_T)}{\sqrt T}
\right]
\tag{\ref{eq:main-CI}}
\end{equation*}
has coverage
$\Pp\{\vartheta^\star\in\mathcal I_{\mathcal D,T}(1-\alpha)\}\to1-\alpha$ as $T \to \infty$.
\end{corollary}

Corollary~\ref{cor:coverage} formalizes this scale cancellation.  The interval
uses only the scalar sequence
$Y_t=\bm v^\top\bx_t$ and requires neither the Jacobian $\bG$ nor the long-run
covariance $\bS$.  This is the main inferential benefit of controlling the entire
partial-sum path.

\paragraph{Expected-Width Benchmark.}
The interval in \eqref{eq:main-CI} has half-width
$c_{1-\alpha,\mathcal D}\mathcal D(B_T)/\sqrt T$.  After removing the common
factor $\sigma_{\bm v}/\sqrt T$, its limiting value is
$c_{1-\alpha,\mathcal D}\mathcal D(\mathcal B)$ due to \eqref{eq:intro-joint-fclt}.  Accordingly, we use the
normalized expected asymptotic half-width
$\ell_\alpha(\mathcal D):=c_{1-\alpha,\mathcal D}
\E\mathcal D(\mathcal B)$ as the benchmark for comparing bridge functionals.

\begin{proposition}[Expected-Width Benchmark]
\label{prop:width-benchmark}
Let $\Phi$ denote the standard normal distribution function.  
Every functional $\mathcal{D}$ covered by Corollary~\ref{cor:coverage} with finite mean satisfies
\[
  \ell_\alpha(\mathcal D)\ge \Phi^{-1}(1-\alpha/2).
\]
\end{proposition}

This lower bound is attained by the textbook Gaussian interval when
$\sigma_{\bm v}$ is known or can be consistently estimated.  The
covariance-free intervals considered in \eqref{eq:main-CI} instead use a nondegenerate random
normalizer and therefore incur asymptotic width inflation.  We use the bound
as a reference for comparing this inflation across the normalizers below.  We
next introduce two concrete choices of $\mathcal D$.  The symbols $P_K$ and
$L_m$ serve as labels for the corresponding normalizer families.

\subsection{Polynomial-Series Normalizer}

Fixed-smoothing inference uses a fixed number of whitened orthonormal-series
projections to obtain standard Student or $F$ limits \citep{sun2013orthonormal}.
We apply this principle to the bridge generated by nonlinear Markov SA and use
polynomial coordinates.  For a fixed integer $K\ge1$, $P_K$ labels the
polynomial-series normalizer based on $K$ projections.  Define
\[
  a_{j,T}:=\int_0^1 B_T(r)r^{j-1}\,\mathrm{d}r,
  \quad j=1,\ldots,K,
  \qquad
  \mathcal D_{P_K,T}^2:=K^{-1}\bm a_T^\top\bm\Sigma_K^{-1}\bm a_T,
\]
where $\bm a_T:=(a_{1,T},\ldots,a_{K,T})^\top$ and
\begin{equation}
  (\bm\Sigma_K)_{jk}
  :=\frac{1}{(j+1)(k+1)(j+k+1)},
  \qquad 1\le j,k\le K.
  \label{eq:series-covariance}
\end{equation}
The matrix $\bm\Sigma_K$ is the covariance of the corresponding polynomial
projections of a standard Brownian bridge.  Multiplication by
$\bm\Sigma_K^{-1/2}$ puts these projections on a common scale.

\begin{proposition}[Fixed-$K$ Series Interval]
\label{prop:series-interval}
Under the conditions of Corollary~\ref{cor:coverage}, for every fixed $K$,
\[
  \frac{\mathcal D_{P_K,T}^2}{\sigma_{\bm v}^2}\Rightarrow\frac{\chi_K^2}{K},
  \qquad
  \frac{\sqrt T(\bar Y_T-\vartheta^\star)}{\mathcal D_{P_K,T}}\Rightarrow t_K.
\]
Here, $t_K$ denotes the Student $t$ distribution with $K$ degrees of freedom.
Consequently, with $t_{K,1-\alpha/2}$ denoting its $(1-\alpha/2)$ quantile, the interval
\begin{equation}
  \left[
  \bar Y_T-t_{K,1-\alpha/2}\frac{\mathcal D_{P_K,T}}{\sqrt T},
  \bar Y_T+t_{K,1-\alpha/2}\frac{\mathcal D_{P_K,T}}{\sqrt T}
  \right]
  \label{eq:confidence-interval}
\end{equation}
has asymptotic coverage $1-\alpha$.  Its normalized expected asymptotic
half-width is
\[
  \ell_{\alpha,K}
  =t_{K,1-\alpha/2}\sqrt{\frac2K}
  \frac{\Gamma((K+1)/2)}{\Gamma(K/2)}
  \longrightarrow \Phi^{-1}(1-\alpha/2)
  \quad\text{as }K\to\infty.
\]
\end{proposition}

\begin{remark}[Online Implementation]
\label{rem:series-online}
Appendix~\ref{app:series-scaler} provides an asymptotically equivalent implementation that uses $O(K)$ operations and $2K+1$ scalar
accumulators per iteration without changing the coverage conclusion.  It is used in our experiments.
\end{remark}

In the main experiments, we use the interval induced by $P_5$, which is based
on five polynomial projections.  Larger fixed values of $K$ produce shorter
asymptotic confidence intervals but can yield less stable finite-sample coverage.
At $95\%$ confidence, the normalized expected asymptotic half-width is $2.446$
for the interval induced by $P_5$, compared with the lower bound $1.960$.  The ablation in
Appendix~\ref{app:series-ablation} compares $P_5$, $P_{10}$, and $P_{20}$ and
supports the choice of $P_5$.  Proposition~\ref{prop:series-interval} covers any
fixed $K$ but does not allow $K$ to increase with $T$.

\subsection{\texorpdfstring{$L_m$}{L-m} Bridge Normalizers}

A second choice uses the full bridge rather than finitely many projections.
For $m\ge1$, $L_m$ labels the bridge-norm normalizer
\[
  \mathcal D_{L_m,T}:=\left\{\int_0^1|B_T(r)|^m\,\mathrm{d}r\right\}^{1/m},
  \qquad
  Z_{L_m,T}(\vartheta):=
  \frac{\sqrt T(\bar Y_T-\vartheta)}{\mathcal D_{L_m,T}}.
\]

\begin{proposition}[Fixed-$m$ Bridge-Norm Interval]
\label{prop:lm-interval}
Under the conditions of Corollary~\ref{cor:coverage}, for every fixed $m\ge1$,
\[
  Z_{L_m,T}(\vartheta^\star)\Rightarrow
  Z_{L_m}:=
  \frac{\mathcal W(1)}
  {\{\int_0^1|\mathcal B(r)|^m\,\mathrm{d}r\}^{1/m}}.
\]
Consequently, if $c_{1-\alpha,L_m}$ is the $(1-\alpha)$ quantile of
$|Z_{L_m}|$, then the interval
\begin{equation}
\left[
    \bar Y_T-c_{1-\alpha,L_m}\frac{\mathcal D_{L_m,T}}{\sqrt T},
    \bar Y_T+c_{1-\alpha,L_m}\frac{\mathcal D_{L_m,T}}{\sqrt T}
  \right]
  \label{eq:lm-confidence-interval}
\end{equation}
has asymptotic coverage $1-\alpha$.
\end{proposition}

\begin{remark}[Online Implementation]
\label{rem:lm-online}
For every fixed even integer $m$, Appendix~\ref{app:lm-baselines} constructs an
online normalizer $\widehat{\mathcal D}_{m,T}$ satisfying
$\widehat{\mathcal D}_{m,T}-\mathcal D_{L_m,T}=o_p(1)$.  It uses $m+2$
scalar running sums and $O(m)$ operations per iteration; in particular, $L_2$
uses four sums and $L_4$ uses six.  Replacing $\mathcal D_{L_m,T}$ by
$\widehat{\mathcal D}_{m,T}$ preserves the coverage conclusion.
\end{remark}

Unlike the series normalizer, which has a standard Student critical value, the
$L_m$ pivot has a nonstandard but nuisance-free distribution whose critical
values can be simulated once from a standard Brownian bridge.  We use the
interval induced by the familiar $L_2$ functional from random-scaling inference
and report the interval induced by $L_4$ as a
sensitivity check \citep{lee2021fast}.  Appendix~\ref{app:lm-baselines} reports
the critical values used for $L_2$ and $L_4$.

\section{Applications}
\label{sec:examples}

This section applies the framework to five applications: three Q-learning
recursions in Sections~\ref{sec:tabular-q}--\ref{sec:soft-q} and two
stochastic-gradient recursions in Sections~\ref{sec:markov-glm}--\ref{sec:lora}.
For each, we define the update and scalar inferential target and state the
resulting guarantee.  
Specifically, $\vartheta^\star$ denotes the application-specific target and $Y_t$ the corresponding scalar data computed from the current iterate.  Both are defined before the formal statement, whose proofs are in
Appendix~\ref{app:examples}.

\subsection{Asynchronous Tabular Q-Learning}
\label{sec:tabular-q}

Consider a discounted finite Markov decision process (MDP) with state space
$\cS$, action space $\cA$, transition law $p(s'\mid s,a)$, bounded reward
$|R_t|\le R_{\max}$, and discount factor $\gamma\in(0,1)$.  A stationary policy
$\pi$ selects $A_t\sim\pi(\cdot\mid S_t)$, after which a random reward is received and a new state is reached:
$R_t\sim p_R(\cdot\mid S_t,A_t)$ and
$S_{t+1}\sim p(\cdot\mid S_t,A_t)$.  The optimal action-value function is
$Q^\star(s,a):=\sup_\pi\E_\pi[\sum_{k=0}^\infty\gamma^kR_{t+k}
\mid S_t=s,A_t=a]$, the largest expected discounted return from taking action
$a$ in state $s$ and following an optimal policy thereafter.  Here $\E_\pi$
denotes expectation over the trajectory generated under $\pi$.

We next describe the asynchronous Q-learning used to
estimate $Q^\star$ from one trajectory \citep{watkins1989learning}.  Define
$d_Q:=|\cS||\cA|$ and fix an ordering of the state--action pairs.  Let
$\bq_t\in\R^{d_Q}$ be the vectorized action-value table at iteration $t$, and
denote its $(s,a)$ coordinate by $q_t(s,a)$.  A fixed behavior policy $\pi_b$
generates the single observed trajectory according to the preceding dynamics.
Thus $\{(S_t,A_t)\}_{t\ge0}$ is a Markov chain on $\cS\times\cA$ with transition probability
$p(s'\mid s,a)\pi_b(a'\mid s')$ from $(s,a)$ to $(s',a')$.
The update is
asynchronous because iteration $t$ changes only the coordinate associated with
the visited pair $(S_t,A_t)$, whereas synchronous Q-learning updates every
state--action coordinate using generative-model observations.  The asynchronous
update is
\begin{equation}
 q_{t+1}(s,a)=q_t(s,a)-\eta_t\mathbf{1}\{(s,a)=(S_t,A_t)\} \left[q_t(s,a)-R_t-\gamma\max_{a'}q_t(S_{t+1},a')\right].
 \label{eq:tabular-q-update}
\end{equation}
To identify the mean field $\bm g$, let $d_{\pi_b}$ be the stationary state distribution.
The diagonal matrix $\bD_{\pi_b}\in\R^{d_Q\times d_Q}$ has entries
$d_{\pi_b}(s)\pi_b(a\mid s)$.  Define the mean-reward vector $\br\in\R^{d_Q}$ by $r(s,a)=\E(R_t\mid S_t=s,A_t=a)$.  Let $\bP\in\R^{d_Q\times|\cS|}$ be the state--action to next-state transition matrix, with entry $p(s'\mid s,a)$, and define $\mathcal M:\R^{d_Q}\to\R^{|\cS|}$ by $(\mathcal M\bq)(s):=\max_a q(s,a)$.  The mean field is
\begin{equation}
  \bg(\bq)=\bD_{\pi_b}\{\bq-\br-\gamma\bP\mathcal M\bq\}.
  \label{eq:q-mean-field}
\end{equation}
The map $\bq \mapsto \br+\gamma\bP\mathcal M\bq$ is the Bellman optimality operator.  Its
nonlinearity comes from the maximization in $\mathcal M$.  
The root $\bq^\star\in\R^{d_Q}$ is the vectorized optimal action-value
function.  For a prespecified nonzero $\bm v\in\R^{d_Q}$ satisfying
$\sigma_{\bm v}^2>0$, we infer
$\vartheta^\star:=\bm v^\top\bq^\star$ from the scalar sequence
$Y_t:=\bm v^\top\bq_t$.  A coordinate vector gives inference for an individual
$Q^\star(s,a)$, while a general $\bm v$ gives a linear combination of optimal
action values.

\begin{assumption}[Shared Q-Learning Conditions]
\label{ass:q-learning}
The behavior policy $\pi_b$ has full support, meaning that
$\pi_b(a\mid s)>0$ for every $(s,a)\in\cS\times\cA$.  The induced
state--action chain $\{(S_t,A_t)\}_{t\ge0}$ is irreducible and aperiodic, and
the step-size $\{\eta_t\}_{t\ge0}$ satisfies
Assumption~\ref{ass:gain}.
\end{assumption}

\begin{proposition}[Tabular Q-learning]
\label{prop:tabular-q}
Under Assumption~\ref{ass:q-learning}, suppose every state $s$ has a unique optimal
action $a^\star(s)$.  Then the iterates in~\eqref{eq:tabular-q-update} satisfy
Assumptions~\ref{ass:local}--\ref{ass:path}.  Theorem~\ref{thm:fclt} applies
with Jacobian $\bG_Q=\bD_{\pi_b}(\bI_{d_Q}-\gamma\bP\bPi^\star)$, where
$(\bPi^\star\bq)(s)=q(s,a^\star(s))$.
Propositions~\ref{prop:series-interval}--\ref{prop:lm-interval} apply to
$\vartheta^\star$ based on $Y_t$.
\end{proposition}

Under the given conditions, classical asynchronous Q-learning theory gives
$\bq_t\to\bq^\star$ almost surely \citep{tsitsiklis1994asynchronous}.  Because
each state has a unique optimal action and the state and action spaces are
finite, the smallest difference between the optimal action value and the
second-largest action value is strictly positive.  Hence, once $\bq_t$ is
sufficiently close to $\bq^\star$, the action attaining
$\max_a q_t(s,a)$ is $a^\star(s)$ in every state.  Throughout this neighborhood,
the maximum reduces to selecting $q(s,a^\star(s))$, so the mean field $\bg$ is
affine.  This local form yields the Jacobian in
Proposition~\ref{prop:tabular-q} and makes the nonlinear remainder vanish.

\subsection{Q-Learning with Linear Function Approximation}
\label{sec:linear-q}

Tabular Q-learning maintains one parameter for every state--action pair.  A
linear approximation instead represents the action value using a bounded
feature vector $\bvarphi(s,a)\in\R^d$ as
$Q_{\btheta}(s,a)=\bvarphi(s,a)^\top\btheta$.  To make the stability requirement
explicit, we consider the projected semi-gradient recursion on a compact convex
set $\cC$:
\begin{equation}
 \btheta_{t+1}=\Pi_{\cC}\!\left[\btheta_t-\eta_t\bvarphi(S_t,A_t)
 \left\{
 \bvarphi(S_t,A_t)^\top\btheta_t-R_t-\gamma\max_{a'}\bvarphi(S_{t+1},a')^\top\btheta_t\right]\right\}.
 \label{eq:linear-q-update}
\end{equation}
Before projection, this is the standard semi-gradient linear Q-learning update
studied under Markov sampling by \citet{chen2022finite}.
The projection is a stability device.  It is eventually inactive if the target $ \btheta^\star$ lies in the interior of $\cC$.
We infer
$\vartheta^\star:=\bm v^\top\btheta^\star$ from
$Y_t:=\bm v^\top\btheta_t$ for a prespecified nonzero $\bm v\in\R^d$
satisfying $\sigma_{\bm v}^2>0$.
Choosing $\bm v=\bvarphi(s,a)$ then gives inference for the value $Q_{\btheta^\star}(s,a)$.

\begin{proposition}[Linear Q-function approximation]
\label{prop:linear-q}
Under Assumption~\ref{ass:q-learning}, suppose:
\begin{enumerate}[leftmargin=*,label=(\roman*)]
  \item the mean field has a root
  $\btheta^\star\in\operatorname{int}(\cC)$;
  \item for every nonzero $\bu\in\R^d$,
  \begin{equation}
    \gamma^2\sum_{s\in\cS}d_{\pi_b}(s)
    \max_{a\in\cA}\{\bvarphi(s,a)^\top\bu\}^2
    <\sum_{s\in\cS}\sum_{a\in\cA}d_{\pi_b}(s)\pi_b(a\mid s)
    \{\bvarphi(s,a)^\top\bu\}^2;
    \label{eq:linear-q-global-condition}
  \end{equation}
  \item the greedy action under $Q_{\btheta^\star}$ is unique at every state.
\end{enumerate}
Then $\btheta_t\to\btheta^\star$ almost surely, the projection is eventually
inactive, and Assumptions~\ref{ass:local}--\ref{ass:path} hold.  Theorem~\ref{thm:fclt}
applies with Jacobian
$\bG_\varphi=\E_{\pi_b}[\bvarphi(S,A)
\{\bvarphi(S,A)-\gamma\bvarphi(S',a^\star(S'))\}^\top]$.
Propositions~\ref{prop:series-interval}--\ref{prop:lm-interval} provide
inference for $\vartheta^\star$ from $Y_t$.
\end{proposition}

Condition~(ii) is the stability condition used to establish convergence of the
linear Q-learning recursion under Markov sampling \citep{chen2022finite}.  It
requires the behavior policy to collect enough information in every feature
direction so that the maximization term in~\eqref{eq:linear-q-update} does not
amplify parameter errors.  The
condition is restrictive, as some
feature--sampling pairs can make off-policy Q-learning diverge, but it is
directly checkable from the features $\{\bvarphi(s,a)\}_{s,a}$ and stationary sampling law $d_{\pi_b}$.  Uniqueness of
the greedy action makes the mean field $\bm g$ linear near $\btheta^\star$, so the
local remainder vanishes and Condition~(iii) holds automatically.

\subsection{Entropy-Regularized Tabular Q-Learning}
\label{sec:soft-q}

Entropy-regularized Q-learning \citep{geist2019regularized} uses the asynchronous update in
Section~\ref{sec:tabular-q}, replacing its hard maximum $\mathcal{M}$ by a smooth
log-sum-exp operator, denoted by $\mathcal{M}_\tau$.
For a temperature $\tau>0$, $\mathcal{M}_\tau$ is defined by
$(\mathcal M_\tau\bq)(s):=\tau\log\!\sum_a\exp\{q(s,a)/\tau\}$.  The update is
obtained from~\eqref{eq:tabular-q-update} by replacing the maximum with
$(\mathcal M_\tau\bq_t)(S_{t+1})$.

Let
$\bq_\tau^\star:=\operatorname{vec}(Q_\tau^\star)$ denote the unique fixed point
of the resulting soft Bellman operator.  For a nonzero
$\bm v\in\R^{d_Q}$ satisfying $\sigma_{\bm v}^2>0$, the inference target and
scalar iterate are
$\vartheta_\tau^\star:=\bm v^\top\bq_\tau^\star$ and
$Y_t:=\bm v^\top\bq_t$, respectively.  To express the Jacobian, define the
softmax selector $\bPi_\tau^\star\in\R^{|\cS|\times d_Q}$ by
\[
  (\bPi_\tau^\star\bm h)(s):=\sum_{a\in\cA}
  \frac{\exp\{Q_\tau^\star(s,a)/\tau\}}
  {\sum_{a'}\exp\{Q_\tau^\star(s,a')/\tau\}}h(s,a).
\]

\begin{proposition}[Entropy-regularized Q-learning]
\label{prop:soft-q}
Under Assumption~\ref{ass:q-learning}, suppose that $\eta_t=\eta(t+1)^{-\kappa}$  for some $\eta>0$ and $\kappa\in(1/2,1)$.
Then Assumptions~\ref{ass:local}--\ref{ass:path} hold.
Consequently, Theorem~\ref{thm:fclt} applies to $\{\bq_t\}$ with
Jacobian $\bG_\tau=\bD_{\pi_b}(\bI_{d_Q}-\gamma\bP\bPi_\tau^\star)$, and
Propositions~\ref{prop:series-interval}--
\ref{prop:lm-interval} provide inference for $\vartheta_\tau^\star$ from $Y_t$.
\end{proposition}

The soft Bellman operator remains a contraction, so standard results give
$\bq_t\to\bq_\tau^\star$ almost surely under
Assumption~\ref{ass:q-learning}.  This example also illustrates a useful
verification principle: once a mean-square convergence rate is available,
smoothness can verify several path conditions at once.  Here,
\citet[Theorem~1(3)]{chen2021lyapunov} gives Condition~(ii) directly.  Since
log-sum-exp is
globally smooth, the remainder defined in Assumption~\ref{ass:local} satisfies
$\norm{\br_g(\bq)}\le C\norm{\bq-\bq_\tau^\star}^2$ locally.  The same
mean-square error bound therefore yields Condition~(iii).  In particular, entropy
regularization removes the unique-optimal-action condition required in
Proposition~\ref{prop:tabular-q}.

\subsection{Markov Generalized Linear Stochastic Gradients}
\label{sec:markov-glm}

Generalized linear models provide a nonlinear stochastic-gradient application
outside reinforcement learning \citep{mccullagh1989generalized}.  Let $(\bz_t,y_t)$ be a geometrically ergodic Markov process, where
$\bz_t\in\R^d$ is the covariate vector and $y_t$ is the scalar response.  The stochastic gradient is
\begin{equation}
  \bH_{\mathrm{glm}}(\btheta,\bz_t,y_t)=
  \{b'(\bz_t^\top\btheta)-y_t\}\bz_t,
  \label{eq:markov-glm-field}
\end{equation}
where $b$ is the cumulant function.  Define the population field
$\bg(\btheta):=\E_\mu\bH_{\mathrm{glm}}(\btheta,\bz,y)$ under the stationary law $\mu$ of the Markov
data, and let $\btheta^\star$ be its root.  For a compact convex set $\cC$, we use the projected recursion
$\btheta_{t+1}=\Pi_{\cC}[\btheta_t-\eta_t\bH_{\mathrm{glm}}(\btheta_t,\bz_t,y_t)]$.  Given a
prespecified nonzero $\bm v\in\R^d$, we infer
$\vartheta^\star:=\bm v^\top\btheta^\star$ from
$Y_t:=\bm v^\top\btheta_t$.

\begin{proposition}[Markov GLM stochastic gradients]
\label{prop:markov-glm}
Suppose $\btheta^\star\in\operatorname{int}(\cC)$ and
$\eta_t=\eta(t+1)^{-\kappa}$ for some $\eta>0$ and
$\kappa\in(1/2,1)$.  Assume the Markov-chain conditions of
Proposition~\ref{prop:v-route} and that, on a neighborhood of $\cC$,
$\bH_{\mathrm{glm}}$ in~\eqref{eq:markov-glm-field} has an extension to
$\R^d$ satisfying~\eqref{eq:h-v-lipschitz}.  Suppose also that the first two derivatives of
$\bH_{\mathrm{glm}}$ with respect to $\btheta$ have analogous weighted
envelopes on a neighborhood of $\cC$.\footnote{That is, for derivative orders
$j=1,2$, their norms are bounded by $C V(\bz,y)^{1/p}$, uniformly over
$\btheta$ in that neighborhood.}  If
$\inf_{\btheta\in\cC}\lambda_{\min}
\{\E_\mu[b''(\bz^\top\btheta)\bz\bz^\top]\}>0$, then
Assumptions~\ref{ass:local}--\ref{ass:path} hold with
$\bG_{\mathrm{glm}}=\E_\mu[b''(\bz^\top\btheta^\star)\bz\bz^\top]$.
Consequently,
Theorem~\ref{thm:fclt} and Propositions~\ref{prop:series-interval}--
\ref{prop:lm-interval} apply to $\vartheta^\star$ based on $Y_t$.
\end{proposition}

The $V$-geometric condition in Proposition~\ref{prop:v-route} follows, for
example, from a standard Foster--Lyapunov drift condition and small-set
minorization under irreducibility and aperiodicity
\citep[Chapter~16]{meyn2012markov}.  The Hessian condition is the usual uniform
strong-convexity condition for a GLM.  Logistic regression with
$b(u)=\log(1+\exp u)$ satisfies the extension condition under the stated
weighted-moment bounds.  Poisson regression with $b(u)=\exp u$ is also covered
on compact $\cC$ when the corresponding weighted exponential moments are
finite and the displayed Hessian lower bound holds.

\subsection{Balanced LoRA and Inference for the Identified Product}
\label{sec:lora}

Low-rank adaptation (LoRA) keeps a pretrained weight matrix $\bW_0$ fixed and trains
only a rank-$r$ update $\bW=\bB\bA$, where
$\bB\in\R^{m\times r}$, $\bA\in\R^{r\times n}$, and $r<\min\{m,n\}$,
so that the adapted layer uses $\bW_0+\bW$ \citep{hu2022lora}.  The factorization is not
unique, as $(\bB,\bA)$ and $(\bB\bQ,\bQ^{-1}\bA)$ represent the same update for every invertible
$\bQ\in\R^{r\times r}$.  Therefore, a confidence interval for an entry of $\bB$ or $\bA$
would depend on an arbitrary factorization.  We instead conduct inference on
the identified fitted update $\bW=\bB\bA$, which is invariant to the choice of
factorization.

\paragraph{Balanced Factor Updates.}
To learn the product from Markov training data, let $\xi_t=(\bz_t,\by_t)$ be the current example, where $\bz_t\in\R^n$ and $\by_t\in\R^m$.  The sample loss and its
gradient with respect to the low-rank product are
\[
 \ell(\bW;\xi_t)=\frac12\norm{\by_t-(\bW_0+\bW)\bz_t}^2,\quad
 \bGamma(\bW,\xi_t):=\{(\bW_0+\bW)\bz_t-\by_t\}\bz_t^\top.
\]
The balanced algorithm has three operations at iteration $t$: compute this product
gradient, update both factors in the same iteration, and re-express the resulting product
with balanced factors.  In formulas,
\begin{equation}
 \begin{aligned}
 \widetilde{\bB}_{t+1}&=\bB_t-\eta_t\bGamma(\bW_t,\xi_t)\bA_t^\top,\\
 \widetilde{\bA}_{t+1}&=\bA_t-\eta_t\bB_t^\top\bGamma(\bW_t,\xi_t),\\
 (\bB_{t+1},\bA_{t+1})&=\mathfrak B(\widetilde{\bB}_{t+1},\widetilde{\bA}_{t+1}).
 \end{aligned}
 \label{eq:balanced-lora-update}
\end{equation}
The balancing map $\mathfrak B$ is explicit.  Write an $r$-term singular-value
decomposition
$\widetilde{\bB}_{t+1}\widetilde{\bA}_{t+1}=\bL\bSigma\bR^\top$, padding
$\bSigma\in\R^{r\times r}$ with zeros and $\bL,\bR$ with orthonormal columns if
the raw product has rank below $r$.  The map returns
$\bB_{t+1}=\bL\bSigma^{1/2}$ and $\bA_{t+1}=\bSigma^{1/2}\bR^\top$.  Hence, it leaves the
network output unchanged while enforcing
$\bB_{t+1}^\top\bB_{t+1}=\bA_{t+1}\bA_{t+1}^\top=\bSigma$.  For rank one, this is simply a
rescaling of the two factors.  Balancing is thus an algorithmic choice of one
representative of the same product, not an identifiability assumption
\citep{castin2026balancedlora}.

\paragraph{Closed Product Recursion.}
The role of balancing is visible algebraically.  Without it, the product update
depends separately on $\bB\bB^\top$ and $\bA^\top\bA$, which vary across factorizations of
the same $\bW$.  If $\bW=\bB\bA$ is balanced, these two matrices are determined by $\bW$:
$\bB\bB^\top=(\bW\bW^\top)^{1/2}$ and
$\bA^\top\bA=(\bW^\top\bW)^{1/2}$.
Define $\bmM_L(\bW):=(\bW\bW^\top)^{1/2}$ and $\bmM_R(\bW):=(\bW^\top\bW)^{1/2}$. Multiplication of the
two raw updates in~\eqref{eq:balanced-lora-update} gives the exact product recursion
\begin{equation}
 \bW_{t+1}=\bW_t-\eta_t\{\bGamma(\bW_t,\xi_t)\bmM_R(\bW_t)+\bmM_L(\bW_t)\bGamma(\bW_t,\xi_t)\}
 +\eta_t^2\bGamma(\bW_t,\xi_t)\bW_t^\top\bGamma(\bW_t,\xi_t).
 \label{eq:lora-product-recursion}
\end{equation}
Although $\bGamma(\bW,\xi)$ is affine in $\bW$, the order-$\eta_t$ field is
nonlinear through the matrix-square-root maps $\bmM_L(\bW)$ and
$\bmM_R(\bW)$.  Its Taylor remainder is controlled by the quadratic occupation
condition stated below.  The last term in~\eqref{eq:lora-product-recursion} is a
separate order-$\eta_t^2$ perturbation caused by updating both
factors.  Proposition~\ref{prop:second-order-perturbation} shows that this term does
not change the averaged-iterate FCLT.  This closed product recursion is the reason
factor nonidentifiability does not prevent inference for $\bW$.

\paragraph{Why the Product Dynamics Are Locally Stable.}
For inference, the expected update must correct small product errors in every
direction that preserves rank $r$.  Let
${\bmM}_{L,\star}=\bmM_L(\bW^\star)$ and
${\bmM}_{R,\star}=\bmM_R(\bW^\star)$.  At $\bW^\star$, the first-order mean
update applied to an allowable perturbation $\bDelta$ is
\begin{equation}
 \mathcal L_\star[\bDelta]
 =\bDelta{\bmM}_{R,\star}+{\bmM}_{L,\star}\bDelta.
 \label{eq:lora-local-operator}
\end{equation}
For every nonzero such $\bDelta$, $ \langle\bDelta,\mathcal L_\star[\bDelta]\rangle_F
 =\norm{\bDelta{\bmM}_{R,\star}^{1/2}}_F^2
 +\norm{{\bmM}_{L,\star}^{1/2}\bDelta}_F^2>0.$
Thus the mean update dampens every small perturbation of the fitted product.
Changes in $\bB$ and $\bA$ that leave $\bB\bA$ unchanged do not appear in
\eqref{eq:lora-local-operator}; factor nonidentifiability therefore creates no
unstable direction for product inference.

\paragraph{Inference for the Fitted Update.}
Fix a matrix $\bC\in\R^{m\times n}$.  We infer the linear summary
$\vartheta_{\bC}^\star:=\langle\bC,\bW^\star\rangle_F$ from the scalar sequence
$Y_t^{\bC}:=\langle\bC,\bW_t\rangle_F$.  The proposition below combines the
local calculation above with the Markov-noise analysis and the order-$\eta_t^2$
term in~\eqref{eq:lora-product-recursion}.  The only trajectory-level condition
retained below is almost-sure convergence of the fitted product.  The proof derives
the localized mean-square error bound used to prove the required accumulated-error
bound.
\begin{proposition}[Balanced LoRA inference]
\label{prop:lora}
Let $\xi_t=(\bz_t,\by_t)$ be an irreducible and aperiodic Markov chain on a finite state
space.  Suppose its stationary law satisfies
$\by=(\bW_0+\bW^\star)\bz+\bvarepsilon$,
$\E_\mu(\bvarepsilon\bz^\top)=\bzero$, and
$\E_\mu(\bz\bz^\top)=\bI_n$, where $\bW^\star$ has rank $r$.  Let
$\eta_t=\eta(t+1)^{-\kappa}$ for some $\eta>0$ and
$\kappa\in(1/2,1)$.  Assume the product iterates from
\eqref{eq:balanced-lora-update} converge almost surely to $\bW^\star$.  Then, for
every fixed $\bC\in\R^{m\times n}$ for which the asymptotic variance of
$Y_t^{\bC}=\langle\bC,\bW_t\rangle_F:=\operatorname{tr}(\bC^\top\bW_t)$ is
positive, Corollary~\ref{cor:smooth-functional} and
Propositions~\ref{prop:series-interval}--
\ref{prop:lm-interval} apply to
$\vartheta_{\bC}^\star=\langle\bC,\bW^\star\rangle_F$ based on $Y_t^{\bC}$.
\end{proposition}

The finite-state assumption can be relaxed to $V$-geometric Markov chains
under conditions analogous to Proposition~\ref{prop:v-route}.
Choosing $\bC=\be_i\be_j^\top$ gives a confidence interval for the $(i,j)$ entry of the
trained low-rank update.  More general $\bC$ covers a prespecified linear summary of a layer.
The inferential target is the fitted product under repeated Markov training data,
rather than either nonidentified factor or a rank-selection decision.

\section{Experiments}
\label{sec:experiments}

We evaluate coverage, interval length, and computational cost in five
experiments corresponding to the applications studied in this paper.
Each experiment uses 250 independent replications, nominal
$1-\alpha=95\%$ intervals, and eight checkpoints.  The single-recursion
methods use inference windows between 5,000 and 100,000 observations; the
online bootstrap is plotted at approximately matched update-function budgets
and uses 500--10,000 observations.  The inferential target is always a prespecified scalar linear
functional of the parameter or value function learned by the algorithm.  Our
primary method is $P_5$ in~\eqref{eq:confidence-interval}, whose moderate
polynomial dimension is chosen to balance interval length and finite-sample
coverage.  We also report $L_2$ and $L_4$ in
\eqref{eq:lm-confidence-interval} as random-scaling comparisons.  The exact target is used
only after an interval has been constructed, to evaluate empirical coverage.
Every implementable interval is computed from a single observed trajectory.

\subsection{Two Baselines}

We compare the random-scaling intervals with two online baselines.  The first
is the online overlapping-batch-means (OBM) estimator introduced for i.i.d.
SGD by \citet{zhu2021online} and analyzed under Markov sampling by
\citet{roy2023online}.  It estimates the long-run variance of the scalar
iterate sequence $Y_1,\ldots,Y_n$.  We call it Markov OBM to emphasize this
analysis.  It does not estimate a transition model.  Set $a_1=1$ and
$a_j=\lfloor2j^\beta\rfloor$ for $j\geq2$, where
$\beta=2/(1-\kappa)$ and $\kappa=0.55$ is the common step-size exponent.  For the unique $j$ such that
$a_j\leq i<a_{j+1}$, let $t_i=a_j$, $l_i=i-t_i+1$, and
$S_i=\sum_{k=t_i}^iY_k$.  With $\bar Y_n=n^{-1}\sum_{i=1}^nY_i$, the estimator
and its normal interval are
\[
  \widehat\Sigma_n^{\mathrm{OBM}}
  =\frac{\sum_{i=1}^n(S_i-l_i\bar Y_n)^2}{\sum_{i=1}^nl_i},
  \qquad
  \bar Y_n\ \mathbin{\pm}\ \Phi^{-1}(1-\alpha/2)
  \sqrt{\widehat\Sigma_n^{\mathrm{OBM}}/n}.
\]
Thus this baseline uses one trajectory but explicitly estimates its long-run
variance from growing overlapping blocks.

The second baseline, the online bootstrap of \citet{ramprasad2021online},
maintains $B=10$
perturbed recursions driven by the same Markov observations.  Bootstrap path
$b$ evolves as
\[
  \widehat\bx_{t+1}^{(b)}=\widehat\bx_t^{(b)}
  -\eta_tW_t^{(b)}\bH(\widehat\bx_t^{(b)},\xi_t),
  \qquad
  W_t^{(b)}\overset{i.i.d.}{\sim}\operatorname{Unif}(1-\sqrt3,1+\sqrt3).
\]
For $\bar Y_n^{(b)}:=n^{-1}\sum_{t=1}^n\bm v^\top\widehat\bx_t^{(b)}$ and
$\bar Y_n^{(\cdot)}:=B^{-1}\sum_{b=1}^B\bar Y_n^{(b)}$, its interval is
\[
  \bar Y_n\ \mathbin{\pm}\ \Phi^{-1}(1-\alpha/2)
  \left\{\frac{1}{B-1}\sum_{b=1}^B
  \bigl(\bar Y_n^{(b)}-\bar Y_n^{(\cdot)}\bigr)^2\right\}^{1/2}.
\]
An infeasible oracle normal interval uses the across-replication standard
deviation and serves only as a diagnostic.  The bootstrap theory covers linear
SA for policy evaluation \citep{ramprasad2021online}, so the nonlinear comparisons are empirical.

\subsection{Q-Learning Applications}

\begin{figure}[th]
  \centering
  \includegraphics[width=\linewidth]{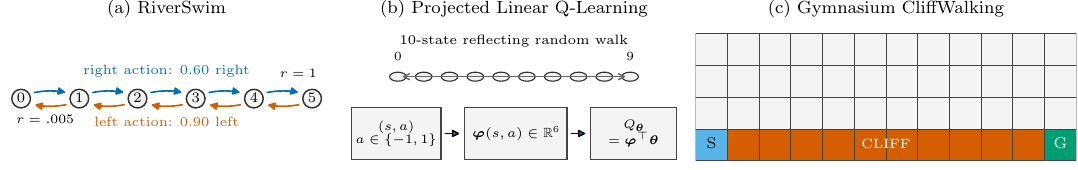}
  \captionof{figure}{Schematics of the three Q-learning experiments. \textbf{(a)} RiverSwim uses a six-state chain with asymmetric left and right dynamics.  \textbf{(b)} The linear
  design maps each state--action pair of a ten-state reflecting chain to six
  features.  \textbf{(c)} CliffWalking uses the $4\times12$ grid, with start (S), cliff cells,
  and goal (G) on the bottom row.}
  \label{fig:q-learning-tasks}
  \vspace{-10pt}
\end{figure}

\paragraph{RiverSwim design.}
RiverSwim has six states $s\in\{1,\ldots,6\}$ and two actions (left and right).  The left
action moves one state to the left with probability $0.9$ and otherwise stays.
The right action moves left, stays, or moves right with probabilities
$(0.05,0.35,0.60)$.  At the boundaries, attempted outward moves remain in place.
Figure~\ref{fig:q-learning-tasks}(a) shows these transitions and the two nonzero
mean rewards: $0.005$ for taking left in state 1 and $1$ for taking right in
state 6.  Centered uniform reward noise has standard deviation $0.2$.  The behavior policy selects left and right with probabilities $(0.3,0.7)$, the discount factor is $\gamma=0.8$, and the step size is
$\eta_t=\min\{1,10(t+1)^{-0.55}\}$.  The target is
$Q^\star(3,\mathrm{right})$.  This design tests asynchronous Q-learning when the
large reward must propagate through one slowly moving trajectory.

\paragraph{Projected linear Q-learning design.}
Let $z_s$ be the ten equally spaced points in $[-1,1]$, and encode the two actions
by $a\in\{-1,1\}$.  The state stays in place with probability $1/2$ and moves one
position left or right with probability $1/4$ each, with reflection at the
endpoints.  Figure~\ref{fig:q-learning-tasks}(b) depicts the reflecting chain and
the feature-based parameterization.  Actions are sampled uniformly.  We use
$\bvarphi(s,a)=(1,z_s,z_s^2,a,az_s,az_s^2)^\top$ and
$\btheta^\star=(0.5,0.3,-0.1,0.4,-0.2,0.15)^\top$.
Define $q^\star(s,a):=\bvarphi(s,a)^\top\btheta^\star$. We construct the mean
reward as
$\bar r(s,a)=q^\star(s,a)-\gamma\E[\max_{a'}q^\star(S',a')\mid S=s]$.
Thus the Bellman solution is exactly representable.  We take $\gamma=0.7$, use
centered uniform reward noise with standard deviation $0.3$, project onto $[-10,10]^6$, and set
$\eta_t=0.7(t+1)^{-0.55}$.  The target is the average of $q^\star(s,a)$ over all
20 state--action pairs.  One can numerically verify that this design satisfies the conditions of
Proposition~\ref{prop:linear-q}.

\paragraph{Entropy-regularized CliffWalking design.}
The third task uses the slippery $4\times12$ CliffWalking transition and reward
table in Gymnasium, the maintained successor to OpenAI Gym
\citep{towers2024gymnasium}.  Figure~\ref{fig:q-learning-tasks}(c) shows its start,
cliff, and goal cells.  The maximum in the Q-learning update is replaced by the
log-sum-exp operator with $\tau=0.05$.  The behavior probabilities for up,
right, down, and left are $(0.35,0.35,0.15,0.15)$.  We set
$\gamma=0.9$ and $\eta_t=\min\{1,20(t+1)^{-0.55}\}$.  The target is the soft
optimal Q-value of taking up in the start state, computed independently by soft
value iteration.  This experiment exercises the smooth Bellman remainder in
Proposition~\ref{prop:soft-q}, rather than relying on the
unique-optimal-action condition.
The three designs therefore satisfy the conditions of
Propositions~\ref{prop:tabular-q}--\ref{prop:soft-q}, respectively.

\FloatBarrier
\begin{figure}[!t]
  \centering
  \setlength{\abovecaptionskip}{5pt}
  \includegraphics[width=0.96\linewidth]{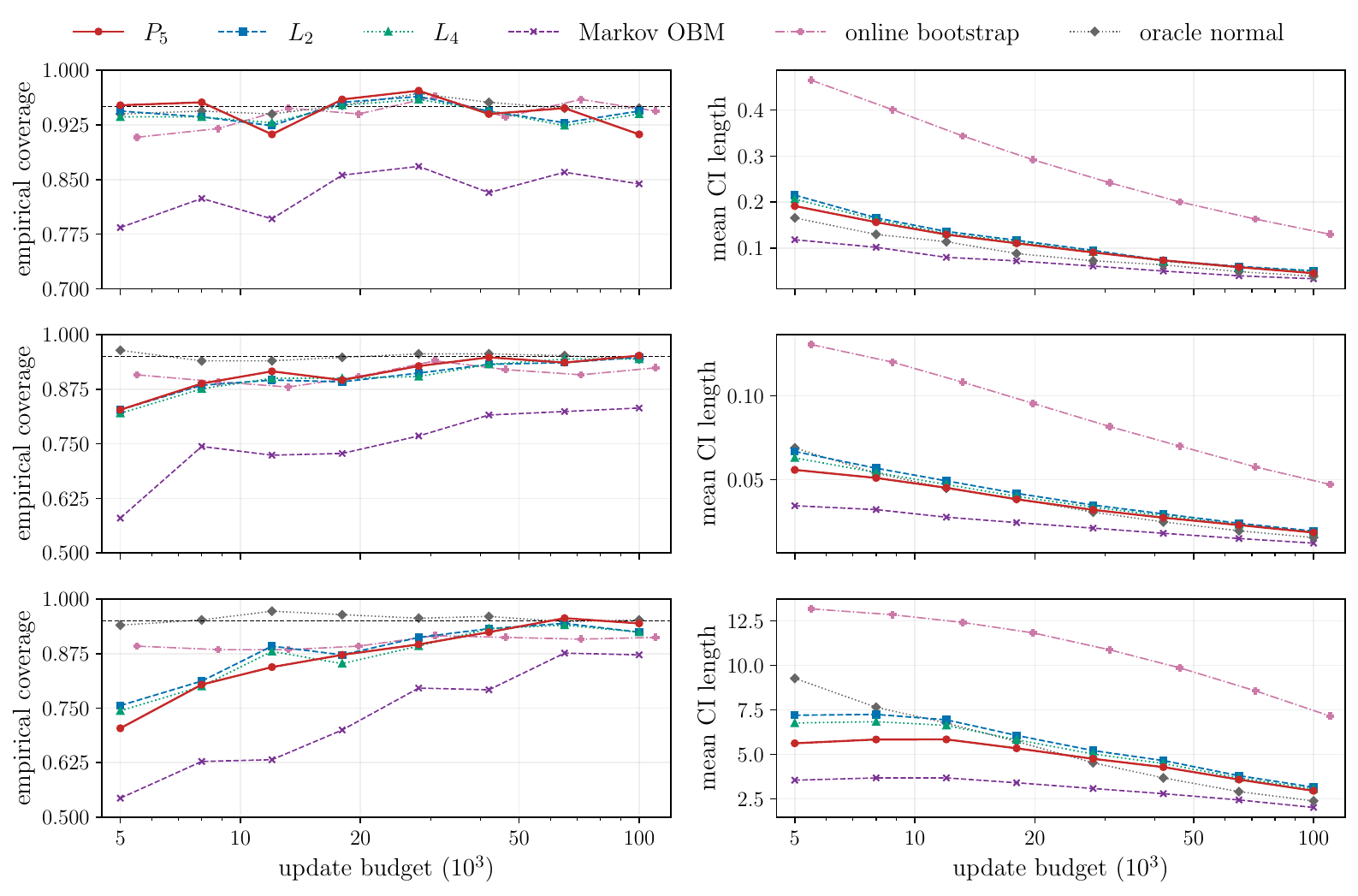}
  \caption{The three Q-learning applications.  Rows correspond to RiverSwim (top),
  projected linear Q-learning (middle), and entropy-regularized CliffWalking (bottom).  Columns show
  empirical coverage (left) and mean interval length (right).}
  \label{fig:q-learning-experiments}
  \vspace{-10pt}
\end{figure}

\paragraph{Results.}
Figure~\ref{fig:q-learning-experiments} shows the intended coverage--length
tradeoff of $P_5$: coverage generally approaches its nominal level while
interval length decreases as the inference window grows.  At $10^5$ updates,
its coverages are $91.2\%$, $95.2\%$, and $94.4\%$ for
RiverSwim, linear Q-learning, and CliffWalking, while its mean intervals are
respectively $10.4\%$, $4.5\%$, and $6.4\%$ shorter than those of $L_2$.
RiverSwim is the exception to the near-nominal coverage pattern because of its
long value-propagation transient.  The $L_2$ coverages are
$94.4\%$, $94.8\%$, and $92.4\%$, while Markov OBM gives
$84.4\%$, $83.2\%$, and $87.2\%$, indicating slower long-run-variance estimation on
these decreasing-step-size paths.  The online-bootstrap coverages are $94.4\%$,
$92.4\%$, and $91.2\%$.  At approximately matched update-function budgets, its intervals are
$2.4$--$2.8$ times as long as those of $P_5$, while requiring ten perturbed
recursions.

\subsection{Nonlinear Markov Logistic SGD}

The covariates are initialized in stationarity and evolve according to
$\bz_{t+1}=0.85\bz_t+\sqrt{1-0.85^2}\,\bvarepsilon_{t+1}$, where
$\bvarepsilon_t\sim N(\bzero,\bI_4)$.  Conditional on $\bz_t$, the response $y_t$ is
Bernoulli with success probability
$\{1+\exp(-\bz_t^\top\btheta^\star)\}^{-1}$, where
$\btheta^\star=(0.6,-0.8,0.5,-0.3)^\top$.  We run logistic SGD projected onto
$[-5,5]^4$ with step size
$\eta_t=\min\{1,1.2(t+1)^{-0.55}\}$ and infer the first coefficient.  This
stationary Gaussian chain is geometrically ergodic with unbounded covariates and
moments of every order, so the conditions of Proposition~\ref{prop:v-route}
can be verified despite the unbounded covariates.  Together with the compact
projection and logistic curvature, these properties satisfy the conditions of
Proposition~\ref{prop:markov-glm}.

The top row of Figure~\ref{fig:nonlinear-applications} shows decreasing interval
lengths and improving coverage.  At $10^5$ updates, $P_5$ covers in $94.4\%$
of replications with an interval $4.3\%$ shorter than $L_2$, which covers in
$93.2\%$.  The oracle, Markov OBM, and online-bootstrap coverages are $94.4\%$,
$84.8\%$, and $88.8\%$.  Thus, in this nonlinear setting, $P_5$ attains near-nominal coverage with a shorter interval than
$L_2$, whereas both online baselines substantially under-cover.

\subsection{Balanced LoRA}

In the last setting, we use balanced rank-one LoRA in matrix regression with three responses and four
covariates.  The target is
$\bW^\star=1.2\bu_\star\bw_\star^\top$, where $\bu_\star$ and $\bw_\star$ normalize
$(0.8,-0.5,0.3)^\top$ and $(0.7,-0.4,0.5,-0.3)^\top$, respectively.  The Markov
covariate $\bz_t$ persists with probability $0.85$.  Otherwise it refreshes
uniformly from $\{\mathord\pm2\be_j:1\leq j\leq4\}$, so its invariant covariance
is $\bI_4$.  Responses satisfy
$\by_t=\bW^\star\bz_t+\bvarepsilon_t$; thus, this experiment takes
$\bW_0=\bzero$.  The independent noise coordinates take
values $\mathord\pm0.2$ equiprobably.  We update both rank-one factors with
$\eta_t=\min\{0.4,0.75(t+1)^{-0.55}\}$ and rescale them to equal Frobenius norm
without changing their product.  The target is
$(\bW^\star)_{11}$.  This construction satisfies the explicit model conditions
of Proposition~\ref{prop:lora}.  Its
inference guarantee is conditional on the assumed almost-sure convergence of
the product iterates.

\FloatBarrier
\begin{figure}[!t]
  \centering
  \setlength{\abovecaptionskip}{5pt}
  \includegraphics[width=0.96\linewidth]{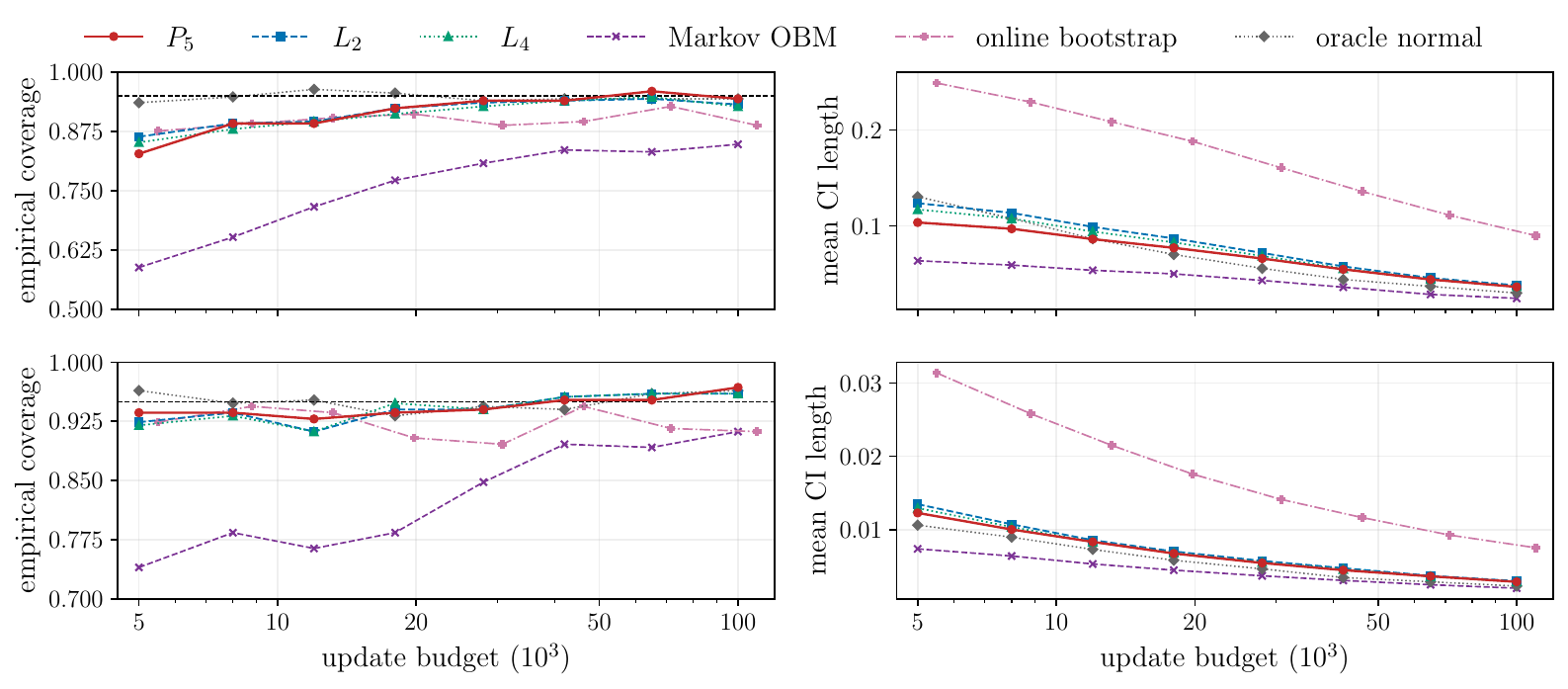}
  \caption{The two nonlinear applications outside Q-learning.  The top row is
  projected logistic SGD with Markov covariates, and the bottom row is balanced
  rank-one LoRA under Markov matrix regression.  Columns show empirical coverage (left)
  and mean interval length (right), using the conventions of
  Figure~\ref{fig:q-learning-experiments}.}
  \label{fig:nonlinear-applications}
  \vspace{-10pt}
\end{figure}

The bottom row of Figure~\ref{fig:nonlinear-applications} shows stable coverage
and decreasing interval length.  At $10^5$ updates, $P_5$ covers in $96.8\%$
with an interval $3.6\%$ shorter than $L_2$, which covers in $96.0\%$.  The
oracle covers in $96.4\%$, while Markov OBM and online bootstrap each cover in
$91.2\%$.

\subsection{Overall Comparison and Computation}

Taken together, the five applications support $P_5$ as a balanced default.  It
shortens the mean interval relative to $L_2$ in every experiment and maintains
coverage close to the nominal level in four of the five, with RiverSwim again
affected by its long transient.  The polynomial-dimension ablation in
Appendix~\ref{app:series-ablation} shows that larger $K$ further shortens the
interval but makes finite-sample coverage less stable.  Thus $K=5$ favors a
coverage--length balance rather than the shortest possible interval.

We finally assess update-function evaluations, storage, and runtime under a
common workload.  Table~\ref{tab:computation} reports a benchmark based on
50,000 tabular-Q observations.  For one scalar target, $P_5$
stores eleven scalar accumulators and adds $18.1\%$ to the
point-estimation runtime, compared with $4.8\%$ for $L_2$ and $7.8\%$ for
Markov OBM.  The online bootstrap instead maintains ten additional parameter
paths, uses 11 times as many evaluations, and takes $2.36$ times the runtime.
Thus $P_5$ costs more than $L_2$ and Markov OBM but remains a constant-memory,
single-recursion procedure and is substantially cheaper than the bootstrap.
Appendix~\ref{app:reproducibility} gives the benchmarking protocol.

\begin{table}[!t]
\centering
\small
\setlength{\tabcolsep}{3pt}
\begin{tabular}{P{0.34\linewidth}rrrr}
\toprule
Method & Seconds & Relative & Function evals & Floats stored \\
\midrule
Point estimate only & 0.165 & $1.00\times$ & 50,000 & 26 \\
$P_5$ & 0.194 & $1.18\times$ & 50,000 & 36 \\
$L_2$ & 0.172 & $1.05\times$ & 50,000 & 29 \\
Markov OBM & 0.177 & $1.08\times$ & 50,000 & 31 \\
Online bootstrap ($B=10$) & 0.388 & $2.36\times$ & 550,000 & 560 \\
\bottomrule
\end{tabular}
\caption{Measured end-to-end runtime for 50,000 tabular-Q observations, together
with exact counts of update-function evaluations and stored values.  The online bootstrap
evaluates the original iterate and 10 perturbed iterates per observation.}
\label{tab:computation}
\vspace{-10pt}
\end{table}

\section{Proof Sketch}
\label{sec:proof-architecture}

This section sketches the proof of Theorem~\ref{thm:fclt}.  We first use the
Poisson equation to isolate the martingale process that determines the
Brownian limit.  We then decompose the difference between this process and the
SA partial sums into four remainders and control each uniformly over
$r\in[0,1]$.  The main difficulty is an endpoint-dependent remainder generated
by the decreasing step sizes.  
The key point is that Lemma~\ref{lem:decreasing-gain} provides the required uniform control and is the key step in passing from a CLT to the FCLT.

\subsection{The Poisson Martingale}

We first isolate the martingale term that determines the Brownian limit.
Recall $\bepsilon_t$ from~\eqref{eq:poisson-innovations}, and define its
target-evaluated counterpart as
$\bepsilon_t^\star:=\cP\bU(\bx^\star,\xi_{t-1})-\bU(\bx^\star,\xi_t)$.
To identify its covariance, define $\bU^\star(\xi):=\bU(\bx^\star,\xi)$ and
\begin{equation}
 \bOmega(\xi):=\cP(\bU^\star\bU^{\star\top})(\xi)
 -(\cP\bU^\star)(\xi)(\cP\bU^\star)(\xi)^\top,
 \qquad \bS:=\int_{\cX}\bOmega(\xi)\,\mu(\mathrm{d}\xi).
 \label{eq:conditional-variance-function}
\end{equation}
The matrix $\bS$ in~\eqref{eq:conditional-variance-function} equals the
long-run covariance in~\eqref{eq:long-run-covariance}.  Both $\bepsilon_t$ and
$\bepsilon_t^\star$ are martingale differences because $\bx_t$ is
$\cF_{t-1}$-measurable.
Let $D([0,1],\R^d)$ denote the space of right-continuous paths with left limits,
equipped with the Skorokhod $J_1$ topology.
The following lemma gives the required martingale approximation.

\begin{lemma}[Poisson martingale approximation]
\label{lem:poisson-martingale}
Under Assumptions~\ref{ass:poisson} and~\ref{ass:path}, the following convergence
holds in $D([0,1],\R^d)$ equipped with the Skorokhod $J_1$ topology:
\begin{equation}
  \left\{\frac1{\sqrt T}\sum_{t=1}^{\lfloor Tr\rfloor}
  \bepsilon_t^\star:0\le r\le1\right\}
  \Rightarrow \bS^{1/2}\bcalW(\cdot),
  \label{eq:root-martingale-fclt}
\end{equation}
and
\begin{equation}
  \max_{1\le n\le T}\frac1{\sqrt T}
  \norm{\sum_{t=1}^n(\bepsilon_t-\bepsilon_t^\star)}
  \xrightarrow{p}0.
  \label{eq:actual-root-martingale}
\end{equation}
Consequently the same FCLT holds with $\bepsilon_t$ in place of $\bepsilon_t^\star$.
Moreover, $\sup_{t\ge1}\E\norm{\bepsilon_t}^p<\infty$.
\end{lemma}
Its proof is in Appendix~\ref{app:poisson}.  Equation~\eqref{eq:root-martingale-fclt} gives the partial-sum limit of
$\bepsilon_t^\star$, while
\eqref{eq:actual-root-martingale} shows that the partial sums of $\bepsilon_t$
and $\bepsilon_t^\star$ are uniformly asymptotically equivalent.  Thus the
serial dependence that remains in the limit is summarized by $\bS$.  We next
connect this martingale limit to the SA recursion by separating the update
noise into $\bepsilon_t$ and several remainders.

\subsection{Martingale--Residual--Coboundary Decomposition}

The Poisson equation gives
$\bH(\bx_t,\xi_t)-\bg(\bx_t)=\bU(\bx_t,\xi_t)-\cP\bU(\bx_t,\xi_t)$.
Adding and subtracting the relevant one-step conditional expectations yields
\begin{align}
  \bH(\bx_t,\xi_t)-\bg(\bx_t)
  =&-\underbrace{\bepsilon_t}_{\text{martingale}}
  +\underbrace{\left\{\frac{\eta_{t+1}}{\eta_t}
  \cP\bU(\bx_{t+1},\xi_t)-\cP\bU(\bx_t,\xi_t)\right\}}_{\bnu_t:\ \text{residual}}
  \nonumber\\
  &+\underbrace{\left\{\cP\bU(\bx_t,\xi_{t-1})
  -\frac{\eta_{t+1}}{\eta_t}\cP\bU(\bx_{t+1},\xi_t)\right\}}_{\bm c_t:\ \text{coboundary}}.
  \label{eq:mrc-decomposition}
\end{align}
This is a step-size-adjusted version of the classical martingale--coboundary
decomposition
\citep{durieu2008comparison,gordin2011functional,liang2010trajectory}.  Its
usefulness comes from the identity
$\eta_t\bm c_t=\eta_t\cP\bU(\bx_t,\xi_{t-1})
-\eta_{t+1}\cP\bU(\bx_{t+1},\xi_t)$, which telescopes in the recursion.

Motivated by this telescoping term, define the auxiliary iterate and its
centered error by
\begin{equation}
  \widetilde\bx_t=\bx_t-\eta_t\cP\bU(\bx_t,\xi_{t-1}),
  \qquad \bDelta_t=\widetilde\bx_t-\bx^\star.
  \label{eq:auxiliary-iterate}
\end{equation}
With this definition, $\widetilde\bx_t$ absorbs the coboundary.  Define also
\begin{equation}
  \br_t=\bg(\bx_t)-\bG\bDelta_t
  =\br_g(\bx_t)+\eta_t\bG\cP\bU(\bx_t,\xi_{t-1}),
  \label{eq:r-definition}
\end{equation}
which collects the local nonlinear remainder and the correction introduced by
the auxiliary iterate.  Substituting~\eqref{eq:mrc-decomposition}
into~\eqref{eq:sa-recursion} gives
\begin{equation}
  \bDelta_{t+1}=(\bI-\eta_t\bG)\bDelta_t
  +\eta_t\{\bepsilon_t-\br_t-\bnu_t\}.
  \label{eq:auxiliary-recursion}
\end{equation}
\begin{lemma}[Negligibility of the Poisson residuals]
\label{lem:poisson-residuals}
Under Assumptions~\ref{ass:local}--\ref{ass:path},
\begin{subequations}
\begin{align}
  \frac1{\sqrt T}\sum_{t=1}^T
  \{\norm{\br_t}+\norm{\bnu_t}\}&\xrightarrow{p}0,
  \label{eq:r-nu-negligible}\\
  \max_{1\le n\le T}\frac1{\sqrt T}
  \norm{\sum_{t=1}^n(\bx_t-\widetilde\bx_t)}&\xrightarrow{p}0,
  \label{eq:auxiliary-negligible}\\
  \frac1{\sqrt T}\max_{1\le t\le T+1}
  \{\norm{\bx_t-\bx^\star}+\norm{\bDelta_t}\}&\xrightarrow{p}0.
  \label{eq:interpolation-control}
\end{align}
\end{subequations}
\end{lemma}
Its proof is in Appendix~\ref{app:decomposed-terms}.
Lemma~\ref{lem:poisson-residuals} shows that the accumulated contributions of
$\br_t$ and $\bnu_t$ are negligible at the $\sqrt T$ scale.  Its other two
conclusions transfer the subsequent analysis of the auxiliary recursion
$\{\widetilde\bx_t\}$ back to the original process $\{\bx_t\}$.  In
particular, \eqref{eq:auxiliary-negligible} makes the partial sums of
$\bx_t-\widetilde\bx_t$ uniformly negligible, while
\eqref{eq:interpolation-control} controls the single-iterate interpolation
terms
$T^{-1/2}(Tr-\lfloor Tr\rfloor)(\bx_{\lfloor Tr\rfloor+1}-\bx^\star)$.
Consequently, the linearly interpolated versions of
$T^{-1/2}\sum_{t=1}^n\bDelta_{t+1}$ and
$T^{-1/2}\sum_{t=1}^n(\bx_t-\bx^\star)$ are uniformly asymptotically
equivalent.  The former is easier to analyze because $\bDelta_t$ satisfies
\eqref{eq:auxiliary-recursion}, with negligible residuals $\br_t$ and
$\bnu_t$.  We study this auxiliary partial-sum process next.

\subsection{Four-Remainder Decomposition}

Having reduced the problem to the auxiliary recursion
in~\eqref{eq:auxiliary-recursion}, we now decompose its partial-sum error into
four remainders.  For $n\ge j$, define
\begin{equation}
  \bm X_j^n=\prod_{i=j}^n(\bI-\eta_i\bG),
  \qquad
  \bm A_j^n=\eta_j\sum_{k=j}^n\bm X_{j+1}^k,
  \label{eq:recursion-matrices}
\end{equation}
where an empty product is $\bI$.  Iterating~\eqref{eq:auxiliary-recursion} and
summing gives, for $n\le T$,
\begin{align}
  \frac1{\sqrt T}\sum_{t=1}^n\bDelta_{t+1}
  -\frac1{\sqrt T}\sum_{j=1}^n\bG^{-1}\bepsilon_j
  ={}&\underbrace{\frac{\bm A_1^n(\bI-\eta_1\bG)\bDelta_1}
  {\sqrt T\,\eta_1}}_{\bm R_{0,T}(n)}
  -\underbrace{\frac1{\sqrt T}\sum_{j=1}^n
  \bm A_j^n(\br_j+\bnu_j)}_{\bm R_{1,T}(n)}
  \nonumber\\
  &+\underbrace{\frac1{\sqrt T}\sum_{j=1}^n
  (\bm A_j^T-\bG^{-1})\bepsilon_j}_{\bm R_{2,T}(n)}
  +\underbrace{\frac1{\sqrt T}\sum_{j=1}^n
  (\bm A_j^n-\bm A_j^T)\bepsilon_j}_{\bm R_{3,T}(n)}.
  \label{eq:four-remainder-decomposition}
\end{align}
The four remainders represent the initialization, the Poisson residuals, the
fixed-endpoint weight approximation, and the change from endpoint $T$ to $n$.
It remains to show that $\max_{1\le n\le T}\norm{\bm R_{k,T}(n)}=o_p(1)$
for $k=0,1,2,3$.

\subsection{Controlling the Remainders: The Decreasing-Step-Size Lemma}


\begin{lemma}[Recursion-matrix properties]
\label{lem:recursion-matrices}
Under Assumptions~\ref{ass:local} and~\ref{ass:gain}, there are constants $c,C>0$ such that
\begin{enumerate}[leftmargin=*,label=(\roman*)]
  \item $\norm{\bm X_j^n}\le C\exp\{-c\sum_{t=j}^n\eta_t\}$;
  \item $\sup_{n\ge j\ge1}\norm{\bm A_j^n}<\infty$ and
  $n^{-1}\sum_{j=1}^n\norm{\bm A_j^n-\bG^{-1}}^2\to0$;
  \item for every $c_0>0$ there is $C_0<\infty$ such that
  \begin{equation}
    \sum_{j=1}^n\eta_{j-1}^2
    \exp\left\{-c_0\sum_{t=j}^n\eta_t\right\}
    \le C_0\eta_n.
    \label{eq:gain-convolution}
  \end{equation}
\end{enumerate}
\end{lemma}
Its proof is in Appendix~\ref{app:recursion-matrices}.  The uniform bound on $\bm A_j^n$ controls $\bm R_{0,T}$ and, together with the
triangle inequality and~\eqref{eq:r-nu-negligible}, $\bm R_{1,T}$.  The
fixed-endpoint partial sums in $\bm R_{2,T}$ form a martingale, so Doob's
maximal inequality, Lemma~\ref{lem:poisson-martingale}, and
Lemma~\ref{lem:recursion-matrices}(ii) give
$\max_{1\le n\le T}\norm{\bm R_{k,T}(n)}\xrightarrow{p}0$ for $k=0,1,2$.

For $\bm R_{3,T}$, an exact rearrangement reduces the problem to
\begin{equation}
  \max_{n\le T}\norm{\bm R_{3,T}(n)}
  \le C\max_{n\le T}
  \frac1{\sqrt T\,\eta_{n+1}}
  \norm{\sum_{j=1}^n\bm X_{j+1}^n\eta_j\bepsilon_j}.
  \label{eq:r3-to-decreasing-gain}
\end{equation}
The following lemma provides the uniform bound needed for the right-hand side.

\begin{lemma}[Uniform Negligibility under Decreasing Step Sizes]
\label{lem:decreasing-gain}
Let $(\bm z_t)$ satisfy
\begin{equation}
  \bm z_0=\bzero,
  \qquad
  \bm z_{t+1}=(\bI-\eta_t\bG)\bm z_t+\eta_t\bm d_t,
  \label{eq:decreasing-gain-recursion}
\end{equation}
where $(\bm d_t,\cF_t)$ is a martingale-difference sequence and
$\sup_t\E\norm{\bm d_t}^p<\infty$ for some $p>2$.  Under Assumptions~\ref{ass:local} and~\ref{ass:gain},
\begin{equation}
  \max_{0\le t\le T}
  \frac{\norm{\bm z_{t+1}}}{\sqrt T\,\eta_{t+1}}
  \xrightarrow{p}0.
  \label{eq:decreasing-gain-negligibility}
\end{equation}
\end{lemma}
Its proof is in Appendix~\ref{app:decreasing-gain}.  To see the correspondence, set $\bm d_0=\bzero$ and
$\bm d_j=\bepsilon_j$ for $j\ge1$.  Iterating
\eqref{eq:decreasing-gain-recursion} gives
\[
  \bm z_{n+1}=\sum_{j=1}^n\bm X_{j+1}^n\eta_j\bepsilon_j.
\]
Thus the right-hand side of~\eqref{eq:r3-to-decreasing-gain} is precisely the
maximum in~\eqref{eq:decreasing-gain-negligibility}, up to the constant $C$.
Lemma~\ref{lem:decreasing-gain} therefore gives
$\max_{n\le T}\norm{\bm R_{3,T}(n)}=o_p(1)$.  Because the coefficients
$\bm X_{j+1}^n$ change with $n$, $\{\bm R_{3,T}(n)\}$ is not a martingale,
making this uniform recursion bound essential.  The lemma requires no symmetry
or diagonalizability of $\bG$; its eigenvalues need only have positive real
parts.

Together with the controls for $\bm R_{0,T}$, $\bm R_{1,T}$, and
$\bm R_{2,T}$, \eqref{eq:four-remainder-decomposition} now gives the desired
representation for the auxiliary partial sums.  Equations
\eqref{eq:auxiliary-negligible} and~\eqref{eq:interpolation-control} then
transfer this representation from $\{\bDelta_t\}$ to the original iterate
process and yield
\begin{equation}
  \max_{1\le n\le T}\frac1{\sqrt T}
  \norm{\sum_{t=1}^n(\bx_t-\bx^\star)
  -\bG^{-1}\sum_{t=1}^n\bepsilon_t}
  \xrightarrow{p}0.
  \label{eq:uniform-asymptotic-linearity}
\end{equation}
This is the required uniform asymptotic linear representation.  By
Lemma~\ref{lem:poisson-martingale}, its leading martingale path converges to
$\bG^{-1}\bS^{1/2}\bcalW(\cdot)$.  Because this limit is continuous and
\eqref{eq:interpolation-control} makes the linear interpolation negligible, the
convergence holds in $C([0,1],\R^d)$ under the uniform norm.  This proves
Theorem~\ref{thm:fclt}.

\section{Conclusion and Discussion}
\label{sec:conclusion}

This paper develops online inference for nonlinear stochastic approximation
from a single Markov trajectory.  Its main result is an FCLT for the partial
sums of the iterate errors over all fractions of the observed trajectory.  This
limit enables self-normalization because the estimation error and its random
scale converge jointly with the same unknown asymptotic standard deviation,
which cancels in their ratio.  The polynomial-series and bridge methods thereby
give constant-memory confidence intervals without long-run covariance
estimation or auxiliary recursions.  The proof combines a Poisson-equation
decomposition with a uniform lemma for decreasing step sizes.

An important next step is to establish a weak-convergence rate for the
partial-sum process under Markov noise.  The present FCLT gives asymptotic
coverage but does not quantify the distance between the finite-sample process
and its Brownian limit.  A quantitative bound could yield
coverage-error guarantees and clarify how initialization, dependence, step
size, and warm-up affect the accuracy of the intervals.

A second direction is to study optimal self-normalizing functionals.  The
experiments show a clear trade-off: increasing $K$ can shorten the
polynomial-series interval but may make finite-sample coverage less stable.
A useful theory would optimize expected length subject to coverage accuracy and
computational cost, and guide data-dependent choices that preserve validity.

A third direction is to extend path-level inference to more complex methods
such as AdamW and Muon.  Existing CLTs for averaged Adam-type recursions
\citep{an2026saadam} and convergence results for Muon updates
\citep{kim2026muon} provide starting points, but FCLTs and single-trajectory
self-normalized inference remain open.

\clearpage
\appendix
\section*{Appendix Contents}
\printappendixcontents
\appendixsection{Poisson Equations and Martingale Approximation}
\label{app:poisson}

This appendix proves Proposition~\ref{prop:v-route} and Lemma~\ref{lem:poisson-martingale}.  Throughout the appendices, constants denoted by $C$ may change from line to line.

\begin{proof}[Proof of Proposition~\ref{prop:v-route}]
For fixed $\bx\in\R^d$, set $\bm h_{\bx}=\bH(\bx,\cdot)-\bg(\bx)$ and denote
$\mu(V^{1/p}):=\int V(\xi)^{1/p}\,\mu(\mathrm{d}\xi)$.  Jensen's inequality
gives $\mu(V^{1/p})\le\{\mu(V)\}^{1/p}<\infty$.  By~\eqref{eq:h-v-lipschitz},
$\norm{\bH(\bx,\cdot)}_{V^{1/p}}\le C(1+\norm{\bx-\bx^\star})$.  Since
$\bg(\bx)=\int\bH(\bx,\xi)\,\mu(\mathrm{d}\xi)$, it follows that
$\norm{\bg(\bx)}\le\mu(V^{1/p})\norm{\bH(\bx,\cdot)}_{V^{1/p}}$.
Because $V^{1/p}\ge1$, the triangle inequality yields
\[
  \norm{\bm h_{\bx}}_{V^{1/p}}
  \le C\{1+\mu(V^{1/p})\}(1+\norm{\bx-\bx^\star})
  \le C(1+\norm{\bx-\bx^\star}),
\]
where the last $C$ absorbs the fixed finite factor $1+\mu(V^{1/p})$.
Moreover, $\mu(\bm h_{\bx})=\bzero$, so centering the transition law by $\mu$
does not change $\cP^k\bm h_{\bx}(\xi)$.  The definition of $V^{1/p}$-geometric bound in
Definition~\ref{def:v-geometric} implies that for any $ k\ge0$,
\begin{equation}
\begin{aligned}
  \norm{\cP^k\bm h_{\bx}(\xi)}
  &\le \norm{\bm h_{\bx}}_{V^{1/p}}
    \norm{\cP^k(\xi,\cdot)-\mu}_{V^{1/p}}
\le C(1+\norm{\bx-\bx^\star})\rho^k V(\xi)^{1/p}.
\end{aligned}
  \label{eq:appendix-resolvent-bound}
\end{equation}
After division by $V(\xi)^{1/p}$ and taking the supremum over $\xi$,
\eqref{eq:appendix-resolvent-bound} shows that
$\sum_{k\ge0}\cP^k\bm h_{\bx}$ converges absolutely in the
$V^{1/p}$-weighted norm.  Its sum is~\eqref{eq:poisson-resolvent}, and hence
\[
  \bU(\bx,\cdot)-\cP\bU(\bx,\cdot)
  =\sum_{k\ge0}\cP^k\bm h_{\bx}
   -\sum_{k\ge1}\cP^k\bm h_{\bx}
  =\bm h_{\bx},
\]
which is~\eqref{eq:poisson-family}.  Each summand is $\mu$-centered, so
$\int_{\cX}\bU(\bx,\xi)\,\mu(\mathrm{d}\xi)=\bzero$.

For $\bx,\by\in\R^d$, define $  \bm h_{\bx,\by}
=\{\bH(\bx,\cdot)-\bH(\by,\cdot)\}
-\{\bg(\bx)-\bg(\by)\}.$
The second part of~\eqref{eq:h-v-lipschitz} and Jensen's inequality under $\mu$ give
$\norm{\bm h_{\bx,\by}}_{V^{1/p}}\le C\norm{\bx-\by}$.  Repeating~\eqref{eq:appendix-resolvent-bound} yields
\begin{equation}
  \norm{\bU(\bx,\cdot)-\bU(\by,\cdot)}_{V^{1/p}}
  \le \frac{C}{1-\rho}\norm{\bx-\by}.
  \label{eq:appendix-u-weighted-lipschitz}
\end{equation}
The weighted bound in~\eqref{eq:appendix-u-weighted-lipschitz} gives
\[
 \left\{\cP\norm{\bU(\bx,\cdot)-\bU(\by,\cdot)}^p(\xi)\right\}^{1/p}
 \le C\{(\cP V)(\xi)\}^{1/p}\norm{\bx-\by},
\]
which proves~\eqref{eq:u-weighted-lipschitz} with
$L(\xi)=C\{1+(\cP V)(\xi)^{1/p}\}$.  The weighted envelope of
$\bU(\bx^\star,\cdot)$ gives its required $p$th moment.  Finally,
$\E(\cP V)(\xi_t)=\E V(\xi_{t+1})$.  Hence, the assumed uniform $V$ moment gives
the $p$th moments of $L(\xi_t)$ and $\bU(\bx^\star,\xi_t)$ required
by~\eqref{eq:poisson-moments}.  
\end{proof}

\begin{proof}[Proof of Lemma~\ref{lem:poisson-martingale}]
First consider the target process $\bepsilon_t^\star$.  Since
$\bU(\bx^\star,\xi_t)$ has conditional expectation
$\cP\bU(\bx^\star,\xi_{t-1})$ given $\cF_{t-1}$,
$(\bepsilon_t^\star,\cF_t)$ is a martingale-difference sequence.  Its conditional covariance is
\begin{equation}
  \E(\bepsilon_t^\star\bepsilon_t^{\star\top}\mid\cF_{t-1})
  =\bOmega(\xi_{t-1}),
  \label{eq:appendix-root-covariance}
\end{equation}
with $\bOmega$ defined in~\eqref{eq:conditional-variance-function}.  The $p>2$ moment in~\eqref{eq:poisson-moments} makes $\bOmega$ integrable.  The Markov-chain ergodic theorem gives
\[
  \frac1n\sum_{t=1}^n\bOmega(\xi_{t-1})\longrightarrow\bS
  \quad\text{almost surely}.
\]
The convergence is uniform after normalization by $T$: if
$\bm Q_n=\sum_{t=1}^n\{\bOmega(\xi_{t-1})-\bS\}$, then
$\bm Q_n/n\to0$ almost surely, and the usual finite-prefix argument gives
\begin{equation}
  \max_{n\le T}\frac{\norm{\bm Q_n}}T\longrightarrow0
  \quad\text{almost surely}.
  \label{eq:appendix-uniform-quadratic-variation}
\end{equation}

We next verify that $\bS$ equals the long-run covariance in
\eqref{eq:long-run-covariance}.  Since the chain is positive Harris recurrent
and aperiodic, the total-variation distance
$\norm{\mathcal L(\xi_t)-\mu}_{\mathrm{TV}}\to0$
\citep[Theorem~13.0.1]{meyn2012markov}.  Applying this convergence to bounded
truncations of $\norm{\bU(\bx^\star,\cdot)}^p$ and then using monotone
convergence shows from~\eqref{eq:poisson-moments} that
$\int\norm{\bU(\bx^\star,\xi)}^p\,\mu(\mathrm{d}\xi)<\infty$.
Under stationarity, put $\bm M_T=\sum_{t=1}^T\bepsilon_t^\star$ and
$\bm b_t^\star=\cP\bU(\bx^\star,\xi_t)$.  The Poisson equation gives
$\sum_{t=1}^T\bH(\bx^\star,\xi_t)=-\bm M_T+\bm b_0^\star-\bm b_T^\star$.
Stationarity and conditional Jensen's inequality imply
$\sup_t\E_\mu\norm{\bm b_t^\star}^2<\infty$, and hence
$T^{-1/2}(\bm b_0^\star-\bm b_T^\star)\to\bzero$ in $L^2$.
Martingale orthogonality and~\eqref{eq:appendix-root-covariance} give
$T^{-1}\E_\mu(\bm M_T\bm M_T^\top)=\int\bOmega\,\mathrm{d}\mu$.
Thus the normalized covariance of the update sum converges to $\bS$.

Conditional Jensen's inequality and~\eqref{eq:poisson-moments} give
$\sup_t\E\norm{\cP\bU(\bx^\star,\xi_t)}^p<\infty$ and hence
$\sup_t\E\norm{\bepsilon_t^\star}^p<\infty$.  For every $\delta>0$,
conditional Markov inequality and this uniform $p$th moment give
\begin{align*}
 &\E\left[\frac1T\sum_{t=1}^T
 \E\{\norm{\bepsilon_t^\star}^2
 \mathbf 1(\norm{\bepsilon_t^\star}>\delta\sqrt T)
 \mid\cF_{t-1}\}\right]
 \le
 \frac{1}{\delta^{p-2}T^{p/2}}
 \sum_{t=1}^T\E\norm{\bepsilon_t^\star}^p
 \le C T^{1-p/2}\longrightarrow0.
\end{align*}
Thus the conditional Lindeberg condition holds in probability.  Equations~\eqref{eq:appendix-root-covariance}--\eqref{eq:appendix-uniform-quadratic-variation} and the multidimensional martingale invariance principle \citep{hall2014martingale,whitt2007proofs} prove~\eqref{eq:root-martingale-fclt}.

It remains to replace $\bepsilon_t^\star$ by $\bepsilon_t$.  Put
$\bm d_t=\bepsilon_t-\bepsilon_t^\star$.  This is again a martingale difference.
Conditional Jensen's inequality,
\eqref{eq:u-weighted-lipschitz}, and~\eqref{eq:path-moments} give
\[
 \E\norm{\bm d_t}^p
 \le C\E\{L(\xi_{t-1})^p\norm{\bx_t-\bx^\star}^p\}\le C.
\]
Together with the preceding bound for $\bepsilon_t^\star$, this proves
$\sup_t\E\norm{\bepsilon_t}^p<\infty$.  This moment bound will be used in the
proof of Lemma~\ref{lem:four-remainders}.

For the uniform replacement, the same conditional argument implies
\[
  \E\norm{\bm d_t}^2
  \le C\E\{L(\xi_{t-1})^2\norm{\bx_t-\bx^\star}^2\}.
\]
We claim that the average on the right tends to zero.  For $M>0$, split according to $L(\xi_{t-1})\le M$.  The first part is bounded by
$M^2\E\norm{\bx_t-\bx^\star}^2$.  For the second part, H\"older's and Markov's inequalities, together with~\eqref{eq:poisson-moments} and~\eqref{eq:path-moments}, give
\[
 \E\!\left[L^2\norm{\bm e_t}^2\mathbf 1\{L>M\}\right]
 \le \{\E(L^p\norm{\bm e_t}^p)\}^{2/p}
      \Pp(L>M)^{(p-2)/p}
 \le C M^{-(p-2)},
\]
where $L=L(\xi_{t-1})$.  Hence,
\[
 \limsup_{T\to\infty}\frac1T\sum_{t=1}^T
 \E\{L(\xi_{t-1})^2\norm{\bx_t-\bx^\star}^2\}
 \le C M^{-(p-2)}.
\]
 Letting $M\to\infty$ proves the claim.  Doob's $L^2$ inequality now yields
\[
 \E\max_{n\le T}\frac1T
 \norm{\sum_{t=1}^n\bm d_t}^2
 \le \frac{4}{T}\sum_{t=1}^T\E\norm{\bm d_t}^2\longrightarrow0.
\]
This proves~\eqref{eq:actual-root-martingale}.  Slutsky's theorem completes the proof.
\end{proof}

\appendixsection{The Decomposed Poisson Terms}
\label{app:decomposed-terms}

We repeatedly use $T^{-1/2}\sum_{t=1}^T\eta_t\to0$.  Indeed, for any fixed
$N<T$, Cauchy--Schwarz gives
$T^{-1/2}\sum_{t=1}^T\eta_t
\le T^{-1/2}\sum_{t=1}^N\eta_t
+(\sum_{t=N+1}^{\infty}\eta_t^2)^{1/2}$.
Letting first $T\to\infty$ and then $N\to\infty$ proves the claim because
$\sum_t\eta_t^2<\infty$.

\begin{proof}[Proof of Lemma~\ref{lem:poisson-residuals}]
Recall from~\eqref{eq:r-definition} that $  \br_t=\br_g(\bx_t)+\eta_t\bG\cP\bU(\bx_t,\xi_{t-1}).$
Conditional Jensen's inequality and~\eqref{eq:u-weighted-lipschitz} give
\[
 \norm{\cP\bU(\bx_t,\xi_{t-1})}
 \le \norm{\cP\bU(\bx^\star,\xi_{t-1})}
 +C L(\xi_{t-1})\norm{\bm e_t}.
\]
The first term has a uniform $p$th moment by conditional Jensen's inequality
and~\eqref{eq:poisson-moments}, while the second term is controlled
by~\eqref{eq:path-moments}.  Therefore,
\begin{equation}
 \sup_t\E\norm{\cP\bU(\bx_t,\xi_{t-1})}^p<\infty.
 \label{eq:appendix-pu-path-moment}
\end{equation}
The same bound with $t$ replaced by $t+1$ will be used below.
Assumption~\ref{ass:path} controls the sum of the norms of the first term in
$\br_t$, and~\eqref{eq:appendix-pu-path-moment} controls the second.  Therefore, by Markov's inequality and~\eqref{eq:gain-conditions},
\begin{equation}
  \frac1{\sqrt T}\sum_{t=1}^T
  \eta_t\norm{\cP\bU(\bx_t,\xi_{t-1})}
  \xrightarrow{p}0.
  \label{eq:appendix-pu-sum}
\end{equation}

For $\bnu_t$, add and subtract $\cP\bU(\bx_{t+1},\xi_t)$:
\begin{align}
 \norm{\bnu_t}
 &\le C L(\xi_t)\norm{\bx_{t+1}-\bx_t}
 +\left|\frac{\eta_{t+1}-\eta_t}{\eta_t}\right|
 \norm{\cP\bU(\bx_{t+1},\xi_t)}.
 \label{eq:appendix-nu-bound}
\end{align}
The recursion gives
$\norm{\bx_{t+1}-\bx_t}=\eta_t\norm{\bH(\bx_t,\xi_t)}$.  Cauchy--Schwarz and the uniform $p>2$ moments of $L$ and $\bH$ show that the expectation of the first term in~\eqref{eq:appendix-nu-bound} is at most $C\eta_t$.  Equation~\eqref{eq:appendix-pu-path-moment} bounds the expectation of the last factor in the second term.  Finally, $ \left|\frac{\eta_{t+1}-\eta_t}{\eta_t}\right|=o(\eta_t)$ by~\eqref{eq:gain-conditions}.  Consequently
\begin{equation}
  \frac1{\sqrt T}\sum_{t=1}^T\E\norm{\bnu_t}
  \le \frac{C}{\sqrt T}\sum_{t=1}^T\eta_t+o(1)\longrightarrow0.
  \label{eq:appendix-nu-sum}
\end{equation}
Combining~\eqref{eq:integrated-nonlinearity},
\eqref{eq:appendix-pu-sum}, and~\eqref{eq:appendix-nu-sum} proves~\eqref{eq:r-nu-negligible}.

By~\eqref{eq:auxiliary-iterate},
\[
 \max_{n\le T}\frac1{\sqrt T}
 \norm{\sum_{t=1}^n(\bx_t-\widetilde\bx_t)}
 \le \frac1{\sqrt T}\sum_{t=1}^T
 \eta_t\norm{\cP\bU(\bx_t,\xi_{t-1})},
\]
so~\eqref{eq:auxiliary-negligible} follows from~\eqref{eq:appendix-pu-sum}.

For every $\delta>0$, the union bound and~\eqref{eq:path-moments} yield
\[
 \Pp\left(\max_{t\le T+1}\norm{\bx_t-\bx^\star}>\delta\sqrt T\right)
 \le \sum_{t=1}^{T+1}
 \frac{\E\norm{\bx_t-\bx^\star}^p}{\delta^pT^{p/2}}
 \le \frac{C(T+1)}{\delta^p T^{p/2}}
 \longrightarrow0
\]
because $p>2$.  Moreover,
$\bDelta_t=(\bx_t-\bx^\star)-\eta_t\cP\bU(\bx_t,\xi_{t-1})$.  The same union-bound argument and~\eqref{eq:appendix-pu-path-moment} control the second term.  This proves~\eqref{eq:interpolation-control}.
\end{proof}

\appendixsection{Recursion-Matrix Properties}
\label{app:recursion-matrices}

\begin{proof}[Proof of Lemma~\ref{lem:recursion-matrices}]
Because every eigenvalue of $\bG$ has positive real part, the Lyapunov theorem
gives $\bQ\succ0$ satisfying $\bG^\top\bQ+\bQ\bG=\bI$.  Set
$\norm{\bx}_{\bQ}=(\bx^\top\bQ\bx)^{1/2}$; this norm is equivalent to the
Euclidean norm.  Expanding in this norm gives
$\norm{(\bI-\eta\bG)\bx}_{\bQ}^2
\le(1-c_1\eta+c_2\eta^2)\norm{\bx}_{\bQ}^2$ for some $c_1,c_2>0$.
Consequently, there are constants $a>0$ and $\bar\eta>0$ such that
\begin{equation}
  \norm{\bI-\eta\bG}_{\bQ}\le1-a\eta
  \le e^{-a\eta},\qquad 0<\eta\le\bar\eta.
  \label{eq:appendix-stable-norm}
\end{equation}
Absorb the finitely many earlier $\eta_t$ into a constant and use equivalence of norms.  This proves part (i).

We next quantify the slow variation of the step sizes.  Put
$\delta_j=\sup_{k\ge j+1}(\eta_{k-1}-\eta_k)/\eta_{k-1}^2$; then
$\delta_j\to0$.  For any $b>0$, \eqref{eq:gain-conditions} implies, for all
$k\ge j\ge j_0(b)$,
\begin{equation}
  \frac{\eta_j}{\eta_k}
  =\prod_{i=j+1}^k\frac{\eta_{i-1}}{\eta_i}
  \le C_b\exp\left\{b\sum_{i=j+1}^k\eta_i\right\}.
  \label{eq:appendix-gain-ratio}
\end{equation}
Indeed, the logarithm of each ratio is at most
$2(\eta_{i-1}-\eta_i)/\eta_{i-1}=o(\eta_{i-1})$, uniformly once $j$ is large.

Let $\bm Y_k=\bm X_{j+1}^k$, so $\bm Y_j=\bI$.  Since
$\bm Y_k-\bm Y_{k+1}=\eta_{k+1}\bG\bm Y_k$ and
$\bG\bm A_j^n=\sum_{k=j}^n\eta_j(\bm Y_k-\bm Y_{k+1})/\eta_{k+1}$,
summation by parts gives
\begin{equation}
 \bG\bm A_j^n=\frac{\eta_j}{\eta_{j+1}}\bI-
 \frac{\eta_j}{\eta_{n+1}}\bm Y_{n+1}+
 \sum_{k=j+1}^n\left(\frac{\eta_j}{\eta_{k+1}}-\frac{\eta_j}{\eta_k}\right)\bm Y_k.
 \label{eq:appendix-A-identity}
\end{equation}
The step-size condition gives
$|\eta_j/\eta_{j+1}-1|\le C\delta_j\eta_j$.
Taking $b=a/2$ in~\eqref{eq:appendix-gain-ratio}, writing
$j_0=j_0(a/2)$, and using
\eqref{eq:appendix-stable-norm} bound the second term by
$C\exp\{-a\sum_{i=j+1}^{n+1}\eta_i/2\}$.  For the third term in~\eqref{eq:appendix-A-identity} (which is a sum), $ 0\le \frac{\eta_j}{\eta_{k+1}}-\frac{\eta_j}{\eta_k}
\le C\delta_j\eta_k\frac{\eta_j}{\eta_{k+1}},$
and hence the norm of the sum is at most $ C\delta_j\sum_{k=j+1}^n
\eta_k\exp\left\{-\frac a2\sum_{i=j+1}^{k}\eta_i\right\}
\le C\delta_j.$
Thus
\begin{equation}
 \norm{\bm A_j^n-\bG^{-1}}
 \le C\left[\delta_j+
 \exp\left\{-\frac a2\sum_{i=j+1}^{n+1}\eta_i\right\}\right]
 \label{eq:appendix-A-approximation}
\end{equation}
for $n\ge j\ge j_0$, and the same calculation gives a uniform bound on $\bm A_j^n$.

To prove part~(ii), first use $(u+v)^2\le2u^2+2v^2$ in
\eqref{eq:appendix-A-approximation}.  The finitely many indices $j<j_0$
contribute at most $C/n$ by the uniform bound on $\bm A_j^n$.  Hence,
\begin{align*}
 \frac1n\sum_{j=1}^n\norm{\bm A_j^n-\bG^{-1}}^2
 &\le \frac Cn+\frac Cn\sum_{j=j_0}^n
 \left[\delta_j^2+
 e^{-a\sum_{i=j+1}^{n+1}\eta_i}\right]\\
 &\le \frac Cn+\frac Cn\sum_{j=j_0}^n\delta_j^2
 +\frac Cn\sum_{j=j_0}^n e^{-a(n-j)\eta_n}\\
 &\le \frac Cn+\frac Cn\sum_{j=j_0}^n\delta_j^2
 +\frac{C}{n(1-e^{-a\eta_n})}
 \le \frac Cn+\frac Cn\sum_{j=j_0}^n\delta_j^2
 +\frac{C}{n\eta_n}\longrightarrow0.
\end{align*}
The second line uses $\sum_{i=j+1}^{n+1}\eta_i\ge(n-j)\eta_n$, and the
last line uses $1-e^{-a\eta_n}\ge c\eta_n$ for large $n$.  The average of
$\delta_j^2$ tends to zero because $\delta_j\to0$, while
$n\eta_n\to\infty$.  This proves part~(ii).

For part (iii), enlarge $j_0$ if needed and take $b=c_0/2$
in~\eqref{eq:appendix-gain-ratio} to obtain, for $j\ge j_0$, $  \eta_{j-1}\le C\eta_n
\exp\left\{\frac{c_0}{2}\sum_{t=j}^n\eta_t\right\}.$
Therefore,
\[
 \sum_{j=j_0}^n\eta_{j-1}^2e^{-c_0\sum_{t=j}^n\eta_t}
 \le C\eta_n\sum_{j=j_0}^n\eta_{j-1}e^{-(c_0/2)\sum_{t=j}^n\eta_t}
 \le C\eta_n.
\]
The last sum is a Riemann-sum bound for an exponential integral and is uniformly bounded.  For the finitely many $j<j_0$, part (i), $t\eta_t\to\infty$, and a larger constant give the same bound.  This proves~\eqref{eq:gain-convolution}.
\end{proof}

\appendixsection{Uniform Negligibility under Decreasing Step Sizes}
\label{app:decreasing-gain}

\begin{lemma}[Diagonalizable case]
\label{lem:diagonal-decreasing-gain}
Under the conditions of Lemma~\ref{lem:decreasing-gain}, if $\bG$ is diagonalizable over $\mathbb C$, then~\eqref{eq:decreasing-gain-negligibility} holds.
\end{lemma}

\begin{proof}[Proof of Lemma~\ref{lem:diagonal-decreasing-gain}]
Let $\bm X_j^n$ be as in~\eqref{eq:recursion-matrices}.  Iteration of~\eqref{eq:decreasing-gain-recursion} gives $  \bm z_{t+1}=\sum_{j=0}^t\bm X_{j+1}^t\eta_j\bm d_j.$
Martingale orthogonality and Lemma~\ref{lem:recursion-matrices}(i),(iii) give
\begin{equation}
  \E\norm{\bm z_{t+1}}^2
  \le C\sum_{j=0}^t\eta_j^2
  e^{-c\sum_{i=j+1}^t\eta_i}
  \le C\eta_{t+1},
  \qquad
  \sup_{t\ge0}\frac{\E\norm{\bm z_{t+1}}^2}{\eta_{t+1}}<\infty.
  \label{eq:appendix-z-second-moment}
\end{equation}

Choose $t_0$ so that all factors after $t_0$ are nonzero.  The stable-product
bound in Lemma~\ref{lem:recursion-matrices}(i) and the step-size-ratio bound
in~\eqref{eq:appendix-gain-ratio} imply
\[
 \sup_{t_0\le t\le T}
 \frac{\norm{\bm X_{t_0}^t\bm z_{t_0}}}{\sqrt T\eta_{t+1}}
 \le \frac{C\norm{\bm z_{t_0}}}{\sqrt T\eta_{t_0}}
 \xrightarrow{p}0.
\]
Thus the fixed prefix is negligible and will be suppressed below.

Fix $m\ge2$ and partition $\{0,\ldots,T\}$ into half-open blocks, separating
each starting value from the martingale increments within the block.  Set
$h_k=\lfloor k(T+1)/m\rfloor$, $0\le k\le m$, so that $h_m=T+1$.  The step-size-ratio argument used in~\eqref{eq:appendix-gain-ratio} and Lemma~\ref{lem:recursion-matrices}(i) give
\begin{equation}
  \sup_{h_k\le t<h_{k+1}}
  \frac{\eta_{h_k}}{\eta_{t+1}}\norm{\bm X_{h_k}^t}\le C.
  \label{eq:appendix-block-propagation}
\end{equation}
Define
\[
 \mathcal A_{T,m}=\bigcap_{k=0}^{m-1}
 \left\{\frac{C\norm{\bm z_{h_k}}}{\sqrt T\eta_{h_k}}<\varepsilon\right\}.
\]
By~\eqref{eq:appendix-z-second-moment}, the union bound, and monotonicity of $\eta_t$,
\begin{equation}
  \Pp(\mathcal A_{T,m}^c)
  \le \frac{Cm}{\varepsilon^2T\eta_T}\longrightarrow0
  \quad\text{for fixed }m.
  \label{eq:appendix-block-starts}
\end{equation}

We next control the martingale differences inside one block.  Diagonalize
$\bG=\bm V\bm D\bm V^{-1}$ over $\mathbb C$.  Norm equivalence and a union bound reduce the problem to a scalar eigenvalue $\lambda$ with $\operatorname{Re}\lambda>0$.  Put
$x_{j+1}^t=\prod_{i=j+1}^t(1-\lambda\eta_i)$ and let $e_j$ be the corresponding coordinate of $\bm V^{-1}\bm d_j$.  We claim that, uniformly in $k$,
\begin{equation}
 \Pp\left(
 \max_{h_k\le t<h_{k+1}}
 \frac1{\eta_{t+1}}
 \left|\sum_{j=h_k}^t x_{j+1}^t\eta_j e_j\right|
 >\varepsilon\sqrt T\right)
 \le C\varepsilon^{-p}m^{-p/2}.
 \label{eq:appendix-block-innovation}
\end{equation}

To prove the claim, fix the terminal index $T$ and use
$x_{j+1}^t=(x_{t+1}^T)^{-1}x_{j+1}^T$.  For all sufficiently large indices,
$\eta_{t+1}|x_{t+1}^T|$ is nondecreasing because
\[
 \frac{\eta_t}{\eta_{t+1}}|1-\lambda\eta_t|
 =\{1+o(\eta_t)\}\{1-\operatorname{Re}(\lambda)\eta_t+O(\eta_t^2)\}\le1.
\]
Set $b_t=|\eta_{t+1}x_{t+1}^T|^{-p}$ and
$Y_t=|\sum_{j=h_k}^t x_{j+1}^T\eta_j e_j|^p$.  Then $b_t$ is nonincreasing and $Y_t$ is a submartingale.  Burkholder's inequality gives, with
$\ell=h_{k+1}-h_k\le (T+1)/m+1$,
\[
  \E Y_t
  \le C\ell^{p/2-1}\sum_{j=h_k}^t
  |x_{j+1}^T\eta_j|^p
  \le C\ell^{p/2-1}\sum_{j=h_k}^t b_j^{-1}.
\]
Here we used the bounded ratio $\eta_j/\eta_{j+1}$ to absorb the difference between
$|x_{j+1}^T\eta_j|^p$ and $b_j^{-1}$ into $C$.  Chow's maximal inequality for weighted submartingales \citep{chow1960martingale}, together with Burkholder's inequality \citep{burkholder1988sharp}, yields
\begin{align*}
 &\Pp\left(\max_{h_k\le t<h_{k+1}}b_tY_t>\varepsilon^pT^{p/2}\right)\\
 &\quad\le \frac{C\ell^{p/2-1}}{\varepsilon^pT^{p/2}}
 \left\{\sum_{t=h_k}^{h_{k+1}-2}(b_t-b_{t+1})
 \sum_{j=h_k}^tb_j^{-1}
 +b_{h_{k+1}-1}\sum_{j=h_k}^{h_{k+1}-1}b_j^{-1}\right\}\\
 &\quad\le \frac{C\ell^{p/2}}{\varepsilon^pT^{p/2}}
 \le C\varepsilon^{-p}m^{-p/2}.
\end{align*}
The expression in braces equals $\ell$: summation by parts reduces it to
$b_{h_k}b_{h_k}^{-1}+\sum_{t=h_k+1}^{h_{k+1}-1}b_tb_t^{-1} =  h_{k+1}-h_k=\ell$ .
We then prove~\eqref{eq:appendix-block-innovation}.

For $h_k\le t<h_{k+1}$,
$\bm z_{t+1}=\bm X_{h_k}^t\bm z_{h_k}
+\sum_{j=h_k}^t\bm X_{j+1}^t\eta_j\bm d_j$.
On $\mathcal A_{T,m}$,~\eqref{eq:appendix-block-propagation} gives
\[
 \max_{0\le k<m}\max_{h_k\le t<h_{k+1}}
 \frac{\norm{\bm X_{h_k}^t\bm z_{h_k}}}{\sqrt T\eta_{t+1}}
 \le \max_{0\le k<m}C\frac{\norm{\bm z_{h_k}}}
 {\sqrt T\eta_{h_k}}<\varepsilon.
\]
After summing over the finitely many eigencoordinates, the innovation bound
in~\eqref{eq:appendix-block-innovation} continues to hold for the vector norm,
with a different constant $C$.  Hence, the block decomposition and a union
bound yield
\begin{align*}
 &\Pp\left(\max_{t\le T}\frac{\norm{\bm z_{t+1}}}
 {\sqrt T\eta_{t+1}}>2\varepsilon\right)\\
 &\quad\le \Pp(\mathcal A_{T,m}^c)
 +\sum_{k=0}^{m-1}\Pp\left(
 \max_{h_k\le t<h_{k+1}}
 \frac1{\eta_{t+1}}
 \norm{\sum_{j=h_k}^t\bm X_{j+1}^t\eta_j\bm d_j}
 >\varepsilon\sqrt T\right)\\
 &\quad\le \frac{Cm}{\varepsilon^2T\eta_T}
 +C\varepsilon^{-p}m^{1-p/2}.
\end{align*}
For fixed $m$, the first term tends to zero because $T\eta_T\to\infty$.
Therefore,
\[
 \limsup_{T\to\infty}
 \Pp\left(\max_{t\le T}\frac{\norm{\bm z_{t+1}}}
 {\sqrt T\eta_{t+1}}>2\varepsilon\right)
 \le C\varepsilon^{-p}m^{1-p/2}.
\]
Letting $m\to\infty$ proves the lemma because $p>2$.
\end{proof}

\begin{lemma}[Stability transfer]
\label{lem:stability-transfer}
Let $\lambda\in\mathbb C$ have positive real part and let
$z_{t+1}=(1-\lambda\eta_t)z_t+\eta_t\omega_t$, with $z_0=0$.  If
$\max_{t\le T}|\omega_t|/(\sqrt T\eta_t)\to_p0$, then
$\max_{t\le T}|z_{t+1}|/(\sqrt T\eta_{t+1})\to_p0$.
\end{lemma}

\begin{proof}[Proof of Lemma~\ref{lem:stability-transfer}]
For sufficiently large $i$,
$|1-\lambda\eta_i|\le\exp\{-c\eta_i\}$ for some $c>0$.  Therefore,
iteration of the recursion, with
$\eta_j\omega_j=\eta_j^2(\omega_j/\eta_j)$, gives
\begin{align*}
 \max_{t\le T}\frac{|z_{t+1}|}{\sqrt T\eta_{t+1}}
 &\le \max_{t\le T}
 \frac1{\eta_{t+1}}
 \sum_{j=0}^t\eta_j^2
 \prod_{i=j+1}^t|1-\lambda\eta_i|
 \max_{j\le T}\frac{|\omega_j|}{\sqrt T\eta_j}
 \le C\max_{j\le T}\frac{|\omega_j|}{\sqrt T\eta_j}.
\end{align*}
The final inequality is~\eqref{eq:gain-convolution}, with finitely many early indices absorbed into $C$.  The assumed convergence proves the result.
\end{proof}

\begin{proof}[Proof of Lemma~\ref{lem:decreasing-gain}]
In this proof, call a scalar or vector sequence negligible if it satisfies the
convergence in~\eqref{eq:decreasing-gain-negligibility}.
View $\bG$ as a complex matrix and take a Jordan decomposition
$\bG=\bm V\bm J\bm V^{-1}$.  Since the number of Jordan blocks is finite and norms are equivalent, it is enough to treat one block
\[
 \bm J=
 \begin{pmatrix}
 \lambda&1&&\\
 &\lambda&\ddots&\\
 &&\ddots&1\\
 &&&\lambda
 \end{pmatrix},
 \qquad \operatorname{Re}\lambda>0.
\]
Let $\widetilde{\bm z}_t=\bm V^{-1}\bm z_t$ and
$\widetilde{\bm d}_t=\bm V^{-1}\bm d_t$.  By
Lemma~\ref{lem:diagonal-decreasing-gain}, the final coordinate is negligible.

Proceed backward by induction.  Suppose the displayed convergence holds for
coordinate $k$.  Coordinate $k-1$ satisfies
\begin{equation*}
 (\widetilde{\bm z}_{t+1})_{k-1}
 =(1-\lambda\eta_t)(\widetilde{\bm z}_t)_{k-1}
 -\eta_t(\widetilde{\bm z}_t)_k
 +\eta_t(\widetilde{\bm d}_t)_{k-1}.
\end{equation*}
Let $\widehat z_t$ solve the same recursion without the middle term.
Lemma~\ref{lem:diagonal-decreasing-gain} makes it negligible.  The difference
$q_t=(\widetilde{\bm z}_t)_{k-1}-\widehat z_t$ obeys
\[
 q_{t+1}=(1-\lambda\eta_t)q_t-\eta_t(\widetilde{\bm z}_t)_k.
\]
The induction hypothesis and Lemma~\ref{lem:stability-transfer} make $q_t$
negligible.  Hence, coordinate $k-1$ is negligible.  Backward induction covers the entire
block and hence all of $\bm J$.  Transforming back by $\bm V$
proves~\eqref{eq:decreasing-gain-negligibility}.
\end{proof}

\appendixsection{Functional Limit Theorem and Self-Normalized Inference}
\label{app:fclt-proof}

\begin{lemma}[Remainder terms]
\label{lem:four-remainders}
Under Assumptions~\ref{ass:local}--\ref{ass:path}, the remainders in~\eqref{eq:four-remainder-decomposition} satisfy
\[
 \max_{1\le n\le T}\norm{\bm R_{k,T}(n)}\xrightarrow{p}0,
 \qquad k=0,1,2,3.
\]
\end{lemma}

\begin{proof}[Proof of Lemma~\ref{lem:four-remainders}]
Lemma~\ref{lem:recursion-matrices}(ii) immediately gives
\[
 \max_{n\le T}\norm{\bm R_{0,T}(n)}
 \le \frac{C\norm{(\bI-\eta_1\bG)\bDelta_1}}{\sqrt T\eta_1}
 \longrightarrow0.
\]
The same uniform bound and~\eqref{eq:r-nu-negligible} give
\[
 \max_{n\le T}\norm{\bm R_{1,T}(n)}
 \le \frac C{\sqrt T}\sum_{j=1}^T
 \{\norm{\br_j}+\norm{\bnu_j}\}\xrightarrow{p}0.
\]

For $\bm R_{2,T}$, the matrices $\bm A_j^T$ are deterministic, so its partial sums form a martingale.  Doob's inequality, the uniform second moment of $\bepsilon_j$, and Lemma~\ref{lem:recursion-matrices}(ii) yield
\begin{align*}
 \E\max_{n\le T}\norm{\bm R_{2,T}(n)}^2
 &\le \frac C T\sum_{j=1}^T
 \norm{\bm A_j^T-\bG^{-1}}^2\E\norm{\bepsilon_j}^2
 \longrightarrow0.
\end{align*}

It remains to prove~\eqref{eq:r3-to-decreasing-gain}.  For $n<T$, expand the definitions and exchange the sums:
\begin{align*}
 \sum_{j=1}^n(\bm A_j^T-\bm A_j^n)\bepsilon_j
 &=\sum_{j=1}^n\sum_{k=n+1}^T
 \bm X_{j+1}^k\eta_j\bepsilon_j
 =\sum_{k=n+1}^T\bm X_{n+1}^k
 \sum_{j=1}^n\bm X_{j+1}^n\eta_j\bepsilon_j\\
 &=\eta_{n+1}^{-1}\bm A_{n+1}^T(\bI-\eta_{n+1}\bG)
 \sum_{j=1}^n\bm X_{j+1}^n\eta_j\bepsilon_j.
\end{align*}
The commutation used here is valid because every factor
$\bI-\eta_i\bG$ is a polynomial in the same matrix $\bG$.
Lemma~\ref{lem:recursion-matrices}(ii) bounds
$\bm A_{n+1}^T(\bI-\eta_{n+1}\bG)$ uniformly.  Set
$\bm z_0=\bm z_1=\bzero$ and
$\bm z_{n+1}=\sum_{j=1}^n\bm X_{j+1}^n\eta_j\bepsilon_j$ for $n\ge1$.
With $\bm d_0=\bzero$ and $\bm d_n=\bepsilon_n$ for $n\ge1$, this is exactly
recursion~\eqref{eq:decreasing-gain-recursion}.  Hence,
Lemma~\ref{lem:decreasing-gain} and~\eqref{eq:r3-to-decreasing-gain} prove that
the maximum of $\bm R_{3,T}$ is $o_p(1)$.
\end{proof}

\begin{proof}[Proof of Theorem~\ref{thm:fclt}]
Iterating~\eqref{eq:auxiliary-recursion} gives
\[
 \bDelta_{t+1}=\bm X_1^t\bDelta_1
 +\sum_{j=1}^t\bm X_{j+1}^t\eta_j
 \{\bepsilon_j-\br_j-\bnu_j\}.
\]
Summing this equality from $t=1$ to $n$ and exchanging finite sums gives~\eqref{eq:four-remainder-decomposition}.  Lemma~\ref{lem:four-remainders} therefore proves
\begin{equation}
 \max_{n\le T}\frac1{\sqrt T}
 \norm{\sum_{t=1}^n\bDelta_{t+1}
 -\bG^{-1}\sum_{t=1}^n\bepsilon_t}\xrightarrow{p}0.
 \label{eq:appendix-delta-linearization}
\end{equation}
Changing $\sum_{t=1}^n\bDelta_{t+1}$ to
$\sum_{t=1}^n\bDelta_t$ contributes only
$(\bDelta_{n+1}-\bDelta_1)/\sqrt T$, uniformly negligible by~\eqref{eq:interpolation-control}.  Replacing $\bDelta_t$ by $\bx_t-\bx^\star$ is uniformly negligible by~\eqref{eq:auxiliary-negligible}.  Thus~\eqref{eq:appendix-delta-linearization} is exactly~\eqref{eq:uniform-asymptotic-linearity} for the step process.

Lemma~\ref{lem:poisson-martingale} and Slutsky's theorem give convergence in $D([0,1],\R^d)$ to
$\bG^{-1}\bS^{1/2}\bcalW(\cdot)$.  The limit is continuous.  The difference between the linearly interpolated iterate path~\eqref{eq:partial-sum-process} and its step version is bounded by $T^{-1/2}\max_{t\le T+1}\norm{\bx_t-\bx^\star}$, which is $o_p(1)$ by~\eqref{eq:interpolation-control}.  The uniform $p>2$
moment of $\bepsilon_t$ and a union bound similarly give
$T^{-1/2}\max_{t\le T+1}\norm{\bepsilon_t}=o_p(1)$.
Thus linear interpolation of the martingale step process also changes it by
$o_p(1)$ uniformly.  The interpolated processes are random elements of
$C([0,1],\R^d)$, and convergence to a continuous limit is therefore convergence
under the uniform norm.  This proves~\eqref{eq:fclt}
and~\eqref{eq:path-linear-representation}.
\end{proof}

\begin{proof}[Proof of Proposition~\ref{prop:second-order-perturbation}]
Use the same Poisson solution, martingale difference $\bepsilon_t$, residual
$\bnu_t$, coboundary, and auxiliary iterate as in
Section~\ref{sec:proof-architecture}.  The definition of $\bnu_t$ is unchanged, but
the increment bound used in~\eqref{eq:appendix-nu-bound} becomes
\[
 \norm{\bx_{t+1}-\bx_t}
 \le \eta_t\norm{\bH(\bx_t,\xi_t)}+\eta_t^2\norm{\bzeta_t}.
\]
H\"older's inequality, the $p$th-moment assumptions on $L(\xi_t)$,
$\bH(\bx_t,\xi_t)$, and $\bzeta_t$, and boundedness of $(\eta_t)$ therefore preserve the bound
$\E\norm{\bnu_t}\le C\eta_t+o(\eta_t)$.  Thus all conclusions of
Lemma~\ref{lem:poisson-residuals} remain valid.

The only change in the auxiliary recursion~\eqref{eq:auxiliary-recursion} is
\begin{equation*}
 \bDelta_{t+1}=(\bI-\eta_t\bG)\bDelta_t
 +\eta_t\{\bepsilon_t-\br_t-\bnu_t+\eta_t\bzeta_t\}.
\end{equation*}
Consequently the four-remainder expansion is unchanged except that
$\br_j+\bnu_j$ in $\bm R_{1,T}$ is replaced by
$\br_j+\bnu_j-\eta_j\bzeta_j$.  Lemma~\ref{lem:recursion-matrices}(ii), Markov's
inequality, and the uniform first moment of $\bzeta_j$ give
\[
 \max_{n\le T}\frac1{\sqrt T}
 \norm{\sum_{j=1}^n\bm A_j^n\eta_j\bzeta_j}
 \le \frac C{\sqrt T}\sum_{j=1}^T\eta_j\norm{\bzeta_j}
 \xrightarrow{p}0
\]
by~\eqref{eq:gain-conditions}.  The other three remainders and the martingale FCLT are
identical to those in the proof of Theorem~\ref{thm:fclt}.  Hence, the same uniform
linear representation and functional limit follow.
\end{proof}

\begin{proof}[Proof of Corollary~\ref{cor:smooth-functional}]
Let
$r_\psi(\bx)=\psi(\bx)-\psi(\bx^\star)
-\bm v^\top(\bx-\bx^\star)$.  By the definitions of
$\Psi_T^\psi$ and $\bPhi_T$,
\[
 \sup_{0\le r\le1}
 \left|\Psi_T^\psi(r)-\bm v^\top\bPhi_T(r)\right|
 \le \frac1{\sqrt T}\sum_{t=1}^{T+1}|r_\psi(\bx_t)|
 \xrightarrow{p}0
\]
by~\eqref{eq:functional-occupation} applied with $T+1$.  Thus the bound covers
both the partial sums and the single endpoint term introduced by linear
interpolation.  Theorem~\ref{thm:fclt} and the continuous mapping theorem give the asserted limit
$\sigma_{\bm v}\mathcal W(\cdot)$ under the uniform norm.
\end{proof}

\begin{proof}[Proof of Corollary~\ref{cor:coverage}]
Apply the continuous linear map $\bphi\mapsto\bm v^\top\bphi$ to Theorem~\ref{thm:fclt}.  The scalar path converges to $\sigma_{\bm v}\mathcal W(\cdot)$, where
$\sigma_{\bm v}^2=\bm v^\top\bG^{-1}\bS\bG^{-\top}\bm v>0$.
The bridge map
\[
 \phi\mapsto\{\phi(r)-r\phi(1):0\le r\le1\}
\]
is continuous in the uniform norm.  Continuity and scale equivariance of $\mathcal D$
therefore give
\[
  \left(\sqrt T(\bar Y_T-\vartheta^\star),\mathcal D(B_T)\right)
  \Rightarrow
  \left(\sigma_{\bm v}\mathcal W(1),
  \sigma_{\bm v}\mathcal D(\mathcal B)\right).
\]
The endpoint $\mathcal W(1)$ is independent of the Brownian bridge
$\mathcal B$.  Positivity of $\mathcal D(\mathcal B)$ makes the
ratio map almost surely continuous, so the continuous mapping theorem proves
\eqref{eq:pivot-limit}.  Conditional on the positive normalizer,
$Z_{\mathcal D}$ has a centered normal distribution with a positive random
scale.  Hence, $|Z_{\mathcal D}|$ has a continuous distribution, and inversion
at its $(1-\alpha)$ quantile proves the coverage statement.
\end{proof}

\subsection{Polynomial-Series Normalizer and Width Benchmark}
\label{app:series-scaler}

\begin{proof}[Proof of Proposition~\ref{prop:series-interval}]
Let $A_j=\int_0^1\mathcal B(r)r^{j-1}\,\mathrm{d}r$.  The FCLT
and continuity of integration give $\bm a_T\Rightarrow\sigma_{\bm v}\bm A$.
The vector $\bm A$ is Gaussian and
\begin{align*}
 \operatorname{Cov}(A_j,A_k)
 &=\int_0^1\!\int_0^1
 \{\min(r,s)-rs\}r^{j-1}s^{k-1}\,\mathrm{d}r\,\mathrm{d}s=\frac{1}{(j+1)(k+1)(j+k+1)},
\end{align*}
which proves~\eqref{eq:series-covariance}.  Since the monomials are linearly
independent, $\bm\Sigma_K$ is positive definite.  Thus
$\bm\Sigma_K^{-1/2}\bm a_T/\sigma_{\bm v}\Rightarrow N(\bzero,\bI_K)$, so
$\mathcal D_{P_K,T}^2/\sigma_{\bm v}^2\Rightarrow\chi_K^2/K$.  Moreover,
$\mathcal W(1)$ is independent of every projection of $\mathcal B$.
Dividing the independent standard normal limit by
$\sqrt{\chi_K^2/K}$ gives the $t_K$ limit and hence the coverage of
\eqref{eq:confidence-interval}.

The limiting series normalizer satisfies
$D\overset{d}=\sqrt{\chi_K^2/K}$, and its critical value is
$t_{K,1-\alpha/2}$.  The standard chi-distribution mean then gives $ t_{K,1-\alpha/2}\E D
=t_{K,1-\alpha/2}\sqrt{\frac2K}
\frac{\Gamma((K+1)/2)}{\Gamma(K/2)}.$
The standard limits for the Student quantile and the gamma ratio prove the
stated convergence to $\Phi^{-1}(1-\alpha/2)$.
\end{proof}

\begin{proof}[Proof of Proposition~\ref{prop:width-benchmark}]
Write $D=\mathcal D(\mathcal B)$ and
$c=c_{1-\alpha,\mathcal D}$.  Recall that
$Z_{\mathcal D}=\mathcal W(1)/D$ and that $c$ is the $(1-\alpha)$ quantile of
$|Z_{\mathcal D}|$.  Since $\mathcal W(1)$ is standard normal and independent
of $D$, conditioning on $D$ gives
\[
  1-\alpha=\Pp(|Z_{\mathcal D}|\le c)
  =\E\!\left[\Pp\{|\mathcal W(1)|\le cD\mid D\}\right]
  =\E G(cD),
\]
where $G(h):=\Pp\{|N(0,1)|\le h\}=2\Phi(h)-1$.
The function $G$ is strictly concave on $(0,\infty)$.  Jensen's inequality
therefore gives
\[
  1-\alpha\le G(c\E D),
  \qquad c\E D\ge \Phi^{-1}(1-\alpha/2),
\]
which is the lower bound.  Equality would require $cD$ to be almost surely
constant, so a nondegenerate random normalizer cannot attain equality.
\end{proof}

\paragraph{Online $P_K$ Implementation.}
To justify Remark~\ref{rem:series-online}, define the step approximation
$\widetilde B_T(r):=T^{-1/2}\sum_{t=1}^{\lfloor Tr\rfloor}(Y_t-\bar Y_T)$.
For $r\in[n/T,(n+1)/T)$, direct subtraction gives
$B_T(r)-\widetilde B_T(r)=(Tr-n)\{Y_{n+1}-\bar Y_T\}/\sqrt T$.
Because $L(\xi)\ge1$, the moment bound in Assumption~\ref{ass:path}(i) implies
$\sup_t\E\norm{\bx_t-\bx^\star}^p<\infty$.  A union bound with $p>2$ yields
$T^{-1/2}\max_{t\le T+1}\norm{\bx_t-\bx^\star}=o_p(1)$, and therefore
\begin{equation}
	\norm{B_T-\widetilde B_T}_\infty
	:=\sup_{0\le r\le1}|B_T(r)-\widetilde B_T(r)|=o_p(1).
	\label{eq:step-bridge-equivalence}
\end{equation}
Let $\widetilde a_{j,T}=\int_0^1\widetilde B_T(r)r^{j-1}\,\mathrm{d}r$,
let $\widetilde{\bm a}_T=(\widetilde a_{1,T},\ldots,\widetilde a_{K,T})^\top$,
and define
$\widetilde{\mathcal D}_{P_K,T}^2=K^{-1}\widetilde{\bm a}_T^\top
\bm\Sigma_K^{-1}\widetilde{\bm a}_T$.
Equation~\eqref{eq:step-bridge-equivalence} gives
$\widetilde{\mathcal D}_{P_K,T}-\mathcal D_{P_K,T}=o_p(1)$ for fixed $K$.
Since $\sum_{t=1}^T(Y_t-\bar Y_T)=0$, integration of the step approximation
gives the exact identity
\begin{equation*}
	\widetilde a_{j,T}
	=-\frac{M_{j,T}-\bar Y_TP_{j,T}}{jT^{j+1/2}},
	\qquad
	M_{j,T}=\sum_{t=1}^Tt^jY_t,
	\quad P_{j,T}=\sum_{t=1}^Tt^j.
\end{equation*}
For numerical stability, define
$m_{j,T}=T^{-(j+1)}M_{j,T}$ and $p_{j,T}=T^{-(j+1)}P_{j,T}$.  They obey
\begin{align*}
	m_{j,T}
	&=\left(\frac{T-1}{T}\right)^{j+1}m_{j,T-1}+\frac{Y_T}{T},\\
	p_{j,T}
	&=\left(\frac{T-1}{T}\right)^{j+1}p_{j,T-1}+\frac1T,
	\qquad
	\widetilde a_{j,T}=-\frac{\sqrt T}{j}(m_{j,T}-\bar Y_Tp_{j,T}).
\end{align*}
Together with the usual update for $\bar Y_T$, these recursions use $2K+1$
scalars and $O(K)$ arithmetic per observation.  The coordinates
$\widetilde a_{j,T}$ above project the bridge onto the monomial basis
$\{1,r,\ldots,r^{K-1}\}$.  For numerical stability, our implementation instead
uses the shifted-Legendre basis of the same polynomial space.  Any change
between these two bases is nonsingular and transforms both
$\widetilde{\bm a}_T$ and its covariance in the same way, leaving
$\widetilde{\bm a}_T^\top\bm\Sigma_K^{-1}\widetilde{\bm a}_T$ unchanged.
Thus the change of basis does not alter the resulting interval.

\subsection{Bridge-Baseline Details}
\label{app:lm-baselines}

For an even integer $m\ge2$, write
$\|f\|_{L_m}:=(\int_0^1|f(r)|^m\,\mathrm{d}r)^{1/m}$.

\begin{table}[t]
\centering
\small
\setlength{\tabcolsep}{5pt}
\begin{tabular}{lrrr@{\hspace{9pt}\vrule\hspace{9pt}}lrrr}
\toprule
Functional & $90\%$ & $95\%$ & $99\%$
& Functional & $90\%$ & $95\%$ & $99\%$ \\
\midrule
$P_5$ & 2.015 & 2.571 & 4.032
& $L_2$ & 5.341 & 6.762 & 10.069 \\
$P_{10}$ & 1.812 & 2.228 & 3.169
& $L_4$ & 4.264 & 5.382 & 7.911 \\
$P_{20}$ & 1.725 & 2.086 & 2.845
& $L_6$ & 3.755 & 4.730 & 6.918 \\
\bottomrule
\end{tabular}
\caption{Two-sided critical values for the six functionals compared in
Figure~\ref{fig:series-ablation}.  The $P_K$ rows use standard Student
quantiles with $K$ degrees of freedom.  The $L_m$ values are estimated from
50,000 Brownian paths on a grid of 1,000 time points.}
\label{tab:critical-values}
\vspace{-10pt}
\end{table}

\begin{lemma}[Online $L_m$ normalizers for even $m$]
\label{lem:online-denominator}
Under the conditions of Theorem~\ref{thm:fclt}, for every fixed even integer
$m\ge2$,
\begin{equation}
  \mathcal D_{L_m,T}:=\|B_T\|_{L_m},\qquad
  \widehat{\mathcal D}_{m,T}:=
  \left\{T^{-1}\sum_{n=1}^T|B_T(n/T)|^m\right\}^{1/m},\qquad
  \widehat{\mathcal D}_{m,T}-\mathcal D_{L_m,T}\xrightarrow{p}0.
  \label{eq:online-denominator-equivalence}
\end{equation}
If $\sigma_{\bm v}^2>0$, replacing $\mathcal D_{L_m,T}$ by
$\widehat{\mathcal D}_{m,T}$ in Corollary~\ref{cor:coverage} preserves its limiting pivot and
coverage conclusion.
\end{lemma}

\begin{proof}[Proof of Lemma~\ref{lem:online-denominator}]
Let $\omega(f,\delta)=\sup_{|r-s|\le\delta}|f(r)-f(s)|$.  On each interval
$((n-1)/T,n/T]$, linear interpolation and the inequality $ \big||a|^m-|b|^m\big|
\le m\max(|a|,|b|)^{m-1}|a-b|$ give
\begin{equation}
 \left|\mathcal D_{L_m,T}^m-\widehat{\mathcal D}_{m,T}^m\right|
 \le m\norm{B_T}_\infty^{m-1}\omega(B_T,T^{-1}).
 \label{eq:appendix-riemann-bound}
\end{equation}
Theorem~\ref{thm:fclt} and the continuous bridge map imply that $B_T$ converges
in $C([0,1])$ to $\sigma_{\bm v}\{\mathcal W(r)-r\mathcal W(1)\}$.  Hence,
$\norm{B_T}_\infty=O_p(1)$ and asymptotic equicontinuity gives
$\omega(B_T,T^{-1})=o_p(1)$.  Equation~\eqref{eq:appendix-riemann-bound} proves
$|\mathcal D_{L_m,T}^m-\widehat{\mathcal D}_{m,T}^m|=o_p(1)$.  The inequality
$|a^{1/m}-b^{1/m}|\le|a-b|^{1/m}$ for $a,b\ge0$ proves
\eqref{eq:online-denominator-equivalence}.
When $\sigma_{\bm v}^2>0$, the common limit of the normalizers is the strictly
positive random variable
$\sigma_{\bm v}\{\int_0^1|\mathcal W(r)-r\mathcal W(1)|^m\,\mathrm{d}r\}^{1/m}$.
Slutsky's theorem therefore permits replacing $\mathcal D_{L_m,T}$ by
$\widehat{\mathcal D}_{m,T}$ in the pivot and interval.
\end{proof}

\paragraph{Online $L_m$ Implementation.}
For Remark~\ref{rem:lm-online}, let $S_n=\sum_{t=1}^nY_t$ and
$R_{\ell,T}^{(m)}=\sum_{n=1}^Tn^{m-\ell}S_n^\ell$ for
$\ell=0,\ldots,m$.  Since $m$ is even, the binomial expansion gives
\[
 \widehat{\mathcal D}_{m,T}^m
 =T^{-1-m/2}\sum_{\ell=0}^m
 \binom m\ell(-\bar Y_T)^{m-\ell}R_{\ell,T}^{(m)}.
\]
Each $R_{\ell,T}^{(m)}$ is updated by adding $T^{m-\ell}S_T^\ell$.  Maintaining
$S_T$ and $R_{\ell,T}^{(m)}$ for $\ell=0,\ldots,m$ therefore requires $m+2$
scalar running sums and $O(m)$ operations per iteration.

\begin{proof}[Proof of Proposition~\ref{prop:lm-interval}]
The map $f\mapsto\{\int_0^1|f(r)|^m\,\mathrm{d}r\}^{1/m}$ is continuous,
scale-equivariant, and positive almost surely at a standard Brownian
bridge.  Corollary~\ref{cor:coverage} therefore gives the stated pivot limit and
coverage.  For fixed even $m$, Lemma~\ref{lem:online-denominator} permits replacing
$\mathcal D_{L_m,T}$ by $\widehat{\mathcal D}_{m,T}$ and gives the constant-memory implementation.
\end{proof}

\appendixsection{Verification of the Examples}
\label{app:examples}

\begin{proof}[Proof of Proposition~\ref{prop:tabular-q}]
Under the fixed full-support behavior policy $\pi_b$, $(S_t,A_t)$ is an irreducible, aperiodic finite Markov chain with stationary mass
$d_{\pi_b}(s)\pi_b(a\mid s)>0$.  Augmenting the state by the conditionally generated reward and next state preserves uniform geometric ergodicity.  With $V\equiv1$, bounded rewards and the Lipschitz continuity of the maximum operator verify~\eqref{eq:h-v-lipschitz}.  Proposition~\ref{prop:v-route} therefore gives Assumption~\ref{ass:poisson}.

We verify boundedness.  Since $\eta_t\to0$, choose a deterministic $t_0$ such
that $0<\eta_t\le1$ for $t\ge t_0$.  Bounded rewards and the finite prefix give
a deterministic $M_0<\infty$ such that $\norm{\bq_{t_0}}_\infty\le M_0$ almost
surely.  Let $M=\max\{M_0,R_{\max}/(1-\gamma)\}$.  If
$\norm{\bq_t}_\infty\le M$ for $t\ge t_0$, then
$|R_t+\gamma\max_a q_t(S_{t+1},a)|\le R_{\max}+\gamma M\le M$.  The updated
coordinate is a convex combination of its old value and this quantity, while
all other coordinates are unchanged.
Induction gives $\sup_{t\ge t_0}\norm{\bq_t}_\infty\le M$; the finite prefix is
also bounded.  Full exploration, the contraction of the Bellman optimality operator, and Assumption~\ref{ass:gain} then imply
$\bq_t\to\bq^\star$ almost surely by the standard asynchronous Q-learning convergence theorem \citep{tsitsiklis1994asynchronous,even2003learning}.

Define the action gap by
\[
 \Delta_{\mathrm{gap}}
 =\min_{s\in\cS}\{q^\star(s,a^\star(s))-
 \max_{a\ne a^\star(s)}q^\star(s,a)\}>0.
\]
When $\norm{\bq-\bq^\star}_\infty<\Delta_{\mathrm{gap}}/3$, the maximizing action is $a^\star(s)$ in every state.  In that neighborhood the mean field~\eqref{eq:q-mean-field} is affine with Jacobian $\bG_Q=\bD_{\pi_b}(\bI_{d_Q}-\gamma\bP\bPi^\star)$.  The matrix
$\bI_{d_Q}-\gamma\bP\bPi^\star$ has nonpositive off-diagonal entries, and its inverse is the nonnegative Neumann series
$\sum_{k\ge0}(\gamma\bP\bPi^\star)^k$.  Left multiplication by the positive diagonal matrix $\bD_{\pi_b}$ preserves the off-diagonal signs and gives inverse
$(\bI_{d_Q}-\gamma\bP\bPi^\star)^{-1}\bD_{\pi_b}^{-1}\ge0$.  Thus
$\bG_Q$ has nonpositive off-diagonal entries and $\bG_Q^{-1}$ is elementwise
nonnegative, so it is a nonsingular $M$-matrix.  Consequently, every eigenvalue of $\bG_Q$ has
positive real part \citep[Chapter~6]{berman1994nonnegative}, verifying
Assumption~\ref{ass:local}.

Almost-sure convergence and boundedness imply~\eqref{eq:time-averaged-l2} by dominated convergence.  Define the almost-surely finite entrance time
$\tau_{\mathrm{gap}}:=\inf\{t\ge1:\norm{\bq_j-\bq^\star}_\infty<
\Delta_{\mathrm{gap}}/3\text{ for every }j\ge t\}$.
Then $\br_g(\bq_t)=\bzero$ for $t\ge\tau_{\mathrm{gap}}$, while boundedness
controls the finitely many earlier terms.  Consequently, $\frac1{\sqrt T}\sum_{t=1}^T\norm{\br_g(\bq_t)}
\le \frac1{\sqrt T}\sum_{t=1}^{\tau_{\mathrm{gap}}-1}
\norm{\br_g(\bq_t)}\longrightarrow0$ almost surely.
Thus~\eqref{eq:integrated-nonlinearity} holds.  All moment conditions
in~\eqref{eq:path-moments} follow from boundedness.  Assumption~\ref{ass:gain}
is imposed.  Therefore, every assumption of Theorem~\ref{thm:fclt} is verified.
\end{proof}

\begin{proof}[Proof of Proposition~\ref{prop:linear-q}]
Bounded features and rewards make the unprojected update globally Lipschitz in
$\btheta$, with a bounded value at $\btheta^\star$.  For the finite irreducible
and aperiodic state--action chain, these properties verify
\eqref{eq:h-v-lipschitz} with $V\equiv1$.  Proposition~\ref{prop:v-route}
therefore gives Assumption~\ref{ass:poisson}.
The corresponding mean field is
\[
 \bg(\btheta)=\E_{\pi_b}\!\left[\bvarphi(S,A)
 \left\{\bvarphi(S,A)^\top\btheta-R
 -\gamma\max_{a'\in\cA}\bvarphi(S',a')^\top\btheta\right\}\right].
\]
Put $\bu=\btheta-\btheta^\star$.  Since $\bg(\btheta^\star)=\bzero$, the
maximum operator bound
$|\max_a\bvarphi(s,a)^\top\btheta\allowbreak-
\max_a\bvarphi(s,a)^\top\btheta^\star|
\le\max_a|\bvarphi(s,a)^\top\bu|$.  Hence, stationarity and
Cauchy--Schwarz give
\[
 \bu^\top\bg(\btheta)
 \ge \E_{\pi_b}[\{\bvarphi(S,A)^\top\bu\}^2]
 -\gamma\sqrt{\E_{\pi_b}[\{\bvarphi(S,A)^\top\bu\}^2]
 \sum_s d_{\pi_b}(s)\max_a\{\bvarphi(s,a)^\top\bu\}^2}>0,
\]
where the strict inequality follows from
\eqref{eq:linear-q-global-condition}.  The right-hand side is continuous and
homogeneous of degree two in $\bu$, so compactness of the unit sphere yields
$\bu^\top\bg(\btheta)\ge c\norm{\bu}^2$ for some $c>0$.  This is precisely the
global stability inequality in
\citet[Proposition~3.1(3)]{chen2022finite}.  The compact projection and the
finite-state Markov noise satisfy the boundedness and averaging conditions of
\citet[Theorem~6.1.1]{kushner2003stochastic}.  Together with the preceding
stability inequality and the nonexpansiveness of Euclidean projection, that
theorem gives $\btheta_t\to\btheta^\star$ almost surely.  Bounded updates and $\eta_t\to0$
make the unprojected increments vanish; since $\btheta^\star$ is interior to
$\cC$, the projection is inactive after a finite random time.

Condition~(iii) of Proposition~\ref{prop:linear-q} and continuity of the linear Q-model fix the greedy action
in a neighborhood of $\btheta^\star$.  In this neighborhood, the mean field is
affine with Jacobian
\[
  \bG_\varphi=
  \E_{\pi_b}\!\left[\bvarphi(S,A)
  \{\bvarphi(S,A)-\gamma\bvarphi(S',a^\star(S'))\}^\top\right].
\]
The preceding global bound now gives
$\bu^\top\bG_\varphi\bu\ge c\norm{\bu}^2$.  Thus the symmetric part of
$\bG_\varphi$ is positive definite, so every
eigenvalue of $\bG_\varphi$ has positive real part.  Moreover,
$\br_g(\btheta)=\bzero$ in this neighborhood.  Compactness and almost-sure
convergence verify~\eqref{eq:path-moments} and
\eqref{eq:time-averaged-l2} by dominated convergence.  After an almost surely
finite random time, the greedy action is fixed and the projection is inactive,
so $\br_g(\btheta_t)=\bzero$.  The finitely many earlier remainder terms are
negligible after division by $\sqrt T$.  Hence,
Assumptions~\ref{ass:local}--\ref{ass:path} hold, and deleting a
finite prefix does not change the FCLT or the interval.
\end{proof}

\begin{proof}[Proof of Proposition~\ref{prop:soft-q}]
The soft Bellman operator
$\mathcal T_\tau\bq=\br+\gamma\bP\mathcal M_\tau\bq$ is a $\gamma$-contraction in
the sup norm because log-sum-exp is one-Lipschitz in that norm.  Consequently it has
a unique fixed point $\bq_\tau^\star$, and the standard asynchronous contraction
argument under full exploration gives $\bq_t\to\bq_\tau^\star$ almost surely
\citep{tsitsiklis1994asynchronous}.  The
boundedness argument in the proof of Proposition~\ref{prop:tabular-q} applies after
replacing $R_{\max}/(1-\gamma)$ by a finite bound that also includes
$\gamma\tau\log|\cA|/(1-\gamma)$.  Hence, all iterates are bounded.  The update at
$\bq_\tau^\star$ is bounded, and the global Lipschitz property of log-sum-exp
makes the update globally Lipschitz in $\bq$.  The finite irreducible and
aperiodic state--action chain then
satisfies Assumption~\ref{ass:poisson} by the uniform-geometric case of
Proposition~\ref{prop:v-route}.

For a perturbation $\bm h\in\R^{d_Q}$, differentiation of the log-sum-exp
operator in Proposition~\ref{prop:soft-q} at $\bq_\tau^\star$ gives $ D\mathcal M_\tau(\bq_\tau^\star)[\bm h]=\bPi_\tau^\star\bm h.$
Therefore, the stationary mean field has Jacobian
$\bG_\tau=\bD_{\pi_b}(\bI_{d_Q}-\gamma\bP\bPi_\tau^\star)$.  The
matrix $\bP\bPi_\tau^\star$ is stochastic, meaning that its entries are
nonnegative and every row sums to one.  The same argument as in the proof of
Proposition~\ref{prop:tabular-q} shows that $\bG_\tau$ has nonpositive
off-diagonal entries and an elementwise nonnegative inverse.  Thus $\bG_\tau$
is a nonsingular $M$-matrix, so every eigenvalue has strictly positive real
part.

The Hessian of log-sum-exp is $\tau^{-1}$ times the covariance matrix of the
softmax probabilities and is uniformly bounded for fixed $\tau>0$.
Consequently, the first-order remainder of the mean field is locally quadratic:
$\norm{\br_g(\bq)}\le C\norm{\bq-\bq_\tau^\star}^2$.  The contraction and
Lipschitz properties above, finite-state geometric mixing, and bounded reward
noise verify the conditions of
\citet[Theorem~1(3)]{chen2021lyapunov}.
Applying that result to a sufficiently late tail of the polynomial-step-size
recursion and using equivalence of norms in $\R^{d_Q}$ give
$\E\norm{\bq_t-\bq_\tau^\star}^2\le
C\log(t+1)/(t+1)^\kappa$.  Consequently,
$T^{-1/2}\sum_{t=1}^T\E\norm{\bq_t-\bq_\tau^\star}^2\to0$ because
$\kappa>1/2$.  Markov's inequality verifies
\eqref{eq:remainder-occupation-route}.  Together with the quadratic remainder
bound, this proves~\eqref{eq:integrated-nonlinearity},
and the same mean-square error bound gives
\eqref{eq:time-averaged-l2}.  Boundedness gives~\eqref{eq:path-moments}.  Thus
Assumptions~\ref{ass:local}--\ref{ass:path} all hold, without requiring unique optimal actions.
\end{proof}

\begin{proof}[Proof of Proposition~\ref{prop:markov-glm}]
Let $\widetilde{\bH}_{\mathrm{glm}}$ be the extension in the proposition.  It
agrees with $\bH_{\mathrm{glm}}$ along the projected recursion, and
Proposition~\ref{prop:v-route} applied to $\widetilde{\bH}_{\mathrm{glm}}$ gives
Assumption~\ref{ass:poisson}.  The derivative bounds stated in
Proposition~\ref{prop:markov-glm} justify differentiation under the stationary expectation
and give
$\bG_{\mathrm{glm}}=\E_\mu[b''(\bz^\top\btheta^\star)\bz\bz^\top]$.
The same lower bound is uniform over $\btheta\in\cC$ and implies
$\lambda_{\min}(\bG_{\mathrm{glm}})>0$.  The derivative bounds also make the
second derivative of $\bg$ uniformly bounded on $\cC$, so we have
$\norm{\br_g(\btheta)}\le C\norm{\btheta-\btheta^\star}^2$ there.

The uniform Hessian lower bound also gives
$(\btheta-\btheta^\star)^\top\bg(\btheta)\ge
c\norm{\btheta-\btheta^\star}^2$ on $\cC$.  The compact projection and the
Poisson decomposition from Proposition~\ref{prop:v-route} verify the
boundedness and correlated-noise averaging conditions of
\citet[Theorem~6.1.1]{kushner2003stochastic}.  That theorem and the preceding
inequality give $\btheta_t\to\btheta^\star$ almost surely.
Compactness of $\cC$ and~\eqref{eq:h-v-lipschitz} give a uniform $p$th moment for the
stochastic gradient.  Since $\sum_t\eta_t^p<\infty$, Borel--Cantelli yields
$\eta_t\norm{\bH_{\mathrm{glm}}(\btheta_t,\bz_t,y_t)}\to0$ almost surely.
The interior condition on $\btheta^\star$ then makes the projection eventually
inactive.

Adapting the mean-square analysis of \citet[Lemma~B.4]{roy2023online} to the
fixed-kernel projected recursion gives
$\E\norm{\btheta_t-\btheta^\star}^2=O(\eta_t)$.  Compactness of $\cC$ and the bounds
assumed in Proposition~\ref{prop:markov-glm} give~\eqref{eq:path-moments},
while the mean-square error bound gives
\eqref{eq:time-averaged-l2} and
$T^{-1/2}\sum_{t=1}^T\E\norm{\btheta_t-\btheta^\star}^2
=O(T^{1/2-\kappa})\to0$.  Markov's inequality verifies
\eqref{eq:remainder-occupation-route}, which, together with the quadratic
remainder bound, proves~\eqref{eq:integrated-nonlinearity}.
\end{proof}

\begin{proof}[Proof of Proposition~\ref{prop:lora}]
The product formed by the first two factor updates has rank at most $r$, so
the balancing map in~\eqref{eq:balanced-lora-update} is well defined.  Write
$\widetilde{\bW}=\bL\bSigma\bR^\top$, where $\bL$ and $\bR$ have $r$
orthonormal columns, padding $\bSigma$ with zeros if needed.  Set
$\bB=\bL\bSigma^{1/2}$ and $\bA=\bSigma^{1/2}\bR^\top$.  Then
$\bB\bA=\widetilde{\bW}$ and $\bB^\top\bB=\bA\bA^\top=\bSigma$.

Suppress the time index and write $\bW=\bB\bA$ and
$\bGamma_\xi=\bGamma(\bW,\xi)$.  Multiplying the first two updates in
\eqref{eq:balanced-lora-update} gives
\begin{align}
 (\bB-\eta\bGamma_\xi\bA^\top)(\bA-\eta\bB^\top\bGamma_\xi)
 =\bW-\eta(\bGamma_\xi\bA^\top\bA+\bB\bB^\top\bGamma_\xi)
 +\eta^2\bGamma_\xi\bW^\top\bGamma_\xi.
 \label{eq:appendix-lora-product}
\end{align}
Balancing preserves this product.  If $\bB^\top\bB=\bA\bA^\top$, then
$\bW\bW^\top=\bB\bA\bA^\top\bB^\top
=\bB(\bB^\top\bB)\bB^\top=(\bB\bB^\top)^2$.
Uniqueness of the positive-semidefinite square root gives
$\bB\bB^\top=\bmM_L(\bW)$ and, similarly,
$\bA^\top\bA=\bmM_R(\bW)$.  Substitution proves
\eqref{eq:lora-product-recursion}.

\paragraph{Local Stability.}
The stationary identities in Proposition~\ref{prop:lora} give
$\overline{\bGamma}(\bW):=\E_\mu\bGamma(\bW,\xi)=\bW-\bW^\star$.
Define the mean field of the order-$\eta_t$ product update by
$\bm F(\bW):=\overline{\bGamma}(\bW)\bmM_R(\bW)
+\bmM_L(\bW)\overline{\bGamma}(\bW)$.
Because $\overline{\bGamma}(\bW^\star)=\bzero$, the derivatives of
$\bmM_L$ and $\bmM_R$ vanish from $D\bm F(\bW^\star)$.  Hence,
$D\bm F(\bW^\star)[\bDelta]=\mathcal L_\star[\bDelta]$, where
$\mathcal L_\star$ is defined in~\eqref{eq:lora-local-operator}.

Let $\bW^\star=\bL\bSigma\bR^\top$ be its rank-$r$ compact SVD, where
$\bSigma:=\diag(\sigma_1,\ldots,\sigma_r)$ and $\sigma_r>0$.  Let
$\bL_\perp$ and $\bR_\perp$ complete orthonormal bases, and denote by
$\mathcal T_\star$ the tangent space of the rank-$r$ matrices at
$\bW^\star$, that is, the set of first-order perturbations that preserve rank.
Every $\bDelta\in\mathcal T_\star$ has a unique representation
\begin{equation}
 \bDelta=\bL\bK\bR^\top+\bL_\perp\bC_L\bR^\top
 +\bL\bC_R\bR_\perp^\top.
 \label{eq:appendix-rank-tangent}
\end{equation}
Here $\bmM_{L,\star}=\bL\bSigma\bL^\top$ and
$\bmM_{R,\star}=\bR\bSigma\bR^\top$, whence
\begin{align*}
 \langle\bDelta,\mathcal L_\star[\bDelta]\rangle_F
 &=\operatorname{tr}(\bDelta^\top\bDelta\bmM_{R,\star})
   +\operatorname{tr}(\bDelta^\top\bmM_{L,\star}\bDelta)
      =\norm{\bDelta\bmM_{R,\star}^{1/2}}_F^2
   +\norm{\bmM_{L,\star}^{1/2}\bDelta}_F^2.
\end{align*}
If this quantity is zero, then $\bDelta\bR=\bzero$ and
$\bL^\top\bDelta=\bzero$, so~\eqref{eq:appendix-rank-tangent} gives
$\bK=\bC_L=\bC_R=\bzero$.  Thus $\mathcal L_\star$ is positive definite on
$\mathcal T_\star$.

\paragraph{Local-Coordinate SA.}
To express the rank-$r$ product recursion in Euclidean coordinates, let
$\mathcal M_r$ denote the set of rank-$r$ matrices.  The three coefficient
matrices in~\eqref{eq:appendix-rank-tangent} have respectively $r^2$,
$(m-r)r$, and $r(n-r)$ free entries, so $\mathcal M_r$ has dimension
$r(m+n-r)$.  Choose
open neighborhoods $\Theta\subset\R^{r(m+n-r)}$ of $\btheta^\star$ and
$\mathcal U\subset\mathcal M_r$ of $\bW^\star$, together with a smooth
one-to-one parametrization $\chi:\Theta\to\mathcal U$ satisfying
$\chi(\btheta^\star)=\bW^\star$.  Thus $\btheta=\chi^{-1}(\bW)$ is only a
local coordinate for the product $\bW$.  The derivative $D\chi(\btheta)$ maps
the coordinate space bijectively onto the tangent space of $\mathcal M_r$ at
$\chi(\btheta)$; below, $[D\chi(\btheta)]^{-1}$ denotes its inverse on that
tangent space.

Define the order-$\eta_t$ update in these coordinates by
\[
 \bH(\btheta,\xi)
 =[D\chi(\btheta)]^{-1}
 \{\bGamma(\chi(\btheta),\xi)\bmM_R(\chi(\btheta))
 +\bmM_L(\chi(\btheta))\bGamma(\chi(\btheta),\xi)\},
\]
and let $\bg(\btheta):=\E_\mu\bH(\btheta,\xi)$ and
$\bG:=D\bg(\btheta^\star)$.  For every
$\bm h\in\R^{r(m+n-r)}$,
\[
 \bG\bm h=[D\chi(\btheta^\star)]^{-1}
 \mathcal L_\star[D\chi(\btheta^\star)\bm h].
\]
Here $D\chi(\btheta^\star)$ converts the coordinate perturbation $\bm h$ into
a tangent perturbation of $\bW^\star$, and its inverse converts the result of
$\mathcal L_\star$ back into coordinates.  Hence, the matrices of
$\bG$ and $\mathcal L_\star$ restricted to $\mathcal T_\star$ differ
only by this change of coordinates and have the same eigenvalues.  Since
$\langle\bDelta,\mathcal L_\star[\bDelta]\rangle_F>0$ for every nonzero
$\bDelta\in\mathcal T_\star$, every eigenvalue of $\mathcal L_\star$ on
$\mathcal T_\star$, and hence every eigenvalue of $\bG$, is positive.
This verifies Assumption~\ref{ass:local}.

Once $\bW_t\in\mathcal U$, Taylor expansion of $\chi^{-1}$ in
\eqref{eq:lora-product-recursion} gives
\begin{equation}
 \btheta_{t+1}=\btheta_t-\eta_t\bH(\btheta_t,\xi_t)
 +\eta_t^2\bzeta_t,
\label{eq:lora-coordinate-drift}
\end{equation}
where $\bzeta_t$ is bounded whenever $\btheta_t$ is restricted to a compact
subset of $\Theta$.  Equation~\eqref{eq:lora-coordinate-drift} is the SA
recursion~\eqref{eq:sa-recursion} plus the second-order perturbation allowed by
Proposition~\ref{prop:second-order-perturbation}.

\begin{lemma}[Localized accumulated squared-error bound]
\label{lem:localized-squared-error}
Let $\{\xi_t\}$ be an irreducible, aperiodic finite-state Markov chain and
$\eta_t=\eta(t+1)^{-\kappa}$ for $\kappa\in(1/2,1)$.  Suppose that
$\bx_t\to\bx^\star$ almost surely and, whenever $\bx_t$ lies in a
neighborhood of $\bx^\star$,
$\bx_{t+1}=\bx_t-\eta_t\bH(\bx_t,\xi_t)
+\eta_t^2\bzeta_t$.  Suppose locally that $\bH(\cdot,\xi)$ is $C^2$ for
every $\xi$, $\bg:=\E_\mu\bH$ satisfies $\bg(\bx^\star)=\bzero$,
$D\bg(\bx^\star)$ has eigenvalues with positive real parts, and
for every compact $K$, some deterministic $C_K<\infty$ satisfies
$\norm{\bzeta_t}\le C_K$ whenever $\bx_t\in K$.
Then
$T^{-1/2}\sum_{t=1}^T\norm{\bx_t-\bx^\star}^2\xrightarrow{p}0$.
\end{lemma}

\paragraph{Path Estimates and Inference.}
The assumed almost-sure convergence $\bW_t\to\bW^\star$ and Weyl's inequality imply that
$\sigma_r(\bW_t)\ge\sigma_r(\bW^\star)/2$ and $\bW_t\in\mathcal U$
eventually.  Hence, $\btheta_t=\chi^{-1}(\bW_t)$ is eventually well defined
and converges to $\btheta^\star$.  On every compact neighborhood of
$\btheta^\star$, $\bH$ is smooth and $\bzeta_t$ is bounded.  Together
with the finite-state Markov chain and the positivity of $\bG$ proved
above, these facts allow us to apply Lemma~\ref{lem:localized-squared-error}
directly to~\eqref{eq:lora-coordinate-drift}.  Therefore,
$T^{-1/2}\sum_{t=1}^T\norm{\btheta_t-\btheta^\star}^2=o_p(1)$.

Finite-state observations and local smoothness give, along the eventually
local coordinate recursion, the bounds in
\eqref{eq:u-weighted-lipschitz}, \eqref{eq:poisson-moments}, and
\eqref{eq:path-moments}.  The localized $O(\eta_t)$ estimate established in
Lemma~\ref{lem:localized-squared-error} gives
\eqref{eq:time-averaged-l2}.  Define
$\br_g(\btheta):=\bg(\btheta)
-\bG(\btheta-\btheta^\star)$.  Since $\bg$ is $C^2$ near
$\btheta^\star$,
$\norm{\br_g(\btheta_t)}\le
C\norm{\btheta_t-\btheta^\star}^2$.  The lemma therefore verifies
\eqref{eq:integrated-nonlinearity}.  Proposition~\ref{prop:second-order-perturbation}
and Theorem~\ref{thm:fclt} therefore apply after the almost-surely finite time
from which $\btheta_t$ remains in the coordinate neighborhood.  Removing this
finite initial segment changes the normalized partial-sum path by $o_p(1)$,
so the FCLT holds for the original partial sums.

Finally, define
$\psi_{\bC}(\btheta):=\langle\bC,\chi(\btheta)\rangle_F$, so that
$Y_t^{\bC}=\psi_{\bC}(\btheta_t)$ and
$\vartheta_{\bC}^\star=\psi_{\bC}(\btheta^\star)$.  Because $\psi_{\bC}$ is
$C^2$, its first-order remainder is bounded by
$C\norm{\btheta_t-\btheta^\star}^2$.  Lemma~\ref{lem:localized-squared-error}
therefore verifies
\eqref{eq:functional-occupation}.  Corollaries~\ref{cor:smooth-functional}
and~\ref{cor:coverage} now give the stated confidence interval.
\end{proof}

\begin{proof}[Proof of Lemma~\ref{lem:localized-squared-error}]
Put $\bm e_t:=\bx_t-\bx^\star$ and, for a sufficiently small compact
neighborhood $K_1$ of $\bx^\star$, define
$\tau_n:=\inf\{t\ge n:\bx_t\notin K_1\}$,
$A_{n,t}:=\{\tau_n>t\}$, and $E_n:=\{\tau_n=\infty\}$.
Almost-sure convergence gives $\Pp(E_n)\to1$.  It suffices to prove
\begin{equation}
\E[\norm{\bm e_t}^2\mathbf1_{E_n}]\le C_n\eta_t
\quad\text{for every fixed $n$ and all $t\ge n$}.
\label{eq:localized-mean-square}
\end{equation}
Indeed, for fixed $n$, Markov's inequality gives
\[
\Pp\!\left(\frac1{\sqrt T}\sum_{t=1}^T\norm{\bm e_t}^2>\varepsilon\right)
\le \Pp(E_n^c)+o(1)
 +\frac{2C_n}{\varepsilon\sqrt T}\sum_{t=n}^T\eta_t
=\Pp(E_n^c)+o(1).
\]
The first $o(1)$ comes from the finite prefix, while the last equality uses
$T^{-1/2}\sum_{t=n}^T\eta_t=O(T^{1/2-\kappa})$.  Letting $n\to\infty$
proves the lemma.

To prove~\eqref{eq:localized-mean-square}, let
$\bG=D\bg(\bx^\star)$ and choose $\bQ\succ0$ satisfying
$\bG^\top\bQ+\bQ\bG=\bI$.  The quantities in
\eqref{eq:auxiliary-iterate}--\eqref{eq:r-definition} and
\eqref{eq:mrc-decomposition}, write
\[
\bDelta_{t+1}=(\bI-\eta_t\bG)\bDelta_t+\eta_t\bepsilon_t
+\eta_t\widetilde\bzeta_t,\qquad
\widetilde\bzeta_t:=-\br_t-\bnu_t+\eta_t\bzeta_t.
\]
On $A_{n,t}$, local boundedness gives
$\norm{\bx_{t+1}-\bx_t}\le C\eta_t$.  Hence,
\eqref{eq:appendix-nu-bound} and
$|\eta_{t+1}/\eta_t-1|=o(\eta_t)$ imply
$\norm{\bnu_t}\le C\eta_t$, while \eqref{eq:r-definition} implies
$\norm{\br_t}\le C\norm{\bm e_t}^2+C\eta_t$.  Since
$\bm e_t=\bDelta_t+O(\eta_t)$, shrinking $K_1$ yields
$\norm{\widetilde\bzeta_t}\le\delta_0\norm{\bDelta_t}+C\eta_t$
for a sufficiently small $\delta_0>0$.

Set $\bm M_t:=\bI-\eta_t\bG$.  The Lyapunov identity gives
$\norm{\bm M_t\bm z}_{\bQ}^2
\le(1-a\eta_t)\norm{\bm z}_{\bQ}^2$ for some $a>0$.
Moreover, $\bDelta_t$ is $\cF_{t-1}$-measurable and
$\E(\bepsilon_t\mid\cF_{t-1})=\bzero$, so their cross term vanishes.
Expanding the square gives
\begin{align*}
\E(\|\bDelta_{t+1}\|_{\bQ}^2\mid\cF_{t-1})
={}&\|\bm M_t\bDelta_t\|_{\bQ}^2
 +\eta_t^2\E(\|\bepsilon_t\|_{\bQ}^2\mid\cF_{t-1})\\
&+2\eta_t\E\!\left[
\left\langle\bm M_t\bDelta_t+\eta_t\bepsilon_t,
\widetilde\bzeta_t\right\rangle_{\bQ}\middle|\cF_{t-1}\right]
 +\eta_t^2\E(\|\widetilde\bzeta_t\|_{\bQ}^2\mid\cF_{t-1})\\
&\le(1-a\eta_t)\|\bDelta_t\|_{\bQ}^2+C\eta_t^2
 +C\eta_t(\|\bDelta_t\|+\eta_t)\|\widetilde\bzeta_t\|
 +C\eta_t^2\|\widetilde\bzeta_t\|^2\\
&
\le(1-b\eta_t)\|\bDelta_t\|_{\bQ}^2+C\eta_t^2
\end{align*}
on $A_{n,t}$ for some $b>0$, after choosing $\delta_0$ small and $t$ large.
Thus, for
$a_t:=\E[\norm{\bDelta_t}_{\bQ}^2\mathbf1_{A_{n,t}}]$, we have
$a_{t+1}\le(1-b\eta_t)a_t+C\eta_t^2$.  Iteration and
Lemma~\ref{lem:recursion-matrices}(i),(iii) give
\[
a_{t+1}\le C e^{-b'\sum_{i=n}^t\eta_i}a_n
+C\sum_{j=n}^t\eta_j^2e^{-b'\sum_{i=j+1}^t\eta_i}
\le C_n\eta_{t+1}.
\]
Since
$\bm e_t=\bDelta_t+\eta_t\cP\bU(\bx_t,\xi_{t-1})$ and $E_n\subset A_{n,t}$,
$\E[\norm{\bm e_t}^2\mathbf1_{E_n}]\le C_n\eta_t$, as required.
\end{proof}

\FloatBarrier
\appendixsection{Additional Experimental Details}
\label{app:series-ablation}

\begin{figure}[!t]
	  \vspace{-10pt}
  \centering
  \includegraphics[width=\linewidth]{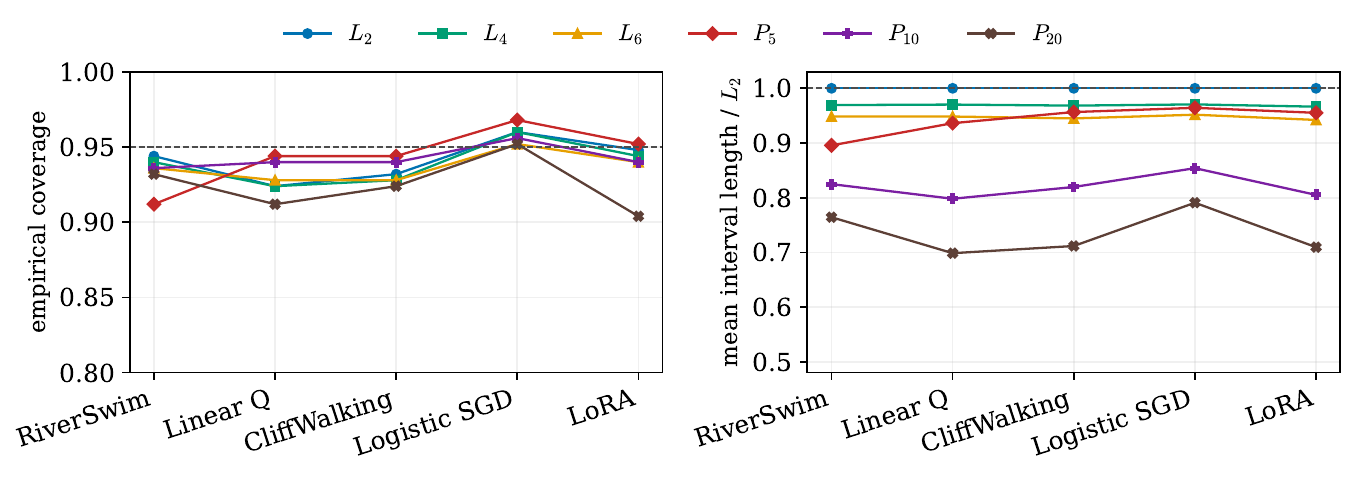}
  \caption{Functional ablation over 250 replications at 100,000 post-warm-up
  iterations.  The left panel reports empirical coverage, and the right panel
  reports mean interval length relative to $L_2$.}
  \label{fig:series-ablation}
  \vspace{-10pt}
\end{figure}

\paragraph{Polynomial-Dimension Ablation.}
Figure~\ref{fig:series-ablation} compares the polynomial-series functionals
$P_5$, $P_{10}$, and $P_{20}$ with the $L_2$, $L_4$, and $L_6$ bridge functionals, leaving
every application design unchanged.  At the final checkpoint, average coverage
across the five applications is $94.4\%$, $94.2\%$, and $92.5\%$ for
$P_5,P_{10},P_{20}$, while the corresponding mean length ratios relative to $L_2$ are
$94.2\%$, $82.0\%$, and $73.5\%$.  Thus additional projections reduce width, as
Proposition~\ref{prop:width-benchmark} suggests, but coverage becomes less
stable for larger $K$, most visibly for $P_{20}$.  We therefore use the
moderate fixed choice $P_5$ in the main experiments.

\phantomsection\label{app:reproducibility}
\paragraph{Simulation Protocol and Baselines.}
Section~\ref{sec:experiments} describes the interval methods and five
data-generating processes.  Here we record the remaining implementation
details.  The bootstrap's eight checkpoints use 500--10,000 observations,
equivalent to 5,500--110,000 update-function evaluations because $B=10$.
Warm-up is 5,000 observations except for CliffWalking, which uses 500,000;
these costs are excluded from the axes, and running sums are then restarted.
The perturbed paths start 1,000 steps earlier, adding 10,000 unplotted updates.
Table~\ref{tab:computation} reports the median of seven Python/NumPy runs on
Apple-arm64 using a pre-generated 50,000-transition tabular-Q path and online
bootstrap multipliers.

\paragraph{Additional Q-Learning Details.}
The linear-Q target uses the average feature vector over its 20 state--action
pairs.  The CliffWalking start state is state 36 in Gymnasium's indexing, and the
target is coordinate $(36,0)$.  Independent soft value iteration gives the value
$-42.614$.  The observed transition outcomes, including terminal flags and cliff
rewards, are sampled directly from Gymnasium's transition table rather than from an aggregated approximation.

\paragraph{LoRA Initialization.}
The initial factors are
$\bB_0=0.45\widetilde{\bu}$ and $\bA_0=0.45\widetilde{\bw}^{\top}$, where
$\widetilde{\bu}$ and $\widetilde{\bw}$ normalize
$\bu_\star+(0.25,0.15,-0.10)^\top$ and
$\bw_\star+(-0.15,0.20,0.10,0.05)^\top$.  Thus the experiment does not begin at the
target factor directions.  The exact scalar target is
$(\bW^\star)_{11}=0.6822$.

\bibliography{bib/stat}

\end{document}